\input amstex
\documentstyle{amsppt}
\NoRunningHeads
 
\magnification=1200
\pagewidth{32pc}
\pageheight{42pc}
\vcorrection{1.2pc}

\define\wh{\widehat}

\define\bc{\Bbb C}

\define\bbz{{\Bbb Z}^{\nu+1}}
\define\bbl{{\hat{\Cal L}_{{\Bbb C}_Q}}}

\define\vep{\varepsilon}

\topmatter
\title A Unified View of some Vertex Operator Constructions
\endtitle
\author Stephen Berman\footnotemark"$^{1}$", Yun
Gao\footnotemark"$^{2}$", Shaobin Tan
\footnotemark"$^{3}$"
\endauthor
\footnotetext"$^1$"{ Department of Mathematics and Statistics, University
of Saskatchewan,
Saskatoon, Canada
S7N 5E6; berman\@math.usask.ca}
\footnotetext "$^2$"{Department of Mathematics and Statistics,
 York University, Toronto, Canada M3J 1P3; ygao\@ yorku.ca}
\footnotetext"$^3$"{ Department of Mathematics, Xiamen University, Xiamen,
China 361005; tans\@jingxian.xmu.edu.cn} 

 
\abstract
We present a general vertex operator construction based on the Fock space 
for an affine Lie algebras of type $A$.
This construction allows us to give a unified treatment for both the
homogeneous and principle
realizations of the affine Lie algebras $\hat{gl}_N$ as well as  for some extended affine Lie algebras 
coordinatized by certain quantum tori.
 
\endabstract
 
\endtopmatter \document
 
\subhead\S 0. Introduction\endsubhead
 
\medskip
 This paper presents a unified view of certain vertex operator constructions for some of the extended affine Lie 
algebras 
(EALA's for short) which are coordinatized by certain quantum tori. Recall that for the affine Kac-Moody Lie algebras 
vertex operator representations were developed in [LW] and [KKLW] for the principal realizations and in [FK],[S] 
in the homogeneous realization. Our motivation comes from the paper [F1] of I. Frenkel, where he presented a unified 
construction for both the principal and homogeneous realizations of the affine Lie algebras of type $A^{(1)}.$ 
This is accomplished by using the affine algebra $\hat{gl}_M$ rather than $\hat {sl}_M.$ Moreover Frenkel used a 
Clifford  algebra structure, which was inherent in his situation, to define a new type of normal ordering which then 
led to his unified view. The Clifford structure had been studied before in the works [F3,4] and [KP].

   The structure theory of EALA's has been developed over the last ten years(see [H-KT], [BGK],[AABGP] and [ABGP]). 
Roughly speaking these Lie algebras are generalizations of both the affine Kac-Moody Lie algebras and the finite 
dimensional simple Lie algebras over the complex numbers which admit Laurent like coordinates  in a finite number 
of variables. It turns out that algebras of different types admit different types of coordinates. For example, 
those of type $A_l$ admit the non-commutative quantum torus as coordinates (see [M],[BGK]).
Representations for these Lie algebras over quantum tori have been constructed in
[JK], [G-KK], [G1,2,3], [BS], [VV].  When $l=1$ there are even Jordan 
algebras which serve as coordinates of EALAs.  The study of representations for this type of
Lie algebra has been initiated in [T1]. 
Perhaps the  examples which have attracted the most attention so far are the toroidal algebras 
which have the commutative associative Laurent polynomials as their coordinates. The toroidal 
algebras have been studied since the mid 80's both in terms of their structure theory as well as their representation 
theory (see for example [F2], [MRY], [Y], [W], [BC], [FJW],[T2]). For our purposes we 
want to mention that vertex operator representations have played a predominant role in much of this work. Indeed, in 
[G1,2], one finds both homogeneous and principal realizations given for many of the EALA's  of type $A$. The principal
realization for those EALAs was also implicitly given by [G-KK] in studying the so-called $\Gamma$-conformal algebras.  
Our goal 
in this work is to unify the various approaches and show how they all follow from the same type of approach. Of course, 
the work in the affine case, namely [F1],  shed light on doing this.

   Working with a standard type of Fock space we are able to define a general type of vertex operator which depends on 
a non-zero scalar from $\bc$ and to then  compute the commutator of two of these operators. This is presented in the 
second section of this paper while, in the first section, we give the basics on the Lie algebras, which are all of 
type $A$, which we will later go on to find representations for. Already in Section one it is evident that there is 
somewhat of a unified picture for these algebras. When we define our vertex operators in Section two the reader will 
see that we are using a Clifford algebra structure to define the normal ordering we are using, just as was done in [F1]. 
In the third section we introduce some Lie algebras associated to certain choices of subgroups, $G$, of non-zero 
complex numbers as well as the choice of a positive integer $M$. These Lie algebras are spanned by the moments of our 
vertex operators and hence, by  construction, we automatically have a representation for this Lie algebra. We show the 
representations we have are completely reducible and find the irreducible components. In the fourth and final section 
we show how certain choices of the group $G$ and the integer $M$ lead to representations of the algebras of section one. 
Thus, we recover both the principal and homogeneous representations for the affine algebras of type $A$ as well as those 
for the EALA's studied in [BS], [G1,2], [G-KK]. It is from seeing the various applications in this fourth section that one 
understands the unification of our treatment. Finally we want to emphasize that this unified treatment would not have 
been possible without first knowing the particular special cases of this result.

\bigskip
 
\subhead \S 1. Preliminaries\endsubhead
 
\medskip
 
In this section we shall review some of the  basics on Lie algebras coordinatized
by quantum tori. We present this from a general point of view which unifies our treatment. For notation we always denote the integers, positive
 and negative integers respectively by ${\Bbb Z}$,
${\Bbb Z}_+$, and ${\Bbb Z}_-$.
 
Let $\frak g$ be any associative $\Bbb C$-algebra with a symmetric
bilinear form $(\cdot,\cdot)$: $\frak g$$\times \frak g$$\to {\Bbb C}$
such that $(xy,z)=(x,yz)$ for $x,y,z\in$ $\frak g$. Let
$A=\oplus_{\alpha\in {\Bbb Z}^{\nu+1}} A_{\alpha}$, $\nu\ge 0$, be any
${\Bbb Z}^{\nu+1}$-graded associative algebra such that
dim$A_{\alpha}<\infty$ for all $\alpha\in {\Bbb Z}^{\nu+1}$. Fix a base
$(x_{i
\alpha})_{i\in I_{\alpha}}$ of $A_{\alpha}$, where $I_{\alpha}$ is the
index set corresponding to the subspace $A_{\alpha}$. Let 
 $d_0,d_1,\cdots,d_{
\nu}$ be degree derivations of $A$ such that $d_ix=\alpha_ix$ for $x\in
A_{\alpha}$, $i=0,1,\cdots,\nu$ and $\alpha=(\alpha_0,
\cdots, \alpha_{\nu})\in {\Bbb Z}^{\nu+1}.$ We define a $\Bbb C$-linear
map $\phi:$ $A\to {\Bbb C}$ by linear extension of
$$
  \phi(x_{i\alpha})=\cases 1 &\text{if}\;\;{\alpha=(0,\cdots, 0)}\\
0 &\text {otherwise}.\endcases\tag 1.1
$$
for $i\in I_{\alpha},$ $\alpha\in {\Bbb Z}^{\nu+1}$. The tensor product
${\frak g}\otimes_{\Bbb C}A$, with the canonical product $(x\otimes a)
(y\otimes b)=xy\otimes ab$ for $x,y\in {\frak g}, a,b\in A$, is also an
associative algebra. Moreover with the commutator product
$[x\otimes a, y\otimes b]_{loop}=(x\otimes a)
(y\otimes b)-
(y\otimes b)(x\otimes a),$ ${\frak g}\otimes_{\Bbb C}A$ forms a ${\Bbb
Z}^{\nu+1}$-graded Lie algebra. We call this algebra a loop type
Lie algebra. Consider the vector space
$$
{\hat {\frak g}}_A:=({\frak g}\otimes_{\Bbb C}A)\oplus {\Cal C}\tag 1.2
$$
where $\Cal C$=$\oplus_{0\le i\le \nu}{\Bbb C}c_i$ is a
$\nu+1$-dimensional
vector space. There is an alternating bilinear map $[\cdot,\cdot]
:$ ${\hat {\frak g}}_A\times {\hat {\frak g}}_A\to {\hat {\frak g}}_A$
determined by the conditions:
$$
[c_i,{\hat {\frak g}}_A]=0
$$
$$
[x\otimes a, y\otimes b]=[x\otimes a, y\otimes b]_{loop}+(x,y)\sum_{0\le
i\le\nu}\phi((d_ia)b)c_i
$$
for $x,y\in {\frak g}, a,b\in A$ and $i=0,1,\cdots, \nu$. It is
straightforward to check that ${\hat {\frak g}}_A$ is a Lie algebra.
Indeed
there is an exact sequence of Lie algebras with canonical maps
$$
0\to \oplus_{0\le i\le\nu}{\Bbb C}c_i\to {\hat {\frak g}}_A\to 
{\frak g}\otimes_{\Bbb C}A\to 0, 
$$
and so we have that ${\hat {\frak g}}_A$ is a central extension of the loop type Lie
algebra ${\frak g}\otimes_{\Bbb C}A$.
 
  Let $M_{\infty}(\bc)=$span$_{\bc}\{E_{ij}|\;\;1\le i,j <\infty\}$ be
the infinite matrix algebra, where $E_{ij}$ is the infinite matrix with a $1$
in the $(i,j)$-entry and zero's elsewhere. We also let 
$M_n(\bc)=$span$_{\bc}\{E_{ij}|\;\;1\le i,j\le n\}$. This subspace
of $M_{\infty}(\bc)$ for $n\ge 1$ is isomorphic to the usual matrix algebra of $n \times n$ 
matrices with entries in $\bc$ .

Let $Q=(q_{ij})$ be a $(\nu+1)\times(\nu+1)$ matrix with entries
$q_{ij}\in
{\Bbb C}^{\times}$ satisfying $q_{ii}=1$ and $q_{ij}=q_{ji}^{-1}$ for
$0\le i,j\le \nu$. The quantum torus associated with the matrix
$Q$ is a unital associative ${\Bbb C}$-algebra
${\Bbb C}_Q :={\Bbb C}_Q[t_0^{\pm 1},\cdots,t_{\nu}^{\pm 1}]$
with generators $t_0^{\pm 1},\cdots,t_{\nu}^{\pm 1}$ and relations
 $
t_it_i^{-1}=t_i^{-1}t_i=1, \quad  t_it_j=q_{ij}t_jt_i,
$
for $0\le i,j\le \nu$. If $Q$ is $2\times 2$ matrix, so then $\nu=1$, 
the matrix $Q=(q_{ij})$ is determined by a single $q=q_{10}$. In this case we often simply
denote
${\Bbb C}_Q={\Bbb C}_Q[t_0^{\pm 1},t_1^{\pm 1}]$ by ${\Bbb C}_{q}$.   
Choose the bilinear
form on $M_n(\Bbb C)$ to be the trace form. Set $A=\bc_Q[t^{\pm
1}_0,\cdots,
t^{\pm 1}_{\nu}]$,
 with the ${\Bbb Z}^{\nu+1}$-gradation $A=\oplus_{\alpha\in
\bbz}A_{\alpha}$, where
the subspace $A_{\alpha}$ is spanned by
$t^{\alpha}=t_0^{\alpha_0}t_1^{\alpha_1}\cdots t_{\nu}^{\alpha_{\nu}}$
for $\alpha=
(\alpha_0,\cdots, \alpha_{\nu})\in\bbz$.      Define $\sigma_Q:$
$\bbz\times\bbz\to \bc$ by
$$
\sigma_Q(\alpha,\beta)=\prod_{0\le i<j\le\nu}q_{ji}^{\alpha_j\beta_i}\tag
1.3
$$
for $\alpha=(\alpha_0,\cdots,\alpha_{\nu})$,
$\beta=(\beta_0,\cdots,\beta_{\nu})\in\bbz$. Then we have
 $
t^{\alpha}t^{\beta}=\sigma_Q(\alpha,\beta)t^{\alpha+\beta}.$

The proof of the following Lemma is clear.

\medskip

\proclaim
{\bf Lemma 1.4} Let $m,n\ge 1$ be integers. Then there is a Lie algebra
isomorphism
$$
(M_m(\bc)\otimes M_n(\bc))^{\wedge}_{\bc_Q}\cong
(M_{mn}(\bc))^{\wedge}_{\bc_Q}
$$
which is given by
$$
E_{ij}\otimes E_{kl}\otimes t^{\alpha}\mapsto E_{(i-1)n+k,(j-1)n+l}\otimes
t^{\alpha}
$$
$$
c_s\mapsto c_s,\;\;\;\; s=0,1,\cdots, \nu
$$
for $\alpha=(\alpha_0,\cdots,\alpha_{\nu})\in\bbz$, $1\le i,j\le m, 1\le
k,l\le n.$
\endproclaim
 
 Let $\bbl$ be the Lie subalgebra of $(M_m(\bc)\otimes
M_n(\bc))^{\wedge}_{\bc_Q} $ generated by elements of the form
$ E_{ij}\otimes E_{kl}\otimes t_0^{\alpha_0(n-1)+l-k}t^{\alpha}$ for $1\le
i,j\le m$, $1\le k,l\le n$ and $\alpha=
(\alpha_0,\cdots,\alpha_{\nu})\in \bbz$. The following result gives the structure of $\bbl .$

\proclaim 
{\bf Proposition 1.5} 
$$
 \bbl\cong(M_m(\bc)\otimes
M_n(\bc))^{\wedge}_{\bc_{Q^{\ast}}}  
$$
 where ${\bc_Q}={\bc_Q}[t_0^{\pm
1},\cdots, t_{\nu}^{\pm 1}]
$ with $Q=(q_{ij})$, and ${\bc_{Q^{\ast}}}={\bc_{Q^{\ast}}}[\tau_0^{\pm
1},\cdots,
\tau_{\nu}^{\pm
1}]
$ with $Q^{\ast}=(q^{\ast}_{ij})$ such that $q^{\ast}_{ij}=q_{ij}$ if
$i,j\not= 0$, and $q^{\ast}_{ij}=q^n_{ij}$ if $i=0$ or $j=0$.
 
\endproclaim
Proof. Define a linear map $f:$ $(M_m(\bc)\otimes
M_n(\bc))^{\wedge}_{\bc_{Q^{\ast}}}  \to \bbl$ by
$$
E_{ij}\otimes E_{kl}\otimes \tau^{\alpha}\mapsto (\prod_{1\le
s\le\nu}q_{s0}^{l\alpha_s})E_{ij}\otimes
E_{kl}\otimes
t_0^{(n-1)\alpha_0+l-k}t^{\alpha}-k\delta_{ij}\delta_{kl}\delta_{\alpha,{
0}}c_0
 $$
$$
c_0\mapsto nc_0,\;\;\;\;\;\; c_s\mapsto c_s, \;\;\;s=1,2,\cdots, \nu
$$
for $\alpha=(\alpha_0,\cdots, \alpha_{\nu})\in\bbz$, $1\le i,j \le m$ and
$1\le k,l\le n$. Let ${\bar
{\alpha}}=(n\alpha_0+l-k,\alpha_1,\cdots,\alpha_{\nu})$, and 
${\bar {\alpha}'}=(n\alpha_0'+l'-k',\alpha_1',\cdots,\alpha_{\nu}')\in
\bbz$. Using the identity
$$
\sigma_{Q}({\bar {\alpha}},{\bar
{\alpha}'})=\sigma_{Q^{\ast}}(\alpha,\alpha')\prod_{1\le j\le
\nu}q_{j0}^{\alpha_j(l'-k')}\tag 1.6
$$
one can easily check that the map $f$ is the desired Lie algebra
isomorphism.  {\qed}

Putting together the two previous results we get the following identification of $\bbl$.

\medskip
\proclaim
 {\bf Corollary 1.7} 
$$
\bbl\cong(M_{mn}(\bc))^{\wedge}_{\bc_{Q^{\ast}}}  
$$
where $Q$ and $Q^{\ast}$ are given in the previous proposition .
{\qed}
\endproclaim
\medskip
 
Let $\xi=\xi_n$ be an $n$-th primitive root of unity and let $E,F\in
M_n(\bc)$ be defined by saying
$$
E=E_{12}+E_{23}+\cdots +E_{n-1,n}+E_{n1}, \quad
F=\sum_{i=1}^nE_{ii}(\xi^{i-1}).
\tag 1.8
$$
Then the following fact is well-known.
 
\proclaim{Lemma 1.9} The set of matrices $\{F^iE^j\}_{1 \leq i, j\leq
n}$, forms a basis of
the 
matrix algebra $M_n(\bc)$ (so a basis of the general linear Lie algebra
${gl}_n(\Bbb{C})$). Moreover,
$$EF=\xi FE, \quad E^n=F^n=Id_n, \tag 1.10$$
and
$$
E_{ij}={\frac 
1 n}\sum^{n-1}_{k=0}\xi^{k(1-i)}F^kE^{j-i},\;\;\;\;F^iE^j=\sum^n_{l=1}
\xi^{i(
l-1)}E_{l,{\overline {l+j}}}\tag 1.11
$$
for $1\le i,j\le n$, where, for notation, we are letting ${\bar l}$ denoted the unique integer, 
$l$, in $\{ 1, 2, \cdots, n \}$ such that ${\bar l}=l($mod $n)$.
\endproclaim
 
\medskip

Note that
$$
E_{ij}\otimes E_{kl}\otimes t_0^{\alpha_0(n-1)+l-k}t^{\alpha}=
\sum_{s=0}^{n-1}\xi^{s(1-k)}E_{ij}\otimes F^sE^{l-k}\otimes
t_0^{\alpha_0n+l-k}t_1^{\alpha_1}\cdots t_{\nu}^{\alpha_{\nu}}
$$
$$
=\sum_{s=0}^{n-1}\xi^{s(1-k)}E_{ij}\otimes F^sE^{\alpha_0'}\otimes
t_0^{\alpha_0'}t_1^{\alpha_1}\cdots t_{\nu}^{\alpha_{\nu}}
$$
where $\alpha_0'=\alpha_0n+l-k$, for $1\le i,j\le m$, $1\le k,l\le n$, and
$\alpha=(\alpha_0,\alpha_1,\cdots, \alpha_{\nu})\in
\bbz.$ From this one sees that the Lie subalgebra $\bbl$ of
$(M_m(\bc)\otimes M_n(\bc))^{\wedge}_{\bc_Q}  $ has a basis
of the form
$$
E_{ij}\otimes F^kE^{l_0}\otimes t_0^{l_0}\cdots
t_{\nu}^{l_{\nu}},\;\;\;\; c_0,\; c_1,\cdots, c_{\nu}\tag 1.12
$$
where $1\le i,j\le m$, $0\le k\le n-1$ and $l_0,l_1,\cdots,l_{\nu}\in
{\Bbb Z}$. Moreover, the commutation relation of $\bbl$ are
determined by
$$
[E_{ij}\otimes F^kE^{\alpha_0}\otimes t^{\alpha},
E_{i'j'}\otimes F^{k'}E^{\alpha_0'}\otimes t^{\alpha'}]\tag 1.13
$$
$$
=\delta_{ji'}\xi^{\alpha_0k'}\sigma_{Q}(\alpha,\alpha')E_{ij'}\otimes 
F^{k+k'}E^{
\alpha_0+\alpha_0'}\otimes t^{\alpha+\alpha'}
$$
$$
-\delta_{j'i}\xi^{\alpha_0'k}\sigma_{Q}(\alpha',\alpha)E_{i'j}\otimes 
F^{k+k'}E^{
\alpha_0+\alpha_0'}\otimes t^{\alpha+\alpha'}
+n\delta_{ji'}\delta_{ij'}\delta_{{\overline {k+k'}},{
0}}\delta_{\alpha+\alpha', 0}\xi^{\alpha_0k'}\sum_{0\le
s\le\nu}\alpha_sc_s
$$
for $1\le i,i',,j,,j'\le m,$ $0\le k,k'\le n-1$,
$\alpha=(\alpha_0,\cdots,\alpha_{\nu}),$
$\alpha'=(\alpha_0',\cdots,\alpha_{\nu}')
\in \bbz,$ as well as the fact that the elements $c_0,\cdots, c_{\nu}$ are central in $\bbl$.
 
 From now on we will identity the Lie algebra $(M_m(\bc)\otimes
M_n(\bc))^{\wedge}_{\bc_Q}$ with $(M_{mn}(\bc))^{\wedge}_{\bc_Q}$,
and also identify $\bbl$ with   $(M_{mn}(\bc))^{\wedge}_{\bc_{Q^{\ast}}}$, where
$Q=(q_{ij}),$ $Q^{\ast}=(q^{\ast}_{ij})$ and as above $q^{\ast}_{ij}=q_{ij}
$ if $i,j\not= 0$, and $q^{\ast}_{ij}=q^n_{ij}$ if $i=0$ or $j=0$. For 
simplicity we will write $a\alpha=(a\alpha_1,\cdots, a\alpha_{\nu})$ for 
$a\in{\Bbb Z}$ and $\alpha\in {\Bbb Z}^{\nu}$, also we will write
$$
q_0^{\alpha}=q_{10}^{\alpha_1}\cdots q_{\nu 0}^{\alpha_{\nu}}
$$
for $q_0=(q_{10},\cdots, q_{\nu 0})\in {\Bbb C}^{\nu}$.

The
Lie algebra structure (1.13) of $\bbl$ can be described by formal
power series identities. For this purpose we let $z,z_1,z_2$ be formal
variables. For $1\le i,j\le m, $  $0\le k\le n-1$ and
$\alpha=(\alpha_1, \cdots, \alpha_{\nu})\in {\Bbb Z}^{\nu}$, we set
$$
X_{ij}^k(\alpha,z)=\sum_{l\in \Bbb Z}(E_{ij}\otimes F^k E^l\otimes
t_0^lt_1^{\alpha_1}\cdots t_{\nu}^{\alpha_{\nu}})z^{-l}\in {\bbl}[[
z,z^{-1}]],\tag 1.14
$$
and $\delta(z)=\sum_{l\in {\Bbb Z}}z^l$, $(D\delta)(z)=\sum_{l\in \Bbb
Z}lz^l$. Then the algebra structure of $\bbl$ is described by the
following lemma.
 
\proclaim
{\bf Lemma 1.15} Let $1\le i,j,i',j'\le m$, $0\le k,k'\le n-1$,
$\alpha=(\alpha_1,\cdots, \alpha_{\nu})$, $\alpha'
=(\alpha_1',\cdots, \alpha_{\nu}')\in {\Bbb Z}^{\nu}$. Then the following
power series identity  is equivalent to (1.13) 
$$
[X_{ij}^k(\alpha, z_1), X_{i'j'}^{k'}(\alpha',
z_2)]=\delta_{ji'}\sigma(\alpha,\alpha')X_{ij'}^{\overline
{k+k'}}(\alpha+\alpha',\xi^{-k'}z_1)\delta({\frac
{\xi^{k'}z_2} {z_1q^{\alpha}_0}})\tag 1.16
$$
$$
-\delta_{j'i}\sigma(\alpha',\alpha)X_{i'j}^{\overline
{k+k'}}(\alpha+\alpha',\xi^{-k}z_1)\delta({\frac
{z_2q_0^{\alpha'}}{\xi^{k}z_1} })
$$
$$
+n\delta_{ji'}\delta_{ij'}\delta_{\overline 
{k+k'},0}\delta_{\alpha+\alpha',
{ 0}}\sigma(\alpha,\alpha')\{(D\delta)({\frac
{\xi^{k'}z_2} {z_1q_0^{\alpha}}})c_0+\delta({\frac
{\xi^{k'}z_2} {z_1q_0^{\alpha}}})\sum_{1\le s\le
\nu}\alpha_sc_s\}
$$
where ${\bar k}=k($mod$ \ n)$ and $ k \in\{0,1,\cdots, n-1\}$.
 \endproclaim

As  very special cases, one  chooses $n=1,$ $\nu=0$, then $\bbl$ is just
the affine algebra ${\widehat {gl}_m}(\bc)$ in the so-called homogeneous
picture; while if one chooses $m=1,$ $\nu=0$ then $\bbl$ is 
the affine algebra ${\widehat {gl}_n}(\bc)$ in the so-called
principal picture. In these two cases, the identity (1.16) can simply be written
as
follows 
$$
[X_{ij}^0(z_1), X_{kl}^0(z_2)]=X_{il}^0(z_1)\delta_{jk}\delta({\frac
{z_2} {z_1}})-
X_{kj}^0(z_2)\delta_{il}\delta({\frac {z_2}
{z_1}})+\delta_{jk}\delta_{il}(D\delta)({\frac {z_2} {z_1}})c_0,
\tag 1.17
$$
for $1\le i,j,k,l\le m$, for the first of these  and
$$
[X^i_{11}(z_1), X^j_{11}(z_2)]=X^{\overline {i+j}}_{11}(z_2)\delta({\frac
{\xi^jz_2} {z_1}})-
X^{\overline {i+j}}_{11}(z_1)\delta({\frac {\xi^iz_1}
{z_2}})+n\delta_{{\overline{i+j}},0}(D\delta)({\frac {\xi^jz_2} {z_1}}
)c_0.\tag 1.18
$$
for $0\le i,j\le n-1$,  ${\bar i}=i($mod$ \ n)$ and $ i \in \{0,1,\cdots, n-1\}$ for the second one.
 
  Moreover if we choose $n=1, \nu=1$, or $m=1, \nu=1$, and write
$q_{10}=q$, then $\bbl$ gives respectively the homogeneous
realization of the Lie algebra ${\widehat {gl}_m}(\bc_q)$, and the
principal
realization of ${\widehat {gl}_n}(\bc_{q^{n}})$. The
algebra structure of these two cases can
be described as follows
$$
  [X_{ij}^0(r,z_1), X_{kl}^0(s,z_2)]  
=  X_{il}^0(r+s,z_1)\delta_{jk}\delta({\frac
{z_2}{q^rz_1}})-
 X_{kj}^0(r+s,z_2)\delta_{il}\delta({\frac {q^sz_2}{z_1}}) \tag 1.19
$$
$$
 +\delta_{il}\delta_{jk}\delta_{r+s,0}((D\delta)({\frac
{z_2}{q^rz_1}})c_0+
 r\delta({\frac {z_2}{q^rz_1}})c_1) $$
 for $1\le i,j,k,l\le m$ and $r,s\in {\Bbb Z}$, for the first and
 $$
  [X^i_{11}(r,z_1), X^j_{11}(s,z_2)]= X^{\overline
{i+j}}_{11}(r+s,\xi^{-j}z_1)\delta({\frac
{\xi^jz_2}{q^rz_1}})-
 X^{\overline {i+j}}_{11}(r+s,\xi^{-i}z_2)\delta({\frac
{q^sz_2}{\xi^iz_1}})
\tag 1.20
 $$
$$
 +n\delta_{r+s,0}\delta_{\overline{i+j},0}((D\delta)({\frac
{\xi^jz_2}{q^rz_1}})c_0+r
 \delta({\frac {\xi^jz_2}{q^rz_1}})c_1)
 $$
 for $1\le i,j\le n$, $r,s\in {\Bbb Z}$, and ${\overline{i+j}}=i+j$(mod
$n$) for the second.
 
Finally if we choose $\nu=1$,  $m,n\ge 1$, write $q_{10}=q$, then
$\bbl$ is isomorphic to the affine Lie algebra
${\widehat {gl}_{mn}}(\bc_{q^n})$, which contains the special cases
mentioned above. The algebra structure is as follows.
$$
[X^k_{ij}(r,z_1),X^{k'}_{i'j'}(r',z_2)]=\delta_{ji'}X_{ij'}^{\overline{k+k'}}
(r+r',\xi^{-k'}z_1)\delta({\frac {\xi^{k'}z_2} {q^rz_1}})\tag 1.21
$$
$$
-\delta_{j'i}X_{i'j}^{\overline{k+k'}}
(r+r',\xi^{-k}z_2)\delta({\frac {q^{r'}z_2} {\xi^{k}z_1}})
$$
$$
+n\delta_{ji'}\delta_{ij'}\delta_{{\overline {k+k'}},0}\delta_{r+r',0}\{
(D\delta)({\frac {\xi^{k'}z_2} {q^rz_1}})c_0+r\delta({\frac {\xi^{k'}z_2}
{q^rz_1}})c_1\}.
$$
In subsequent sections we are going to give
 irreducible representations for a class of Lie algebras which include the
Lie algebras mentioned  above.
 
\medskip
\medskip
 
\subhead \S 2.  Fock Space and Vertex Operators \endsubhead
 
\medskip
 
In this section, we shall define the Fock space we need and  construct a family
of vertex
operators acting on it. Then we go on to derive the commutation relations between these vertex operators in various
situations.
Some of these commutation relations were implicitly worked out in [G1].
However, we will use the
ideas from [F1] to tie a Clifford algebra structure to our vertex operators. This makes our approach
very natural and concise.
 
Let $\vep_1,\dots,\vep_M$ $(M\ge 1)$ be symbols. We form lattices
$$
\Gamma_M=\oplus^M_{i=1}{\Bbb Z}\vep_i,\;\;\;\;Q_M=\oplus^{M-1}_{i=1}{\Bbb
Z}(\vep_i-\vep_{i+1}),
\tag 2.1 $$
with a symmetric bilinear form $(\vep_i,\vep_j)=\delta_{ij}$. We also
extend this bilinear form to the ${\Bbb
C}$-vector space
$$ H_M:={\Bbb C}\otimes \Gamma_M.\tag 2.2 $$
 
 For each $k \in \Bbb Z$ we take a copy of  $H_M$ with basis labeled by  $\vep_i(k)$  for $1\leq i \leq M,
k\in \Bbb Z$. That is, $\vep_i(k)$, is to be a  copy of $\vep_i$. We form a Lie algebra
$$
{\Cal H}_M={\hbox {span}}_{\Bbb C}\{\vep_i(k),\; c|1\le i\le
M,\;k\in{\Bbb Z}\},
\tag 2.3 $$
with the Lie product
$$[\alpha(k), \beta(l)]=k(\alpha,\beta)\delta_{k+l,0}c, \tag 2.4 $$
for $\alpha,\beta\in H_M,$
$k,l
\in {\Bbb Z}$, where $c$ is a central element. Let
$${\Cal H}^{\pm}_M={\hbox
{span}}\{\vep_i(k)|k\in{\Bbb Z}_{\pm},
1\le i\le M\}. \tag 2.5$$
 Then $${\hat {\Cal H}}_M={\Cal H}^{+}_M+{\Bbb C}c+{\Cal H}^{-}_M \tag
2.6$$
forms a Heisenberg
subalgebra of ${\Cal H}_M$.
 
Let ${\Cal S}({\Cal H}^{-}_M)$ be the symmetric algebra over the abelian
algebra ${\Cal H}^-_M$ and lat
$${\Bbb C}[\Gamma_M]:= \oplus_{\alpha\in \Gamma_M}{\Bbb
C}e^{\alpha}, \tag 2.7$$
 be the group algebra over $\Gamma_M$ twisted by a 2-cocycle so that 
$e^{\alpha}e^{\beta}=\epsilon(\alpha,
 \beta)e^{\alpha+\beta}$ for $\alpha,\beta\in\Gamma_M$. The cocycle 
$$\epsilon:\; \Gamma_M\times \Gamma_M\to \{\pm 1\}, \tag 2.8$$ is defined
by setting
$$\epsilon(\vep_i,\vep_j)=1  \text{ if }
i\le j,  \quad \epsilon(\vep_i,\vep_j)=-1 \text{ if }
i>j, \tag 2.9$$
 and $$\epsilon(\sum_im_i\vep_i,
\sum_jn_j\vep_j)=\prod_{i,j}(\epsilon(\vep_i,\vep_j))^{m_in_j}, \tag
2.10$$
for $m_i,n_j\in{\Bbb Z}$. One can easily check the following result.
 
\proclaim{Lemma 2.11} $\epsilon$ is bi-multiplicative on $\Gamma_M$.
Moreover.
$$\epsilon(\alpha,\alpha)=(-1)^{\frac{(\alpha,\alpha)}{2}},\;\;
\epsilon(\alpha,\beta)\epsilon(\beta,\alpha)
=(-1)^{(\alpha,\beta)}$$ for $\alpha,\beta\in Q_M$.
\endproclaim

Now we define the Fock space
$$
V_M={\Cal S}({\Cal H}^-_{M})\otimes {\Bbb C}[\Gamma_M]
\tag 2.12$$
which affords representations for both the Lie algebra ${\Cal H}_M$ and
the group algebra ${\Bbb C}[\Gamma_M
]$ with the following actions:
 
\medskip
 
$\;\;\;\;\vep_i(k).u\otimes e^{\beta}=k({\frac
{\partial}{\partial\vep_i(-k)}}u)\otimes e^{\beta},$ for $k\in{\Bbb
Z}_+$,
 
$\;\;\;\;\vep_i(k).u\otimes e^{\beta}=(\vep_i(k)u)\otimes e^{\beta},$ for
$k\in{\Bbb
Z}_-$,
 
$\;\;\;\;\vep_i(0).u\otimes e^{\beta}=(\vep_i,\beta)u\otimes e^{\beta},$
 
$\;\;\;\;c.u\otimes e^{\beta}=u\otimes e^{\beta},$ and
$e^{\alpha}.u\otimes e^{\beta}=\epsilon(\alpha,\beta)
u\otimes e^{\alpha+\beta},$
 
\medskip
 
\noindent for $\alpha,\beta\in\Gamma_M,$ $1\le i\le M$, and $u\in {\Cal
S}({\Cal H}^-_{M})$.
 For $\alpha\in\Gamma_M$, we define (we are using the standard notation from [FLM])
$$
\alpha(z)=\sum_{k\in{\Bbb Z}}\alpha(k)z^{-k}\in ({\hbox
{End}}V_M)[[z,z^{-1}]]
\tag 2.13$$
and
$$
E^{\pm}(\alpha,z)={\hbox {exp}}(\sum_{k\in{\Bbb Z}_{\pm}}{\frac
{\alpha(k)} {k}}z^{-k})\in ({\hbox {
End}}V_M)[[z,z^{-1}]].
\tag 2.14$$
Then the follow lemma is straightforward.
 
\proclaim{Lemma 2.15} For $\alpha,\beta\in\Gamma_M$, $a,b\in {\Bbb
C}^{\times}:={\Bbb
C}\setminus\{0\}$, we have
 
\medskip
 
$\;\;\;\;E^{\pm}(0,az)=1,$ $[\alpha(k), E^{+}(\beta,az)]=0$, if $k\ge 0$,
  
$\;\;\;\;[\alpha(z_1), E^{\pm}(\beta,
az_2)]=-(\alpha,\beta)E^{\pm}(\beta, az_2)\sum_{k\in{\Bbb Z}_{\mp}}({\frac {az_2}{z_1}})^k$,
 
$\;\;\;\;E^{\pm}(\alpha, az)E^{\pm}(\beta, bz)={\hbox
{exp}}(\sum_{k\in{\Bbb Z}_{\pm}}{\frac 1 k}(a^{-k}
\alpha(k)+b^{-k}\beta(k))z^{-k}),$
 
$\;\;\;\;E^{+}(\alpha, az_1)E^{-}(\beta, bz_2)=E^{-}(\beta,
bz_2)E^{+}(\alpha, az_1)(1-{\frac {bz_2}{az_1}})
^{(\alpha,\beta)}.$
\endproclaim

Let $v=\alpha_1(-1)^{k_1}\cdots\alpha_r(-r)^{k_r}\otimes e^{\beta}\in
V_M,$ we define a degree operator
$d_0$ of $V_M$ by setting
$$d_0v=(-\sum^r_{i=1}ik_i-{\frac 1 2}(\beta,\beta))v.\tag 2.16$$
 If $a$ is any non-zero
complex number we define operators
$$
z^{\alpha}.u\otimes e^{\beta}=z^{(\alpha,
\beta)}u\otimes
e^{\beta},\;\;\;\;\;
a^{\alpha}.u\otimes
e^{\beta}=a^{(\alpha, \beta)}u\otimes e^{\beta}\tag
2.17$$
for $\alpha,\beta\in \Gamma_M$, $u\in {\Cal S}({\Cal H}^-_M).$ Then $a^{\alpha}$ is just
the evaluation map, at $a$, 
of
the operator $z^\alpha$.
The following result is well-known.
 
\proclaim{Lemma 2.18}
$$\align &
[d_0, E^{\pm}(\alpha, az)]=-D_zE^{\pm}(\alpha, az)=(\sum_{k\in {\Bbb
Z}_{\pm}}\alpha(k)(
az)^{-k})E^{\pm}(\alpha, az), \\
& [\alpha(0), z^{\beta}]=0,\;\;\;\;\;\; z^{\alpha}
e^{\beta}=z^{(\alpha,\beta)}e^{\beta}
z^{\alpha}
\endalign $$
for $\alpha, \beta\in\Gamma_M$, $a\in {\Bbb C}^{\times},$ where
$D_z={\frac d {dz}}$.
 
\endproclaim
 
\medskip
 
We will have need to raise some of the  complex numbers which arise in our construction below to various powers and care must be taken with this. We thus set up the notation we use for this now.  For any complex number $a\not= 0$, there is a unique real number 
 $\theta\in [0,2\pi)$ such that $a=|a|e^{\theta\sqrt {-1}}$. We
define 
$$
Lna=\theta\sqrt{-1}+\ln |a|
$$
Viewing ${\Bbb
C}^{\times}=\Bbb{C}\setminus\{0\}$ as multiplicative group, we
call a subgroup $G$ of $
{\Bbb C}^{\times}$ {\bf  admissible}  if $G=T\times F$, where 
$T=<\xi>$ is a cyclic group of finite order $|T|$ and
$F=<q_j|j\in J>$ is a free abelian group with free generators $q_j, j \in J$. 
For $a=\xi^{-n_0}q^{n_1}_{i_1}\cdots q^{n_k}_{i_k}\in G$, where 
$n_0,n_1,\cdots, n_k\in {\Bbb Z}, i_1,\cdots, i_k\in J$, $0\le
n_0\le |T|-1$. We define
$$
a^r=e^{r(-n_0Ln\xi+n_1Lnq_{i_1}+\cdots +n_kLnq_{i_k})}\tag 2.19
 $$
for $r\in \Bbb C$.

   Recall  the
definition of limit of formal power series from [FLM]. Let $V$ be a
vector space over $\Bbb C$. Let
$$f(z_1,z_2)=\sum_{i,j\in{\Bbb Z}}a_{ij}z^i_1z^j_2\in V[[z^{\pm
1}_1,z^{\pm 1}_2]], $$
we say the limit, $\lim_{
z_2\to z_1}f(z_1,z_2)$, exists if, for any $l\in {\Bbb Z}$,
$a_{i,l-i}=0$
whenever $|i|>>0$, and
 write $$
\lim_{z_2\to z_1}f(z_1,z_2)=f(z_1,z_1)=\sum_{l\in {\Bbb Z}}(\sum_{i\in
{\Bbb
Z}}a_{i,l-i})z_1^l.\tag 2.20$$
For our purposes we
 need another notion of limit as well. Let $f(x,z)=\sum_{i\in {\Bbb
Z}}C_i(x)z^i$, where
$C_i(x)=\sum_jc_{ij}(x)v_j\in
{\Bbb C}(x)\otimes V$, and $j$ runs over a finite set for each 
fixed $i\in
{\Bbb Z}$, and
$c_{ij}(x)\in \Bbb{C}(x)$ are complex rational  functions. We say
the limit,
$\lim_{x\to
a}f(x,z)$, exists if the function $c_{ij}(x)$ has a usual limit
at the point
$a\in {\Bbb C}$ for all $i,j\in {\Bbb
Z}$, and write $$\lim_{x\to
a}f(x,z)=f(a,z)=\sum_i(\sum_jc_{ij}(a)v_j)z^j. \tag 2.21$$
 
\proclaim{Lemma 2.22}  For $1\le i\le M$, we have
$$\align & \lim_{a\to 1}{\frac 1 {1-a}}(a^{-\vep_i}-1)=\vep_i(0), \tag
2.23\\
& \lim_{a\to 1}{\frac 1 {1-a}}(\sum_{k\in {\Bbb Z}_{\pm}}{\frac
{\vep_i(k)} k}(az)^{-k}-
\sum_{k\in {\Bbb Z}_{\pm}}{\frac {\vep_i(k)} k}z^{-k})=\sum_{k\in {\Bbb
Z}_{\pm}}\vep_i(k)z^{-k}.
\tag 2.24\endalign$$
\endproclaim
\demo{Proof} For any $v=u\otimes e^{\beta}\in V_M$, let
$m=-(\vep_i,\beta)\in {\Bbb Z}$. Then
$$
{\frac 1 {1-a}}(a^{-\vep_i}-1).v={\frac {a^m-1}{1-a}}v,
$$
 and
$$
\lim_{a\to 1}{\frac 1 {1-a}}(a^{-\vep_i}-1).v=-mv=\vep_i(0).v
$$
as required. The second identity is clear. \qed \enddemo
 
\proclaim{Corollary 2.25}
$$
\lim_{a\to 1}{\frac 1 {1-a}}(E^{\pm}(\vep_i, az)-E^{\pm}(\vep_i,
z))=E^{\pm}(\vep_i, z)
\sum_{k\in {\Bbb Z}_{\pm}}\vep_i(k)z^{-k}.
$$
\endproclaim
\demo{Proof} Note that
$$
E^{\pm}(\vep_i, az)-E^{\pm}(\vep_i, z)
=\sum^{\infty}_{l=1}{\frac 1 {l!}}[(\sum_{k\in {\Bbb Z}_{\pm}}
{\frac {\vep_i(k)} {k}}(az)^{-k})^l-(\sum_{k\in {\Bbb Z}_{\pm}}
{\frac {\vep_i(k)} {k}}z^{-k})^l],
$$
and $A^l-B^l=(A-B)\sum^{l-1}_{j=0}A^{l-1-j}B^j$, we obtain by applying
the previous lemma
$$
\lim_{a\to 1}{\frac 1 {1-a}}(E^{\pm}(\vep_i, az)-E^{\pm}(\vep_i, z))
$$
$$
=\sum_{l=1}^{\infty}{\frac 1 {l!}}[l(\sum_{k\in {\Bbb Z}_{\pm}}
{\frac {\vep_i(k)} {k}}z^{-k})^{l-1}]\sum_{k\in {\Bbb
Z}_{\pm}}\vep_i(k)z^{-k}
=E^{\pm}(\vep_i, z)\sum_{k\in {\Bbb Z}_{\pm}}\vep_i(k)z^{-k}.\qed
$$
\enddemo
 
\proclaim{Corollary 2.26}
$$
\lim_{a\to 1}{\frac 1 {1-a}}(E^{\pm}(-\vep_i, z)E^{\pm}(\vep_i, az)-1)=
\sum_{k\in {\Bbb Z}_{\pm}}\vep_i(k)z^{-k}.
$$
\endproclaim
\demo{Proof} This follows from the fact that
$$
E^{\pm}(-\vep_i, z)E^{\pm}(\vep_i, az)-1=E^{\pm}(-\vep_i,
z)(E^{\pm}(\vep_i, az)-E^{\pm}(\vep_i, z)),
$$
and the previous corollary.\qed
\enddemo
 
\medskip
 
For $\alpha\in\Gamma_M$, we define
$$
X(\alpha,z)=E^-(-\alpha,
z)E^+(-\alpha,z)e^{\alpha}z^{\alpha}z^{(\alpha,\alpha)/2}.
\tag 2.27$$
 
We may formally write
$$X(\alpha,z)=\sum_{k\in {\Bbb
Z}+(\alpha,\alpha)/2}x_k(\alpha)z^{-k},\tag 2.28$$
 where $x_k(
\alpha)\in {\hbox {End}}(V_M)$ for $k\in {\Bbb Z}+(\alpha,\alpha)/2$. It
is known from [F1]
that
, if
$(\alpha,\alpha)=1$, the operators $\{x_k(\alpha),x_k(-\alpha)|k\in {\Bbb
Z}+{\frac 1 2}\}$ generate a
Clifford algebra with the relations
 
$$
\{x_k(\alpha),x_{-l}(-\alpha)\}=\delta_{kl},\;\{x_k(\alpha),
x_{l}(\alpha)\}=0,\;
\{x_k(-\alpha),x_{l}(-\alpha)\}=0
\tag 2.29
$$
 
\noindent for all $k,l\in {\Bbb Z}+{\frac 1 2}$. Related to this Clifford
structure,
we define the following  normal ordering (see [F1] and [G3]):
 
$$
:x_k(\vep_i)x_{-l}(-\vep_j):=x_k(\vep_i)x_{-l}(-\vep_j)-\delta_{ij}\delta_{kl}
\theta(k) \tag
2.30
$$
 
\noindent for $k,l\in {\Bbb Z}+{\frac 1 2},$ $1\le i,j\le M$, where
$\theta(k)=0$ if $k<0$, $\theta(k)=1$ if
$k>0$. By applying (2.15) and (2.18), one can easily prove the following
result
 
\proclaim{Lemma 2.31} For $1\le i, j\le M$, and $a\in{\Bbb C}^{\times}$,
we have
$$
 :X(\vep_i,z_1)X(-\vep_j,az_2):
 $$
 $$
= (1-{\frac {az_2}{z_1}})^{-\delta_{ij}}{\frac
{(az_2)^{\delta_{ij}/2}}
{z_1^{\delta_{ij}/2}}} \big(\epsilon(\vep_i,\vep_j)z_1^{\frac
{(\vep_i,\vep_i-\vep_j)}{2}}\cdot (az_2)^{-\frac
{(\vep_j,\vep_i-\vep_j)}{2}}
$$
$$
 \cdot e^{\vep_i-\vep_j}z_1^{\vep_i}(az_2)^{-\vep_j}
E^-(-\vep_i,z_1)E^-(\vep_j,az_2)E^+(-\vep_i,z_1)E^+
(\vep_j,az_2)-\delta_{ij}
). 
$$ 
In particular, if $i\not= j$, then
$$ \align & :X(\vep_i,z_1)X(-\vep_j,az_2):=\epsilon(\vep_i,\vep_j)z_1^{\frac 1
2}(az_2)^{\frac 1 2}
e^{\vep_i-\vep_j}z_1^{\vep_i}(az_2)^{-\vep_j}\tag 2.32 \\
& \;\;\;\;\;\;\;\;\;\;\cdot
E^-(-\vep_i,z_1)E^-(\vep_j,az_2)E^+(-\vep_i,z_1)E^+(\vep_j,az_2),
\endalign $$
and, if $i=j$, then
$$ \align &
(1-{\frac {az_2}{z_1}}):X(\vep_i,z_1)X(-\vep_i,az_2): \tag 2.33 \\
= & {\frac {(az_2)^{\frac 1 2}}{z_1^{\frac 1 2}}}\big(({\frac {z_1}
{az_2}})^{\vep_i}
E^-(-\vep_i,z_1)E^-(\vep_i,az_2)E^+(-\vep_i,z_1)E^+(\vep_i,az_2)-1\big).
\endalign $$
\endproclaim
 
\proclaim{Proposition 2.34} For $1\le i,j\le M$, and $a\in {\Bbb
C}^{\times}$, we have
$$ \align
 & :X(\vep_i,z)X(-\vep_j,az): \\
=& \cases & \epsilon(\vep_i,\vep_j)z^{\frac 1 2}(az)^{\frac 1 2}
e^{\vep_i-\vep_j}z^{\vep_i}(az)^{-\vep_j}
E^-(-\vep_i,z)E^-(\vep_j,az)E^+(-\vep_i,z)E^+(\vep_j,az) \\
& \text{if}\;\; i\not=
j\\
&
\vep_i(z)
\;\;\;\;\;\;\;\;\;\;\;\;\;\;\;\;\;\;\;\;\;\;\;\;\;\;\;\;\;\;\;\;\;
\text{if}\;\;  i=j, a=1\\
& {\frac
{a^{1/2}}{1-a}}(a^{-\vep_i}E^-(-\vep_i,z)E^-(\vep_i,az)E^+(-\vep_i,z)
E^+(\vep_i,az)-1)\\
& \text{if}\;\; i=j, a\not= 1.\endcases
\endalign $$
 \endproclaim
 
\demo{Proof} Taking the limit $z_2 \to z_1$ in (2.32) and (2.33) gives
the first and third
identities. The second identity follows from Lemma 2.22, Corollary 2.26 and 
the
third identity by taking
the limit $a\to 1$. \qed \enddemo
 
\noindent{\bf Remark 2.35.} Note that the second identity in Proposition
2.34 was given in
[F1].
 
\medskip
 
\proclaim{Definition 2.36} 
For $a\in {\Bbb C}^{\times}$, $1\le i,j\le M$, we define
 $
X_{ij}(a,z)=:X(\vep_i,z)X(-\vep_j,az):
$
\endproclaim
 \medskip
 
Now we can state our main theorem of this section.
 
\proclaim{Theorem 2.37} 
For $a,b\in {\Bbb C}^{\times}$ and $1\le i,j,k,l\le M$, we have
 
 (i) if $ab\not= 1$, then
$$\align &
[X_{ij}(a,z_1), X_{kl}(b,z_2)]=X_{il}(ab,z_1)\delta_{jk}\delta({\frac {z_2}
{az_1}})-X_{kj}(ab,z_2)\delta_{il}\delta({\frac {z_1} {bz_2}}) \tag 2.38
\\
& \;\;\;\;\;\;\;\;+{\frac {a^{\frac 1 2}b^{\frac 1
2}}{1-ab}}\delta_{il}\delta_{jk}(\delta({\frac {z_2}
{az_1}})
-\delta({\frac {z_1} {bz_2}}))c,
\endalign $$
 
(ii) if $ab=1$, then
$$ \align
 & [X_{ij}(a,z_1),
X_{kl}(b,z_2)]=(X_{il}(1,z_1)\delta_{jk}-X_{kj}(1,z_2)\delta_{il})\delta({\frac
{z_2} {az_1}})
+\delta_{il}\delta_{jk}(D\delta)({\frac {z_2} {az_1}})c.\tag 2.39 \endalign $$
 
\endproclaim
 
The proof of Theorem 2.37 will be carried out in several steps. 
In what follows we will
freely use the following two lemmas. (2.40) can be found in [FLM] and
[K],  and (2.43) can be found
in [J], [G1, 2] or [BS].
 
\proclaim{Lemma 2.40} Let $Y(z_1,z_2)$ be a formal power series in
$z_1,z_2$ with coefficients
in a vector
space, such that  $\lim_{z_2\to z_1}f(z_1,z_2)$ exists. Then
$$
Y(z_1,z_2)\delta({\frac {az_1}{z_2}})=Y(z_1,az_1)\delta({\frac
{az_1}{z_2}}), \tag 2.41
$$
$$
Y(z_1,z_2)(D\delta)({\frac {z_2}{az_1}})=Y(z_1,az_1)(D\delta)({\frac
{z_2}{az_1}})-
(D_{z_2}Y(z_1,z_2))\delta({\frac {z_2}{az_1}}),\tag 2.42
$$
for $a\in \Bbb{C}^\times$.
\endproclaim

\proclaim{Lemma 2.43} Suppose $a,b\in {\Bbb C}^{\times}$. Then
$$
(1-{\frac {z_2}{az_1}})^{-1}(1-{\frac {bz_2}{z_1}})^{-1}-{\frac
{az_1}{z_2}}{\frac {z_1}{
bz_2}}(1-{\frac {z_1}{bz_2}})^{-1}(1-{\frac {az_1}{z_2}})^{-1}
$$
$$
=\cases (1-ab)^{-1}{\frac {az_1}{z_2}}(\delta({\frac
{az_1}{z_2}})-\delta({\frac {z_1}{bz_2}})) &
\text{if}\;\; ab\not=
1\\
{\frac {az_1}{z_2}}(D\delta)({\frac {z_2}{az_1}}) &\text{if}\;\;
ab=1.\endcases$$
 
\endproclaim
 
\medskip
 
Now we divide the proof for Theorem 2.37 into four different cases.

\medskip
 
{\bf Case 1.} $i\not= j,$ $k\not= l$.
 
\medskip
 
We obtain, by applying (2.15) and (2.18), that
$$
\align & [X_{ij}(a, z_1),X_{kl}(b,
z_2)]\tag 2.44\\
= & \epsilon(\vep_i,\vep_j)\epsilon(\vep_k,\vep_l)a^{\frac 1 2}b^{\frac 1
2}
e^{\vep_i-\vep_j}e^{\vep_k-\vep_l}a^{-\vep_j}b^{-\vep_l}z_1^{\vep_i-\vep_j}
z_2^
{\vep_k-\vep_l}z_1z_2
 \\
& \cdot E^-(-\vep_i,z_1)E^-(-\vep_k,z_2)E^-(\vep_j,az_1)E^-(\vep_l,bz_2)
\\
& \cdot E^+(-\vep_i,z_1)E^+(-\vep_k,z_2)
E^+(\vep_j,az_1)E^+(\vep_l,bz_2)P(z_1,z_2)
\endalign $$
where
$$\align &
P(z_1,z_2)\\
=&a^{-(\vep_j,\vep_k-\vep_l)}z_1^{(\vep_i-\vep_j,\vep_k-\vep_l)}(1-{\frac
{z_2}{z_1}})^{\delta_
{ik}}(1-{\frac {z_2}{az_1}})^{-\delta_
{jk}}(1-{\frac {bz_2}{z_1}})^{-\delta_
{il}}(1-{\frac {bz_2}{az_1}})^{\delta_
{jl}} \endalign
$$
$$
-(-1)^{(\vep_i-\vep_j,\vep_k-\vep_l)}b^{-(\vep_l,\vep_i-\vep_j)}
z_2^{(\vep_k-\vep_l,\vep_i-\vep_j)}(1-
{\frac {z_1}{z_2}})^{\delta_
{ik}}(1-{\frac {z_1}{bz_2}})^{-\delta_
{il}}(1-{\frac {az_1}{z_2}})^{-\delta_
{jk}}(1-{\frac {az_1}{bz_2}})^{\delta_
{jl}}
$$
$$
=a^{-(\vep_j,\vep_k-\vep_l)}z_1^{(\vep_i-\vep_j,\vep_k-\vep_l)}(1-{\frac
{z_2}{z_1}})^{\delta_
{ik}}(1-{\frac {bz_2}{az_1}})^{\delta_
{jl}}\left((1-{\frac {z_2}{az_1}})^{-\delta_
{jk}}(1-{\frac {bz_2}{z_1}})^{-\delta_
{il}}\right.
$$
$$
-\left.(-1)^{\delta_{il}+\delta_{jk}}({\frac
{z_2}{az_1}})^{-\delta_{jk}}({\frac {bz_2}{z_1}})^{-\delta_{il}}
(1-{\frac {az_1}{z_2}})^{-\delta_
{jk}}(1-{\frac {z_1}{bz_2}})^{-\delta_
{il}}\right).
$$
Applying Lemma 2.43 we have two subcases.
 
\medskip
 
{\it Subcase 1}. If $ab\not=1$, then
$$
P(z_1,z_2)=\cases 0 &\text{if}\;\;i\not= l, j\not= k\\
{\frac 1{(1-ab)z_1z_2}}(\delta({\frac {z_2}{az_1}})-\delta({\frac
{z_1}{bz_2}})) &\text{if}\;\;i=l, j=k\\
z_1^{-1}\delta({\frac {z_1}{bz_2}}) &\text{if}\;\;i=l, j\not= k\\
(az_1)^{-1}\delta({\frac {z_2}{az_1}}) &\text{if}\;\;i\not= l, j=
k.\endcases
$$
Now (2.38) follows from (2.44) and Lemma 2.40.
 
\medskip
 
{\it Subcase 2}. If $ab=1$, then
$$
P(z_1,z_2)=a^{-(\vep_j,\vep_k-\vep_l)}z_1^{(\vep_i-\vep_j,\vep_k-\vep_l)}
(1-{\frac {z_2}{z_1}})^{\delta_
{ik}}(1-{\frac {bz_2}{az_1}})^{\delta_
{jl}}((1-{\frac {z_2}{az_1}})^{-\delta_
{jk}-\delta_{il}})
$$
$$
-(-1)^{\delta_{il}+\delta_{jk}}({\frac
{z_2}{az_1}})^{-\delta_{jk}-\delta_{il}}
(1-{\frac {az_1}{z_2}})^{-\delta_
{jk}-\delta_{il}}
$$
$$
=\cases 0 &\text{if}\;\;i\not= l, j\not= k\\
{\frac 1{z_1z_2}}(D\delta)({\frac {z_2}{az_1}}) &\text{if}\;\;i=l, j=k\\
z_1^{-1}\delta({\frac {z_2}{az_1}}) &\text{if}\;\;i=l, j\not= k\\
(az_1)^{-1}\delta({\frac {z_2}{az_1}}) &\text{if}\;\;i\not= l, j=
k,\endcases
$$
which yields (2.39) by applying Lemma 2.40.
 
\medskip
 
{\bf Case 2}. $i= j,$ $k\not= l$, and $ab\not= 1$.
 
\medskip
 
If $a=1$, then
$$
X_{ij}(a, z_1)=\vep_i(z_1).
$$
Since
$$\align
& [\vep_i(z_1),
e^{\vep_k-\vep_l}]=(\delta_{ik}-\delta_{il})e^{\vep_k-\vep_l},\\
& [\vep_i(z_1),
E^{\pm}(-\vep_k,z_2)]=\sum_{n\in{\Bbb Z}_{\pm}}\delta_{ik}({\frac
{z_1}{z_2}})^nE^{\pm}(-\vep_k,z_2),
\endalign $$
and
$$
[\vep_i(z_1),E^{\pm}(\vep_l,bz_2)]=-\sum_{n\in{\Bbb
Z}_{\pm}}\delta_{il}({\frac {z_1}{bz_2}})^n
E^{\pm}(\vep_l,bz_2),
$$
 we have
$$
  [X_{ij}(a, z_1),X_{kl}(b, z_2)]=[\vep_i(z_1),X_{kl}(b, z_2)]=
X_{kl}(b, z_2)
$$
$$\cdot\left(\delta_{ik}-\delta_{il}+
\delta_{ik}\sum_{n\in{\Bbb Z}_{-}}({\frac {z_1}{z_2}})^n+
\delta_{ik}\sum_{n\in{\Bbb Z}_{+}}({\frac
{z_1}{z_2}})^n-\delta_{il}\sum_{n\in{\Bbb Z}_{-}}
({\frac {z_1}{bz_2}})^n
-\delta_{il}\sum_{n\in{\Bbb Z}_{+}}({\frac {z_1}{bz_2}})^n\right)
$$
$$ \align
= & X_{kl}(b, z_2)(\delta_{ik}\delta({\frac
{z_2}{z_1}})-\delta_{il}\delta({\frac
{z_1}{bz_2}})) \\
= &X_{kl}(b, z_1)\delta_{jk}\delta({\frac {z_2}{z_1}})-X_{kj}(b,
z_2)\delta_{il}\delta({\frac
{z_1}{bz_2}}),
\endalign $$
as needed.
 
If $a\not= 1$, then
$$
X_{ij}(a, z_1)={\frac
{a^{1/2}}{1-a}}(a^{-\vep_i}E^-(-\vep_i,z_1)E^-(\vep_i,az_1)
E^+(-\vep_i,z_1)E^+(\vep_i,az_1)-1),
$$
Applying (2.15) and (2.18), we get
$$\align
& [X_{ij}(a, z_1),X_{kl}(b, z_2)]\tag 2.45 \\
=& {\frac {a^{1/2}}{1-a}}b^{\frac 1
2}\epsilon(\vep_k,\vep_l)e^{\vep_k-\vep_l}
a^{-\vep_i}b^{-\vep_l}z_2^{\vep_k-\vep_l}z_2 \\
& \cdot
E^-(-\vep_i,z_1)E^-(\vep_i,az_1)E^-(-\vep_k,z_2)E^-(\vep_l,bz_2)\\
& \cdot E^+(-\vep_i,z_1)E^+(\vep_i,az_1)
E^+(-\vep_k,z_2)E^+(\vep_l,bz_2)Q(z_1,z_2),
\endalign $$
where
$$\align &
Q(z_1,z_2)=a^{-(\vep_j,\vep_k-\vep_l)}(1-{\frac {z_2}{z_1}})^{\delta_
{ik}}(1-{\frac {z_2}{az_1}})^{-\delta_
{ik}}(1-{\frac {bz_2}{z_1}})^{-\delta_
{il}}(1-{\frac {bz_2}{az_1}})^{\delta_
{il}}\endalign
$$
$$
-(1-
{\frac {z_1}{z_2}})^{\delta_
{ik}}(1-{\frac {z_1}{bz_2}})^{-\delta_
{il}}(1-{\frac {az_1}{z_2}})^{-\delta_
{ik}}(1-{\frac {az_1}{bz_2}})^{\delta_
{il}}
$$
$$
=a^{\delta_{il}-\delta_{ik}}(1-{\frac {z_2}{z_1}})^{\delta_
{ik}}(1-{\frac {bz_2}{az_1}})^{\delta_
{il}}\left((1-{\frac {z_2}{az_1}})^{-\delta_
{ik}}(1-{\frac {bz_2}{z_1}})^{-\delta_
{il}}\right.
$$
$$
\left.-(-1)^{\delta_{ik}+\delta_{il}}({\frac {az_1}{z_2}})^{\delta_{ik}}
({\frac {z_1}{bz_2}})^{\delta_{il}}(1-{\frac {z_1}{bz_2}})^{-\delta_
{il}}(1-{\frac {az_1}{z_2}})^{-\delta_
{ik}}\right)
$$
$$
=\cases 0 &\text{if}\;\;i\not= k, i\not= l\\
{\frac {1-a}a}\delta({\frac {z_2}{az_1}}) &\text{if}\;\;i= k, i\not= l\\
(a-1)\delta({\frac {z_1}{bz_2}}) &\text{if}\;\;i\not= k, i= l\endcases
$$
 We thus have (2.38) by applying Lemma 2.40.
 
\medskip
 
{\bf Case 3}. $i= j,$ $k\not= l$, and $ab= 1$.
 
\medskip
 
 As we did in Case 2, we  treat the case $a=1$ and $a\not= 1$
separately. First for
$a=1$ (so $b=1$),
 then
 $$
 X_{ij}(a,z_1)=\vep_i(z_1),
 $$
 $$
 X_{kl}(b,z_2)=\epsilon(\vep_k,\vep_l)
e^{\vep_k-\vep_l}z_2^{\vep_k-\vep_l}z_2
E^-(-\vep_k,z_2)E^-(\vep_l,z_2)E^+(-\vep_k,z_2)E^+(\vep_l,z_2),
$$
and we have
$$\align &
[X_{ij}(a,z_1),X_{kl}(b,z_2)]=[\vep_i(z_1),X_{kl}(1,z_2)]\\
=& X_{kl}(1,z_2)(\delta_{ik}\delta({\frac
{z_2}
{z_1}})-\delta_{il}\delta({\frac {z_2}
{z_1}}))\\
=&(X_{il}(1,z_1)\delta_{jk}-X_{kj}(1,z_2)\delta_{il})\delta({\frac {z_2}
{z_1}}),\endalign
$$
as expected.
 
Next assume that $a\not= 1$. By applying similar arguments as  in
Case 2, we have
$$
\align & [X_{ij}(a,z_1),X_{kl}(b,z_2)]\tag 2.46\\
=& {\frac 1
{1-a}}\epsilon(\vep_k,\vep_l)e^{\vep_k-\vep_l}a^{-\vep_i}a^{\vep_l}
z_2^{\vep_k-\vep_l}z_2\\
& \cdot E^-(-\vep_i,z_1)
 E^-(\vep_i,az_1)E^-(-\vep_k,z_2)E^-(\vep_l,bz_2)
\\
& \cdot E^+(-\vep_i,z_1)E^+(\vep_i,az_1)
 E^+(-\vep_k,z_2)E^+(\vep_l,bz_2)R(z_1,z_2)
\endalign $$
where
$$\align &
R(z_1,z_2)=a^{\delta_{il}-\delta_{ik}}(1-{\frac
{z_2}{z_1}})^{\delta_{ik}}(1-{\frac {bz_2}{az_1}})^{
\delta_{il}}\left((1-{\frac
{z_2}{az_1}})^{-\delta_{il}-\delta_{ik}}\right.
\endalign$$
$$
\left.-(-1)^{\delta_{il}+\delta_{ik}}({\frac
{az_1}{z_2}})^{\delta_{il}+\delta_{ik}}
(1-{\frac {az_1}{z_2}})^{-\delta_{il}-\delta_{ik}}\right)
$$
$$
=\cases 0 &\text{if}\;\;i\not= k, i\not= l\\
(b-1)\delta({\frac {z_2}{az_1}}) &\text{if}\;\;i= k, i\not= l\\
(a-1)\delta({\frac {z_2}{az_1}}) &\text{if}\;\;i\not= k, i= l.\endcases
$$
Substituting this into (2.46) and applying Lemma 2.40, we get (2.39).
 
\medskip
 
{\bf Case 4}. $i=j, k=l$.
 
\medskip
 
Note that
$$
X_{ij}(a,z_1)
$$
$$
=\cases {\frac
{a^{1/2}}{1-a}}(a^{-\vep_i}E^-(-\vep_i,z_1)E^-(\vep_i,az_1)E^+(-\vep_i,
z_1)E^+(\vep_i,az_1)-1) &
\text{if}\;\; i=j, a\not= 1\\
\vep_i(z_1) &\text{if}\;\;  i=j, a=1
\endcases
$$
$$
X_{kl}(b,z_2)
$$
$$
=\cases {\frac
{b^{1/2}}{1-b}}(b^{-\vep_k}E^-(-\vep_k,z_2)E^-(\vep_k,bz_2)E^+(-\vep_k,
z_2)E^+(\vep_k,bz_2)-1) &
\text{if}\;\; k=l, b\not= 1\\
\vep_k(z_2) &\text{if}\;\;  k=l, b=1.
\endcases
$$
 
We consider three subcases. First we assume that $a=b=1$, then
$$[X_{ij}(a,z_1),
X_{kl}(b,z_2)]  = [\vep_i(z_1),\vep_k(z_2)] =\delta_{ik}(D\delta)({\frac
{z_2}{z_1}}),$$ as desired.
 
Next we assume $a= 1$, $b\not= 1$, then
$$
[X_{ij}(a,z_1),X_{kl}(b,z_2)]=[\vep_i(z_1),X_{kl}(b,z_2)]
$$
$$
={\frac
{b^{1/2}}{1-b}}b^{-\vep_k}E^-(-\vep_k,z_2)E^-(\vep_k,bz_2)E^+(-\vep_k,
z_2)E^+(\vep_k,bz_2)
$$
$$\cdot\left(\delta_{ik}\sum_{n\in {\Bbb Z}_{-}}({\frac
{z_2}{z_1}})^n
+
\delta_{ik}\sum_{n\in {\Bbb Z}_{+}}({\frac
{z_2}{z_1}})^n-\delta_{ik}\sum_{n\in {\Bbb Z}_{-}}(
{\frac
{z_1}{bz_2}})^n-
\delta_{ik}\sum_{n\in {\Bbb Z}_{+}}({\frac {z_1}{bz_2}})^n\right)
$$
$$
={\frac
{b^{1/2}}{1-b}}b^{-\vep_k}E^-(-\vep_k,z_2)E^-(\vep_k,bz_2)E^+(-\vep_k,
z_2)E^+(\vep_k,bz_2)\delta_{ik}(\delta({\frac {z_2}{z_1}})-\delta({\frac
{z_1}{bz_2}}))
$$
$$
=X_{il}(ab,z_1)\delta_{jk}\delta({\frac {z_2}
{az_1}})-X_{kj}(ab,z_2)\delta_{il}\delta({\frac {z_1} {bz_2}})
+{\frac {a^{\frac 1 2}b^{\frac 1
2}}{1-ab}}\delta_{il}\delta_{jk}(\delta({\frac {z_2} {az_1}})
-\delta({\frac {z_1} {bz_2}})),
$$
as required.
 
Finally we assume $a\not= 1$, $b\not= 1$, then
$$
[X_{ij}(a,z_1),X_{kl}(b,z_2)]\tag 2.47
$$
$$
={\frac {a^{1/2}}{1-a}}{\frac {b^{1/2}}{1-b}}a^{-\vep_i}b^{-\vep_k}
E^-(-\vep_i,z_1)E^-(-\vep_k,z_2)E^-(\vep_i,az_1)
$$
$$
\cdot E^-(\vep_k,bz_2)E^+(-\vep_i,z_1)E^+(-\vep_k,z_2)E^+(\vep_i,az_1)E^+(\vep_k,bz_2
S(z_1,z_2),
$$
where
$$\align &
S(z_1,z_2)=(1-{\frac {z_2}{z_1}})^{\delta_{ik}}(1-{\frac
{bz_2}{z_1}})^{-\delta_{ik}}
(1-{\frac {z_2}{az_1}})^{-\delta_{ik}}
(1-{\frac {bz_2}{az_1}})^{\delta_{ik}}
\endalign $$
$$
-(1-{\frac {z_1}{z_2}})^{\delta_{ik}}(1-{\frac
{z_1}{bz_2}})^{-\delta_{ik}}
(1-{\frac {az_1}{z_2}})^{-\delta_{ik}}
(1-{\frac {az_1}{bz_2}})^{\delta_{ik}}
$$
$$
=(1-{\frac {z_2}{z_1}})^{\delta_{ik}}(1-{\frac
{bz_2}{az_1}})^{\delta_{ik}}
\left((1-{\frac {bz_2}{z_1}})^{-\delta_{ik}}
(1-{\frac {z_2}{az_1}})^{-\delta_{ik}}\right.
$$
$$
\left.-({\frac {az_1}{z_2}})^{\delta_{ik}}({\frac
{z_1}{bz_2}})^{\delta_{ik}}
(1-{\frac {z_1}{bz_2}})^{-\delta_{ik}}
(1-{\frac {az_1}{z_2}})^{-\delta_{ik}}\right)
$$
$$
=\cases (1-{\frac {z_2}{z_1}})(1-{\frac {bz_2}{az_1}}){\frac
{az_1}{z_2}}(D\delta)({\frac {z_2}{az_1}})
 &\text{if}\;\; i=k, ab= 1\\
0 &\text{if}\;\;i\not=k\\
{\frac 1 {1-ab}}(1-{\frac {z_2}{z_1}})(1-{\frac {bz_2}{az_1}}){\frac
{az_1}{z_2}}(\delta
({\frac {z_2}{az_1}})-\delta
({\frac {z_1}{bz_2}})) &\text{if}\;\; i=k, ab\not= 1.
\endcases
$$
Substitute $S(z_1, z_2)$ back into (2.47) and apply Lemma 2.40 to get
(2.38) and (2.39).
 This now completes the proof of Theorem 2.37.
 
\medskip
 
\noindent{\bf Remark 2.48} If $ab\not= 1$, we have
$$
\delta({\frac {z_2} {az_1}})
-\delta({\frac {z_1} {bz_2}})=\sum_{m\in{\Bbb Z}}(1-(ab)^m)({\frac {z_2}
{az_1}})^m,
$$
thus
$$
{\frac {a^{\frac 1 2}b^{\frac 1
2}}{1-ab}}\delta_{il}\delta_{jk}(\delta({\frac {z_2} {az_1}})
-\delta({\frac {z_1} {bz_2}}))=a^{\frac 1 2}b^{\frac 1
2}\delta_{il}\delta_{jk}
\sum_{m\in{\Bbb Z}}{\frac {1-(ab)^m}{1-ab}}({\frac {z_2} {az_1}})^m
$$
$$
=a^{\frac 1 2}b^{\frac 1 2}\delta_{il}\delta_{jk}(\sum_{m\in{\Bbb
Z}_+}(\sum^{m-1}_{s=0}(ab)
^s)
({\frac {z_2} {az_1}})^m+\sum_{m\in {\Bbb
Z}_-}(-\sum^{-m}_{s=1}(ab)^{-s})
({\frac {z_2} {az_1}})^m),
$$
which gives
$$
\lim_{b\to a^{-1}}{\frac {a^{\frac 1 2}b^{\frac 1
2}}{1-ab}}\delta_{il}\delta_{jk}
(\delta({\frac {z_2} {az_1}})
-\delta({\frac {z_1} {bz_2}}))
$$
$$
=\delta_{il}\delta_{jk}(\sum_{m\in{\Bbb Z}_+}m
({\frac {z_2} {az_1}})^m+\sum_{m\in {\Bbb Z}_-}m({\frac {z_2}
{az_1}})^m)=
\delta_{il}\delta_{jk}
(D\delta)({\frac {z_2} {az_1}}).
$$
This indicates that the second identity of Theorem 2.37 can be obtained
from the first one by   
taking the limit
as $b\to a^{-1}$.
 
\medskip
 
\subhead \S 3. Lie Algebras and Representations \endsubhead
 
 \medskip

In this section we are going to define a class of Lie algebras from our
vertex operators which will correspond to  admissible subgroups of $\bc^{\times}$. Indeed for some
choice of the admissible 
group $G$ and positive integer $M$, the Lie algebra ${\Cal G}(G,M)$
(defined below) of operators,
which act on the
Fock space $V_M$, will give realizations of some infinite dimensional Lie
algebras studied in Section 1.  This will include the affine
algebra ${\wh {gl_M}}({\Bbb C})$ in both the principal and homogeneous
pictures as well as  some Lie algebras
with quantum torus coordinates. 
Towards this end, we first introduce some new notation for the vertex operators constructed in the
proceeding section.
 
\medskip
 
\proclaim{Definition 3.1}  For $a,b\in {\Bbb C}^{\times}$, $1\le i,j\le
M$, we set
$X_{ij}(a,b,z):=X_{ij}(
 a^{-1}b, az)$, and write
 $$
 X_{ij}(a,b,z)=\sum_{k\in {\Bbb Z}}x_{ij}(k,a,b)z^{-k} \tag 3.2 $$
 where $x_{ij}(k,a,b)\in$ End$V_M$.
\endproclaim
 
With this notation Theorem 2.37 can be re-written as follows.
 
\medskip
 
\proclaim{Theorem 3.3}  
Let  $a_1,a_2,b_1,b_2\in {\Bbb C}^{\times}$, and $1\le
i,j,k,l\le M$. We have
 
\medskip
 
 (i) if $a_1a_2\not= b_1b_2$, then
$$\align &
[X_{ij}(a_1,b_1,z_1), X_{kl}(a_2,b_2,z_2)]\\
=& X_{il}(a_1,{\frac
{b_1b_2}{a_2}},z_1)\delta_{jk}\delta({\frac {a_2z_2} {b_1z_1}})
-
X_{kj}(a_2,{\frac {b_1b_2}{a_1}},z_2)\delta_{il}
\delta({\frac {a_1z_1} {b_2z_2}})\\
& +\frac{(a_1^{-1}b_1)^{\frac 1 2}(a_2^{-1}b_2)^{frac 1
2}}{1-a_1^{-1}b_1a_2^{-1}b_2}
\delta_{il}\delta_{jk}(\delta({\frac {a_2z_2} {b_1z_1}})
-\delta({\frac {a_1z_1} {b_2z_2}}))c,
\endalign $$
 
(ii) if $a_1a_2= b_1b_2$, then
$$
[X_{ij}(a_1,b_1,z_1), X_{kl}(a_2,b_2,z_2)]
$$
$$
=X_{il}(a_1,{\frac {b_1b_2}{a_2}},z_1)\delta_{jk}\delta({\frac {a_2z_2}
{b_1z_1}})
-
X_{kj}(a_2,{\frac {b_1b_2}{a_1}},z_2)\delta_{il}\delta({\frac {a_1z_1}
{b_2z_2}})
+
\delta_{il}\delta_{jk}(D\delta)({\frac {a_2z_2} {b_1z_1}})c.$$
 
\endproclaim
 
\medskip

Fix an integer $M\ge 1$ and an admissible subgroup $G$ of ${\Bbb
C}^{\times}$.
 Let
${\Cal G}(G,M)$ be the vector space spanned by $c$ and all of the
coefficients of the vertex operators
$X_{ij}(a,b,z)$ for all $1\le i,j\le M$, and $a,b\in G$. Then
we have the following result.
 
\medskip
 
\proclaim{Theorem 3.4} ${\Cal G}(G,M)$
 forms a Lie algebra of operators acting on the Fock space  $V_M$.
Moreover
 $$
 V_M=\oplus_{k\in{\Bbb Z}}V_M^{(k)}
 $$
 where $V_M^{(k)}=e^{k\vep_M+Q_M}\otimes{\Cal S}({\Cal H}^-_M)$, and
$V_M^{(k)}$ is an irreducible
 ${\Cal G}(G,M)$-module.
 \endproclaim
 
\demo{Proof}  It is obvious from Theorem 3.3 that ${\Cal G}(G,M)$ is a Lie algebra and that $V_M^{(k)}$ is a ${\Cal
G}(G,M)$-module.
To see it is
irreducible,
 we note that the Heisenberg algebra ${\wh {\Cal H}}_M\subset$ ${\Cal
H}_M$, and
 ${\Cal H}_M$ is spanned by the coefficient operators of the vertex
operators $X_{ii}(1,1,z)$ for $
 1\le i\le M$. This then implies that,
if $W$ is a non-zero submodule of $V_M^{(k)}$, we can choose
a non-zero element of the form $v=e^{k\vep_M+\alpha}\otimes 1\in
W$ for some $\alpha\in Q_M$.
 Moreover, it is easy to check that
 $$
 x_{ij}(n_{ij}-1,1,1).v=\epsilon(\vep_i,\vep_j)\epsilon(\vep_i-\vep_j,
 k\vep_M+\alpha)e^{k\vep_M+\alpha+\vep_i-\vep_j}
 $$
 for all $1\le i\not= j\le M$, where $n_{ij}=(\vep_j-\vep_i,
k\vep_M+\alpha)\in {\Bbb
Z}$. Therefore $
 e^{k\vep_M+\beta}\otimes 1\in W$ for any $\beta\in Q_M$. This
thus gives $W=V_M^{(k)}
 $ as needed.\qed
 
\enddemo
 
 \medskip
 
\noindent{\bf Remark 3.5} Note that the coefficients of the vertex
operators $X_{ij}(a, z)$ and $X_{ij}(a, bz)$ for any given $a, b\in G$
span the same space. Thus $\Cal{G}(G, M)$ is spanned by $c$ and the
coefficients of the operators $X_{ij}(a, z)$ for $1\leq i, j\leq M, a\in
G$. Therefore, it follows from Theorem 4.25 in [G1] that $\Cal{G}(G, M)$
is an affinization of the Lie algebra $gl_M(\Cal{R}[t, t^{-1}; \tau])$,
where
$\Cal{R}= \Bbb{C}[G]$ is the group algebra and $\Cal{R}[t, t^{-1}; \tau]$
is the skew Laurent
polynomial ring.
 
Recall definition (2.19). We extend the cocycle map $\epsilon:$ $H_M\times
H_M\to\{\bc^{\times}\}$ by
defining
$$
\epsilon(\sum r_i\vep_i,\sum
s_i\vep_i)=\prod_{i,j}(\epsilon(\vep_i,\vep_j))^{r_is_j}\tag 3.6
$$
for $r_i,s_i\in \bc$. It is obvious that
$$
\epsilon(\alpha+\beta,\gamma)=\epsilon(\alpha,\gamma)\epsilon(\beta,\gamma),\;\
\;\;\epsilon(\alpha, \beta+\gamma)
=\epsilon(\alpha,\beta)\epsilon(\alpha,\gamma)
$$
for $\alpha,\beta,\gamma\in H_M$. Moreover, if we restrict $\epsilon$ to
$\Gamma_M\times\Gamma_M$, then $\epsilon$ gives us
the 2-cocycle defined in previous section.
 
  Let $\alpha_i=\vep_i-\vep_{i+1},$ $i=1,\cdots, M-1$, and
$\alpha_M=\vep_1+\cdots+\vep_M$. Then $Q_M=\oplus_{i=1}^{M-1}{\Bbb
Z}\alpha_i$, and $\Gamma_M=\oplus_{i=1}^{M}{\Bbb Z}\alpha_i$. Let
$Q^0_M=\{\alpha\in\bc\otimes_{\Bbb Z}Q_M|\;\;(\alpha, Q_M)\in {\Bbb Z}
\}$ be the dual of the lattice $Q_M$, and set
$$
L^0_M=\{\alpha\in H_M=\bc\otimes_{\Bbb Z}\Gamma_M|\;\;(\alpha, Q_M)\in
{\Bbb Z}
\}.\tag 3.7
$$
Then $L^0_M=Q^0_M\oplus {\bc}\alpha_M$. Let $I=L^0_M/Q_M$, then we have a
$Q_M$-coset decomposition of $L^0_M$
$$
L^0_M=\oplus_{i\in I}(\lambda_i+Q_M)\tag 3.8
$$
for some $\lambda_i\in L^0_M.$

\medskip
 \proclaim
{\bf Proposition 3.9} The Lie algebra of operators, ${\Cal G}(G,M)$,  acts on
the space $V^{(\lambda_i)}_M=$ ${\Cal S}({\Cal H}^-_M)\otimes
{\bc}[Q_M+\lambda_i]$, and $V^{(\lambda_i)}_M$ affords an irreducible
representation of ${\Cal G}(G,M)$ for $i\in I=L^0_M/Q_M$.
Moreover $V^{(\lambda_i)}_M\cong V^{(\lambda_j)}_M$ if and only if $i=j$.
 \endproclaim

\demo{Proof} For $u\otimes e^{\alpha+\lambda_i}\in V^{(\lambda_i)}_M$, $u\in
{\Cal S}({\Cal H}^-_M),$
$\alpha\in Q_M$, we note that
$$
z^{\beta}.(u\otimes
e^{\alpha+\lambda_i})=z^{(\beta,\alpha+\lambda_i)}u\otimes
e^{\alpha+\lambda_i},\;\;\;\;a^{\gamma}.(u\otimes 
e^{\alpha+\lambda_i})=
a^{(\gamma,\alpha+\lambda_i)}u\otimes e^{\alpha+\lambda_i}
$$
$$
e^{\beta}.(u\otimes
e^{\alpha+\lambda_i})=\epsilon(\beta,\alpha+\lambda_i)u\otimes
e^{\alpha+\beta+\lambda_i}
$$
where $(\beta,\alpha+\lambda_i)\in{\Bbb Z},$ $
(\gamma,\alpha+\lambda_i)\in {\bc},$ $\epsilon(\beta,\alpha+\lambda_i)
\in{\bc}^{\times}$ for $\beta\in Q_M,$ $\gamma\in \Gamma_M$ and $a\in G$.
This implies that $X_{ij}(a,b,z).(
u\otimes e^{\alpha+\lambda_i})\in V^{(\lambda_i)}_M[[z,z^{-1}]],$ and so
the Lie algebra ${\Cal G}(G,M)$
acts on the space $V^{(\lambda_i)}_M$. The irreducibility of 
$V^{(\lambda_i)}_M$, for $i\in I$, follows from a similar argument
as  in Theorem 3.4. The last part of this proposition is clear. \qed
 \enddemo

\bigskip
 
\subhead \S 4. Applications \endsubhead
 
 \medskip

 In this section we assume  the admissible subgroup $G$ has the form
$G=T\times F\subset 
{\bc}^{\times}$, where $T=<\xi>$ is generated
by a root of unity $\xi$, and $F$ is a free group with a finite number generators.
First, let $G=\{1\}$, and $M\ge 2$ be any integer. Then the Lie algebra
${\Cal G}(G,M)$ is generated
by the coefficients of the vertex operators $X_{ij}(1,1,z)$ for all $1\le
i,j\le M$. Moreover
from Theorem 3.3 , we see that
$$\align &
[X_{ij}(1,1,z_1),
X_{kl}(1,1,z_2)] \tag 4.1\\
=& X_{il}(1,1,z_1)\delta_{jk}\delta({\frac {z_2}{z_1}})-
X_{kj}(1,1,z_2)\delta_{il}\delta({\frac
{z_1}{z_2}})+\delta_{il}\delta_{jk}(D\delta)({\frac
{z_2}{z_1}}) \endalign$$
for $1\le i,j\le M$. Comparing (4.1) with  (1.17), we obtain the following
result which was originally due to
[F1], see also [FK] and [S].
 
\medskip
 
\proclaim{Corollary 4.2} Let $G=\{1\}$ and $M\ge 2$. Then ${\Cal G}(G,M)$
gives a
representation of the affine
algebra ${\wh {gl_M}}({\Bbb C})$ in the homogeneous picture on the Fock
space $V_M$, and the representation
is given by the mapping:
$$\align &
E_{ij}\otimes t^k_0\mapsto x_{ij}(k,1,1), \\
& c_0\mapsto c
\endalign $$
for $1\le i,j\le M$ and $k\in {\Bbb Z}$.
\endproclaim
 
\medskip
 
Next we choose $G$ to be a cyclic group of order $N\ge 2$ with generator
$\xi=\xi_N$, and take $M=1$.  Note that
$$X_{11}(\xi^i,\xi^j,z)=X_{11}(\xi^{i-j-1}, \xi^{-1}, \xi^{j+1}z)$$
 for
$0\le i,j\le N-1$. This
implies that the Lie algebra ${\Cal G}(G,M)$ is generated
by the coefficients of the vertex operators
$X_{11}(\xi^{i-1},\xi^{-1},z)$ for $0\le
i\le N-1$.
 From Theorem 3.3 we have
 $$
 [X_{11}(\xi^{i-1},\xi^{-1},z_1), X_{11}(\xi^{j-1},\xi^{-1},z_2)]
 $$
 $$
 =\cases X_{11}(\xi^{i-1},\xi^{-j-1}, z_1)\delta({\frac
{\xi^{j}z_2}{z_1}})-
 X_{11}(\xi^{j-1},\xi^{-i-1}, z_2)\delta({\frac
{\xi^{i}z_1}{z_2}})
+
 (D\delta)({\frac {\xi^{j}z_2}{z_1}})c\\
 \;\;\;\;\;\;\;\;\;\;\;\;\;\;\;\;\;\;\;\;\;\;\;\;\;\;\;\;\text{if}\;
i+j=0(\text{mod}N)\\
 X_{11}(\xi^{i-1},\xi^{-j-1}, z_1)\delta({\frac
{\xi^{j}z_2}{z_1}})-
 X_{11}(\xi^{j-1},\xi^{-i-1}, z_2)\delta({\frac
{\xi^{i}z_1}{z_2}})+\\
 \;\;\;\;\;\;\;\;\;\;\;\;\;\;\;\;\;\;\;\;\;\;\;
{\frac {e^{-{\frac {i+j}2}Ln\xi}}{\xi^{i+j}-1}}
 (\delta({\frac {\xi^{j}z_2}{z_1}})-\delta({\frac {\xi^{i}z_1}{z_2}}))
 \;\;\;\;\;\text{if}\; i+j\not=
 0(\text{mod}N).\endcases
 $$
 Recalling the definition of $X_{ij}(a,b,z)$, we have
 $$
 X_{11}(\xi^{i-1},\xi^{-j-1},
z_1)\delta({\frac{\xi^{j}z_2}{z_1}})=
 X_{11}(\xi^{i-1},\xi^{-j-1},
\xi^{j}z_2)\delta({\frac{\xi^{j}z_2}{z_1}})
$$
$$
 =X_{11}(\xi^{i+j-1},\xi^{-1}, z_2)\delta({\frac
{\xi^{j}z_2}{z_1}}),
 $$
while
 $$
 X_{11}(\xi^{j-1},\xi^{-i-1}, z_2)\delta({\frac
{\xi^{i}z_1}{z_2}})
=
 X_{11}(\xi^{i-1},\xi^{-j-1}, \xi^{j}z_1)\delta({\frac
{\xi^{i}z_1}{z_2}})
 $$
$$
=X_{11}(\xi^{i+j-1},\xi^{-1}, z_1)\delta({\frac
{\xi^{i}z_1}{z_2}})
 $$
 for $0\le i,j\le N-1$. Therefore we get
 $$
 [X_{11}(\xi^{i-1},\xi^{-1},z_1), X_{11}(\xi^{j-1},\xi^{-1},z_2)] \tag
4.3
 $$
 $$
 =\cases X_{11}(\xi^{i+j-1},\xi^{-1}, z_2)\delta({\frac
{\xi^{j}z_2}{z_1}})-
 X_{11}(\xi^{i+j-1},\xi^{-1}, z_1)\delta({\frac
{\xi^{i}z_1}{z_2}})+
 (D\delta)({\frac {\xi^{j}z_2}{z_1}})c\\
  \;\;\;\;\;\;\;\;\;\;\;\;\;\;\;\;\;\;\;\;\;\;\;\;\;\;\;\;\text{if}\;
i+j=0(\text{mod}N)\\
 X_{11}(\xi^{i+j-1},\xi^{-1}, z_2)\delta({\frac
{\xi^{j}z_2}{z_1}})-
 X_{11}(\xi^{i+j-1},\xi^{-1}, z_1)\delta({\frac
{\xi^{i}z_1}{z_2}})+\\
 \;\;\;\;\;\;\;\;\;\;\;\;\;\;\;\;\;\;\;\;\;\;\;\;
{\frac {e^{-{\frac
 {i+j}2}Ln\xi}}{\xi^{i+j}-1}}
 (\delta({\frac {\xi^{j}z_2}{z_1}})-\delta({\frac {\xi^{i}z_1}{z_2}})
  \;\;\;\;\text{if}\; i+j\not=
 0(\text{mod}N)c.\endcases
 $$
 Comparing this with the identity (1.18), we obtain
 the following result which was originally due to [F1] and [KKLW].
 
\medskip
 
 \proclaim{Corollary 4.4}  Let $M=1$, and let $G$ be the group generated by
$\xi$, where $\xi$
is a $N$-th
 primitive root of unity for $N\ge 2$. Then the Lie algebra ${\Cal
G}(G,M)$ gives a representation
 of the affine algebra ${\wh {gl_N}}({\Bbb C})$ on the Fock space $V_1$
in the principal picture. The
 representation is given by the mapping
 $$\align &
 F^iE^k\otimes t^k_0\mapsto
 \cases x_{11}(k, \xi^{i-1},\xi^{-1})+
{\frac {e^{{\frac i
2}Ln\xi}}{\xi^i-1}}\delta_{k0}c &\text{if}\; 1\le i\le
 N-1\\x_{11}(k, \xi^{-1},\xi^{-1}) &\text{if}\; i=0\endcases, \\
& c_0\mapsto
{\frac cN}
  \endalign $$
  for $k\in {\Bbb Z}$.
  \endproclaim
 
  \medskip
 
Next we choose $M\ge 2$ and $G=<q>$ where $q\not= 0$ is not a root of
unity. Note that
$X_{ij}(a,b,z)=X_{ij}(1,a^{-1}b,az)$
for $a,b\in G$. We see that the Lie
algebra ${\Cal G}(G,M)$ is generated by the coefficients of the vertex
operators of the form $X_{ij}(1,
q^r,z)$ for all $r\in {\Bbb Z}$ and $1\le i,j\le M$. We apply Theorem 3.3
to obtain
 
\medskip
 
(i) if $r+s=0$, then
$$
[X_{ij}(1,q^r,z_1),X_{kl}(1,q^s,z_2)] \tag 4.5
$$
$$
= X_{il}(1,q^{r+s},z_1)\delta_{jk}\delta({\frac
{z_2}{q^rz_1}})-
X_{kj}(1,q^{r+s},z_2)\delta_{il}\delta({\frac {z_1}{q^sz_2}})
+\delta_{il}\delta_{jk}(D\delta)({\frac {z_2}{q^rz_1}})c,
$$
 
(ii) if $r+s\not= 0$, then
$$\align &
[X_{ij}(1,q^r,z_1),X_{kl}(1,q^s,z_2)]\tag 4.6\\
=& X_{il}(1,q^{r+s},z_1)\delta_{jk}\delta({\frac
{z_2}{q^rz_1}})-
X_{kj}(1,q^{r+s},z_2)\delta_{il}\delta({\frac {z_1}{q^sz_2}})\\
&+{\frac {q^{\frac
{r+s}2}}{1-q^{r+s}}}\delta_{il}\delta_{jk}(\delta({\frac {z_2}{q^rz_1}})-
\delta({\frac {z_1}{q^sz_2}}))c.
\endalign $$
Comparing the above two identities with the identity (1.19), we derive the
following result which
was given in [G1].
 
\medskip
 
\proclaim{Corollary 4.7}  Let $M\ge 2$, and $G=<q>$ be the group
generated by  $q\not= 0$
and
$q$ is not a root
of unity. Then the Lie algebra  ${\Cal G}(G,M)$ of operators, acting on
$V_M$, gives a representation of
the Lie algebra ${\wh {gl_M}}({\Bbb C}_q)$.
 The representation is given by the mapping
$$
E_{ij}\otimes t^m_0t^r_1\mapsto \cases x_{ij}(m, 1, q^r)+{\frac {q^{\frac
{r}2}}{1-q^{r}}}\delta_{ij}\delta_{m0}c
 &\text{if}\; r\not= 0\\
 x_{ij}(m,1,1) &\text{if}\; r=0\endcases
 $$
 $$
 c_0\mapsto c,\;\;\;\;\;\;c_1\mapsto 0\;\;\;\;\;\;\;\;\;\;\;\;\;\;\;
 $$
 for $m,r\in {\Bbb Z}$.\
\endproclaim
 
\medskip
 
\noindent{\bf Remark 4.8}. The representation of the Lie algebra ${\wh
{gl_M}}({\Bbb C}_Q)$
given in (4.7) is called
 the homogeneous realization. This is because of the fact that we are using
the homogeneous gradation. Moreover the
 algebra ${\Cal G}(G,M)$ contains a subalgebra of ${\Cal G}(<1>,M)$ which
is
generated by the operators $
 x_{ij}(m,1,1)$ for $1\le i,j\le M$ and $m\in {\Bbb Z}$, and it is clear
that this subalgebra is nothing but the
 affine algebra ${\wh {gl_M}}({\Bbb C})$ in the homogeneous picture.
 
\medskip
 
 Similarly,  we may have the  principal realization of the Lie algebra
${\wh {gl_N}}({\Bbb
C}_Q)$. For this purpose, we
 choose the group $G=<\xi,q>$ where $q\not= 0$ is not a root of unity and
$\xi$ is a $N$-th primitive
 root of unity. Let $M=1$. Then the Lie algebra ${\Cal G}(G,M)$ is
generated by the coefficients of the
 vertex operators of the form $X_{11}(\xi^{i-1},\xi^{-1}q^r,z)$ for all
$r\in {\Bbb Z}$ and $0\le i\le N-1
 $. From Theorem 3.3 we have
 
\medskip
 
 (i) if $r+s=0$ and ${\overline {i+j}}=0$(mod$N$), then
 $$\align &
 [X_{11}(\xi^{i-1},\xi^{-1}q^r,z_1), X_{11}(\xi^{j-1},\xi^{-1}q^s,z_2)]\\
=&
 X_{11}(\xi^{i+j-1},\xi^{-1}q^{r+s},\xi^{-j}z_1)\delta({\frac
{\xi^jz_2}{q^rz_1}}) \\
 & -
 X_{11}(\xi^{i+j-1},\xi^{-1}q^{r+s},\xi^{-i}z_2)\delta({\frac
{\xi^iz_1}{q^sz_2}})
 +(D\delta)({\frac {\xi^jz_2}{q^rz_1}})c,
\endalign  $$
 
 (ii) if $r+s\not= 0$ or ${\overline {i+j}}\not= 0$(mod$N$), then
 $$\align &
 [X_{11}(\xi^{i-1},\xi^{-1}q^r,z_1), X_{11}(\xi^{j-1},\xi^{-1}q^s,z_2)]\\
=&
 X_{11}(\xi^{i+j-1},\xi^{-1}q^{r+s},\xi^{-j}z_1)\delta({\frac
{\xi^jz_2}{q^rz_1}})
 -
 X_{11}(\xi^{i+j-1},\xi^{-1}q^{r+s},\xi^{-i}z_2)\delta({\frac
{\xi^iz_1}{q^sz_2}})
 \\
& +{\frac {e^{{\frac {i+j} 2}Ln\xi}q^{\frac {r+s}2}}{\xi^{i+j}-q^{r+s}}}
 (\delta({\frac {\xi^jz_2}{q^rz_1}})-\delta({\frac {\xi^iz_1}{q^sz_2}}))c.
 \endalign $$
 Thus if we write
 $$
 {\bar X}_{11}(\xi^{i-1},\xi^{-1}q^r,z)=\cases
X_{11}(\xi^{-1},\xi^{-1},z) &\text{if}\; r=0, i=
 0(\text{mod}N)\\
 X_{11}(\xi^{i-1},\xi^{-1}q^r,z)+
 {\frac {e^{{\frac i 2}Ln\xi}q^{\frac {r}2}}{\xi^{i}-q^{r}}}c 
&\text{otherwise}\endcases
 $$
 for $r\in {\Bbb Z}$ and $0\le i\le N-1$, then
 the above two identities can be written as one identity
 $$
 [{\bar X}_{11}(\xi^{i-1},\xi^{-1}q^r,z_1),
{\bar X}_{11}(\xi^{j-1},\xi^{-1}q^s,z_2)] \tag 4.9
$$
$$
=
 {\bar X}_{11}(\xi^{i+j-1},\xi^{-1}q^{r+s},\xi^{-j}z_1)\delta({\frac
{\xi^jz_2}{q^rz_1}})
  -
 {\bar X}_{11}(\xi^{i+j-1},\xi^{-1}q^{r+s},\xi^{-i}z_2)\delta({\frac
{\xi^iz_1}{q^sz_2}}) 
$$
$$
 +\delta_{r+s,0}\delta_{{\overline {i+j}},0}(D\delta)({\frac
{\xi^jz_2}{q^rz_1}})c.
$$
 Comparing this with the identity (1.20), we get the following result which
was given in [BS] for the
$N=2$ case and in [G2] for arbitrary $N$.
 
 \medskip
 
 \proclaim{Corollary 4.10} Let $M=1$ and $G=<\xi,q>$ be an admissible
subgroup
 of ${\Bbb C}^{\times}$
 generated
by $q$
with
 $q\not= 0$  not a root of unity, and $\xi$ is a $N$-th primitive
 root of unity for $N\ge 2$. Then the Lie algebra ${\Cal G}(G,M)$ of
operators, acting on $V_1$, gives a representation of
 the algebra ${\wh {gl_N}}({\Bbb C}_{q^N})$.
The representation is given by
the mapping
$$
F^iE^{m}\otimes t^m_0t^r_1\mapsto \cases x_{11}(m,\xi^{i-1},\xi^{-1}q^r
)+
\delta_{m,0}{\frac {e^{{\frac i 2}Ln\xi}q^{\frac {r}2}}{\xi^i-q^{r}}}c
 &\text{if}\; i\not= 0(\text{mod}N)\;\;\text{or}\;\;r\not= 0\\
 x_{11}(m,\xi^{-1},\xi^{-1}) &\text{if}\;
i=0(\text{mod}N)\;\;\text{and}\;\;r=0\endcases
 $$
 $$
 c_0\mapsto {\frac cN},\;\;\;\;\;\;c_1\mapsto
0\;\;\;\;\;\;\;\;\;\;\;\;\;\;\;
 $$
 for $r\in {\Bbb Z}$ and $0\le i\le N-1$.\endproclaim
 
 \medskip
 
 \noindent{\bf Remark 4.11}. In general, let $M,N\ge 2$ be integers,  $
G=<\xi,q_1,\cdots, q_{\nu}> 
$
an admissible subgroup of
${\Bbb
C}^{\times}$ with finitely
many generators, where $q_1,\cdots, q_{\nu}$ are the free
generators of $G$ and $\xi$ is
 an $N$-th root of unity.  Then the Lie algebra ${\Cal G}(G,M)$  of
operators, acting on the
 Fock space $V_M$, gives a representation to the Lie algebra ${\wh
{gl_{MN}}}({\Bbb C}_Q)$ where the quantum
 torus ${\Bbb C}_Q={\Bbb C}_Q[t_0^{\pm 1},t_1^{\pm 1},\cdots,t_{\nu}^{\pm
1}]$ is determined by the matrix
 $Q=(q_{ij})_{(\nu+1)\times(\nu+1)}$ with $q_{i0}=q_i^N, q_{0i}=q^{-N}_i$
for $1\le i\le\nu$, and $q_{ij}=1$
 for all other values of $i,j$.
 
 In particular,
if $\nu=1$, that is $G=<\xi,q>$, then the Lie algebra ${\Cal G}(G,M)$
  is generated by the coefficient operators of the vertex operators $
X_{ij}(\xi^{k-1},\xi^{-1}q^r,z)$ for $1\le i,j\le M$, $1\le k\le N-1$ and 
$r\in \Bbb Z$. Moreover, from Theorem 3.3, we have

(i) if ${\overline {k+k'}}\not= 0$ or $r+s\not= 0$, then
$$\align &
[X_{ij}(\xi^{k-1},\xi^{-1}q^r,z_1), X_{i'j'}(\xi^{k'-1},\xi^{-1}
q^{r'},z_2)]\\
=& \delta_{ji'}X_{ij'}(\xi^{k-1},\xi^{-1-k'}q^{r+r'},z_1)\delta({\frac{\xi^{k'}z_2} {q^rz_1}})\\
&-\delta_{j'i}X_{i'j}(\xi^{k'-1},\xi^{-1-k}q^{r+r'},z_2))\delta({\frac
{\xi^{k}z_1} {q^{r'}z_2}})\\
&+{\frac {e^{{\frac {k+k'}{2}}Ln\xi}q^{\frac {r+r'} 2}}{\xi^{k+k'}-q^{r+r'}}}\delta_{ji'}
\delta_{j'i}(\delta({\frac{\xi^{k'}z_2} {q^rz_1}})-\delta({\frac{\xi^{k}z_1} {q^{r'}z_2}}))c\\
=&\delta_{ji'}X_{ij'}(\xi^{k+k'-1},\xi^{-1}q^{r+r'},\xi^{-k'}z_1)\delta({\frac{\xi^{k'}z_2} {q^rz_1}})\\
&-\delta_{j'i}X_{i'j}(\xi^{k+k'-1},\xi^{-1}q^{r+r'},\xi^{-k}z_2))\delta({\frac{\xi^{k}z_1} {q^{r'}z_2}})\\
&+{\frac {e^{{\frac {k+k'}{2}}Ln\xi}q^{\frac {r+r'} 2}}{\xi^{k+k'}-q^{r+r'}}}\delta_{ji'}\delta_{j'i}
(\delta({\frac{\xi^{k'}z_2} {q^rz_1}})-\delta({\frac{\xi^{k}z_1} {q^{r'}z_2}}))c,
\endalign
$$
and

(ii) if ${\overline {k+k'}}= 0$ and $r+s= 0$, then
$$\align &
[X_{ij}(\xi^{k-1},\xi^{-1}q^r,z_1), X_{i'j'}(\xi^{k'-1},\xi^{-1}q^{r'},z_2)]\\
=& \delta_{ji'}X_{ij'}(\xi^{k-1},\xi^{-1-k'}q^{r+r'},z_1)\delta({\frac  
{\xi^{k'}z_2} {q^rz_1}})\\
&-
\delta_{j'i}X_{i'j}(\xi^{k'-1},\xi^{-1-k}q^{r+r'},z_2))\delta({\frac
{\xi^{k}z_1} {q^{r'}z_2}})
+
\delta_{ji'}\delta_{j'i}(D\delta)({\frac
{\xi^{k'}z_2} {q^rz_1}})c\\
=&\delta_{ji'}X_{ij'}(\xi^{-1},\xi^{-1},\xi^{-k'}z_1)\delta({\frac
{\xi^{k'}z_2} {q^rz_1}})\\
&-
\delta_{j'i}X_{i'j}(\xi^{-1},\xi^{-1},\xi^{-k}z_2))\delta({\frac
{\xi^{k}z_1} {q^{r'}z_2}})
+\delta_{ji'}\delta_{j'i}(D\delta)({\frac
{\xi^{k'}z_2} {q^rz_1}})c.
\endalign $$
Set
$$
{\overline X}_{ij}(\xi^{k-1},\xi^{-1}q^r,z)
$$
$$= 
 \cases X_{ij}(\xi^{k-1},\xi^{-1}q^r,z )+
\delta_{ij}{\frac {e^{{\frac k 2}Ln\xi}q^{\frac
{r}2}}{\xi^k-q^{r}}}c
 &\text{if}\; k\not= 0\text{ (mod}N)\;\;\text{or}\;\;r\not= 0\\
 X_{ij}(\xi^{-1}, \xi^{-1},z) &\text{if}\; k=0\text{
(mod}N)\;\;\text{and}\;\;r=0\endcases
 $$
Then we have
$$
[{\overline X}_{ij}(\xi^{k-1},\xi^{-1}q^r,z_1),
{\overline X}_{i'j'}(\xi^{k'-1},\xi^{-1}q^{r'},z_2)]
$$
$$
=\delta_{ji'}{\overline
X}_{ij'}(\xi^{k+k'-1},\xi^{-1}q^{r+r'},\xi^{-k'}z_1)\delta({\frac
{\xi^{k'}z_2} {q^rz_1}})
$$
$$-
\delta_{j'i}{\overline
X}_{i'j}(\xi^{k+k'-1},\xi^{-1}q^{r+r'},\xi^{-k}z_2))\delta({\frac
{\xi^{k}z_1} {q^{r'}z_2}})
+\delta_{ji'}\delta_{j'i}\delta_{\overline
{k+k'},0}\delta_{r+r',0}(D\delta)({\frac
{\xi^{k'}z_2} {q^rz_1}})c.
$$
Comparing this identity with  identity (1.21), we get

\proclaim
{\bf Corollary 4.12}  Let $M,N\ge 2$, and let $G=<\xi,q>$ be an admissible
subgroup
 of ${\Bbb C}^{\times}$
 generated
by $q$
with
 $q\not= 0$  not a root of unity, and $\xi$ an $N$-th primitive
 root of unity. Then the Lie algebra ${\Cal G}(G,M)$ of
operators, acting on $V_1$, gives a representation of
 the algebra ${\wh {gl_{MN}}}({\Bbb C}_{q^N})$, and the representation is
given
by
the mapping
$$
E_{ij}\otimes F^kE^{m}\otimes t^m_0t^r_1
$$
$$
\mapsto \cases
x_{ij}(m,\xi^{k-1},\xi^{-1}q^r
)+
\delta_{m,0}\delta_{ij}
{\frac {e^{{\frac k 2}Ln\xi}q^{\frac
{r}2}}{\xi^i-q^{r}}}c
 &\text{if}\; k\not= 0(\text{mod}N)c\;\;\text{or}\;\;r\not= 0\\
 x_{ij}(m,\xi^{-1},\xi^{-1}) &\text{if}\;
k=0(\text{mod}N)\;\;\text{and}\;\;r=0\endcases
 $$
 $$
 c_0\mapsto {\frac cN},\;\;\;\;\;\;c_1\mapsto
0\;\;\;\;\;\;\;\;\;\;\;\;\;\;\;
 $$
 for $r\in {\Bbb Z}$ and $0\le k\le N-1$.\endproclaim
 
 \medskip

The Lie algebra ${\Cal G}(G,M)$ given in the previous corollary contains 
two interesting 
  subalgebras which give representations to the Lie algebras ${\wh
{gl_{N}}}({\Bbb C}_{q^N})$ and
 ${\wh {gl_{M}}}({\Bbb C}_q)$ . Moreover we will see that
these two subalgebras 
 contain subalgebras that give representations to the affine algebras
${\wh {gl_{N}}}({\Bbb C})$ of level $M$ and
 ${\wh {gl_{M}}}({\Bbb C})$ of level $N$ respectively.
 Indeed, for $a,b\in G=<\xi,q>$, let
 $$
 Y(a,b,z)=\sum_{k=1}^MX_{kk}(a,b,z) \tag 4.13
 $$
 and formally write
$$Y(a,b,z)=\sum_{m\in {\Bbb Z}}y(m,a,b)z^{-m}.\tag 4.14$$
 
 Let ${\Cal L}_1$ be the Lie algebra
 generated by all of the coefficients of $Y(a,b,z)$ for $a,b\in G$.
 We note that $$Y(\xi^iq^r,\xi^jq^s,z)=
 Y(\xi^{i-j-1},\xi^{-1}q^{s-r},\xi^{j+1}q^{r}z),\tag 4.15$$
 so ${\Cal L}_1$ is
indeed generated by the coefficients
 of the vertex operators with the form $Y(\xi^{i-1},\xi^{-1}q^{r},z)$ for
$r\in {\Bbb Z}$ and $0\le i\le
 N-1
 $. Moreover, applying Theorem 3.3, we have, if $r+s\not= 0$ or
${\overline {i+j}}\not= 0$(mod$N$)
 
 $$ 
 [Y(\xi^{i-1},\xi^{-1}q^{r},z_1),Y(\xi^{j-1},\xi^{-1}q^{s},z_2)]\tag
4.16
$$
$$
=
 \sum_{k=1}^M[X_{kk}(\xi^{i-1},\xi^{-1}q^{r},z_1),X_{kk}(\xi^{j-1},
 \xi^{-1}q^{s}
,z_2)] 
 $$
 $$
 =\sum_{k=1}^M\left(X_{kk}(\xi^{i+j-1},\xi^{-1}q^{r+s},\xi^{-j}z_1)
 \delta({\frac 
{\xi^jz_2}{q^rz_1}})-
 X_{kk}(\xi^{i+j-1},\xi^{-1}q^{r+s},\xi^{-j}z_2)\delta({\frac
{\xi^iz_1}{q^sz_2}})\right.
 $$
 $$
 \left.+{\frac {e^{{\frac {i+j} 2}Ln\xi}q^{\frac
{r+s}2}}{\xi^{i+j}-q^{r+s}}}(\delta({\frac {\xi^jz_2}{q^rz_1}})-
 \delta({\frac {\xi^iz_1}{q^sz_2}}))\right)
 $$
 $$
 =Y(\xi^{i+j-1},\xi^{-1}q^{r+s},\xi^{-j}z_1)\delta({\frac
{\xi^jz_2}{q^rz_1}})-
 Y(\xi^{i+j-1},\xi^{-1}q^{r+s},\xi^{-j}z_2)\delta({\frac
{\xi^iz_1}{q^sz_2}})
 $$
 $$
 +{\frac {e^{{\frac {i_j}2}Ln\xi}q^{\frac
{r+s}2}}{\xi^{i+j}-q^{r+s}}}(\delta({\frac {\xi^jz_2}{q^rz_1}})-
 \delta({\frac {\xi^iz_1}{q^sz_2}})),
 $$
 while if $r+s= 0$ and ${\overline {i+j}}= 0$(mod$N$), then
 $$
 [Y(\xi^{i-1},\xi^{-1}q^{r},z_1),Y(\xi^{j-1},\xi^{-1}q^{s},z_2)]
\tag 4.17
$$
$$
=
 \sum_{k=1}^M[X_{kk}(\xi^{i-1},\xi^{-1}q^{r},z_1),X_{kk}(\xi^{j-1},\xi^{-1}
 q^s
,z_2)]
$$
 $$
 =\sum_{k=1}^M\left(X_{kk}(\xi^{-1},\xi^{-1},\xi^{-j}z_1)\delta({\frac
{\xi^jz_2}{q^rz_1}})-
 X_{kk}(\xi^{-1},\xi^{-1},\xi^{-j}z_2)\delta({\frac {\xi^iz_1}{q^sz_2}})
 +(D\delta)({\frac {\xi^jz_2}{q^rz_1}})\right)
 $$
 $$
 =Y(\xi^{-1},\xi^{-1},\xi^{-j}z_1)\delta({\frac {\xi^jz_2}{q^rz_1}})-
 Y(\xi^{-1},\xi^{-1},\xi^{-j}z_2)\delta({\frac {\xi^iz_1}{q^sz_2}})
 +M(D\delta)({\frac {\xi^jz_2}{q^rz_1}}).
 $$
 Therefore, if we define
 $$
 {\bar Y}(\xi^{i-1},\xi^{-1}q^{r},z)=\cases
Y(\xi^{i-1},\xi^{-1}q^{r},z)+M
 {\frac {e^{{\frac i2}Ln\xi}q^{\frac {r}2}}{\xi^{i}-q^{r}}}c
&\text{if}\;r\not= 0\;\text{or}\;{\bar i}\not=
 0\\
 Y(\xi^{-1},\xi^{-1},z) &\text{if}\;r= 0\;\text{and}\;{\bar i}=
 0,\endcases
 $$
 then we can rewrite the  two identities (4.16) and (4.17) into just
one identity
 $$
 [{\bar Y}(\xi^{i-1},\xi^{-1}q^{r},z_1),{\bar
Y}(\xi^{j-1},\xi^{-1}q^{s},z_2)]\tag 4.18
$$
$$
 ={\bar Y}(\xi^{-1},\xi^{-1},\xi^{-j}z_1)\delta({\frac
{\xi^jz_2}{q^rz_1}})
 {\bar Y}(\xi^{-1},\xi^{-1},\xi^{-j}z_2)\delta({\frac
{\xi^iz_1}{q^sz_2}})
 +M(D\delta)({\frac {\xi^jz_2}{q^rz_1}})c.
 $$
 Therefore we have the following result
 
 \medskip
 
 \proclaim{Proposition 4.19}  The Lie algebra ${\Cal L}_1$ of operators
acting on $V_M$ gives a
representation of the
 Lie algebra ${\wh {gl_{N}}}({\Bbb C}_{q^N})$, 
 and the
 representation
 is given by the mapping
 $$
 F^iE^{m}\otimes t^m_0t^r_1\mapsto \cases y(m,\xi^{i-1},\xi^{-1}q^r )+
M{\frac {e^{{\frac i2}Ln\xi}q^{\frac {r}2}}{\xi^i-q^{r}}}\delta_{m,0}c
 &\text{if}\; i\not= 0\text{ (mod}N)\;\;\text{or}\;\;r\not= 0\\
 y(m,1,1) &\text{if}\; i=0\text{ (mod}N)\;\;\text{and}\;\;r=0\endcases
 $$
 $$
 c_0\mapsto Mc,\;\;\;\;\;\;c_1\mapsto 0\;\;\;\;\;\;\;\;\;\;\;\;\;\;\;
 $$
 for $r\in {\Bbb Z}$ and $0\le i\le N-1$.\endproclaim
 
 \medskip
 
 Recall from (4.7), the Fock space $V_M$ affords a representation of the
Lie algebra ${\Cal G}(<q>,M)\subset
 $ ${\Cal G}(<\xi,q>,M)$, where $\xi, q$ are given in (4.10), and
$${\Cal G}(<q>,M)=\text{ span }\{c \text{ and } x_{ij}(m,1,q^r)|\text{
for
} m,r\in {\Bbb Z}, 1\le i,j
 \le M\}. $$
 Now we define a subalgebra of ${\Cal G}(<q>,M)\subset$ ${\Cal
G}(<\xi,q>,M)$
 $$
 {\Cal L}_2=\text{span}\{c\;\; \text{and}\;\; x_{ij}(Nm,1,q^r)|\;
\text{for} \;\;
 m,r\in {\Bbb Z}, 1\le i,j
 \le M\}.
 $$
 Then we have
 
 \medskip
 
 \proclaim{Proposition 4.20} $ {\Cal L}_2$ forms a Lie subalgebra of
${\Cal G}(<q>,M)$, and $ {\Cal
L}_2$
 is also isomorphic to
 ${\Cal G}(<q>,M)$ via the isomorphism given by
 $$
 x_{ij}(m,1,q^r)\mapsto x_{ij}(Nm,1,q^r),\;\;\;\;\;\;c\mapsto Nc.
 $$
 Therefore $ {\Cal L}_2$ gives
 a  representation of ${\wh{gl_{N}}}({\Bbb C}_{q^N})$.
 \endproclaim
 
 \medskip
 
 \proclaim{Proposition 4.21}  For $m,n,r,s,\in {\Bbb Z}$ and $i\not= 0$
(mod $N$), $1\le k\not= l\le
M$, we
have
 $$\align &
 [y(m,\xi^{i-1},\xi^{-1}q^r), x_{kl}(Nn,1,q^s)]
\tag 4.22\\
=& (q^{rNn}-q^{sm})x_{kl}(m+Nn,\xi^{i-1},\xi^{-1}q^{r+s}).
 \endalign
$$
\endproclaim
 \demo{Proof} We apply Theorem 3.3 to obtain
 $$\align &
 [Y(\xi^{i-1},\xi^{-1}q^r, z_1), X_{kl}(1,q^s,z_2)]=
 [\sum_{j=1}^MX_{jj}(\xi^{i-1},\xi^{-1}q^r, z_1), X_{kl}(1,q^s,z_2)]\\
=&\sum_{j=1}^M\{X_{jl}(\xi^{i-1},\xi^{-1}q^{r+s},z_1)\delta_{jk}\delta({\frac
{z_2}{\xi^{-1}q^rz_1}})
-X_{kj}(1,\xi^{-i}q^{r+s},z_2)\delta_{jl}\delta({\frac
{\xi^{i-1}z_1}{q^sz_2}})\}\\
=& X_{kl}(\xi^{i-1},\xi^{-1}q^{r+s},z_1)\delta({\frac
{z_2}{\xi^{-1}q^rz_1}})
-X_{kl}(1,\xi^{-i}q^{r+s},z_2)\delta({\frac {\xi^{i-1}z_1}{q^sz_2}})\\
=& X_{kl}(\xi^{i-1},\xi^{-1}q^{r+s},z_1)\delta({\frac
{z_2}{\xi^{-1}q^rz_1}})
- X_{kl}(\xi^{i-1},\xi^{-1}q^{r+s},q^{-s}z_1) \delta({\frac
{\xi^{i-1}z_1}{q^sz_2}}).
\endalign
$$
This then gives
$$
[y(m,\xi^{i-1},\xi^{-1}q^r),
x_{kl}(n,1,q^s)]=\xi^{-n}(q^{rn}-q^{sm}\xi^{in})
x_{kl}(m+n,\xi^{i-1},\xi^{-1}q^{r+s}).
$$ which immediately implies (4.22).\qed
 
\enddemo
 
 \medskip
 
\noindent{\bf Remark 4.23} Let ${\Cal G}_i\subset {\Cal L}_i\subset {\Cal
G}(<\xi,q>,M)$, $i=1,2$, be such
that
 $$
 {\Cal
G}_1=\text{span}\{c\;\; \text{and}\;\; y(m,\xi^{i-1},\xi^{-1})|\; \text{
for }
 m\in {\Bbb Z}, 0\le i
 \le N-1\}, 
 $$
 $$
 {\Cal G}_2=\text{span}\{c\;\; \text{and}\;\; x_{ij}(Nm,1,1)|\; \text{
for }
\;\;
 m\in {\Bbb Z}, 1\le i,j
 \le M\}.
 $$
  Then the two subalgebras ${\Cal G}_1, {\Cal G}_2$ of ${\Cal
G}(<\xi,q>,M)$
  respectively give representations of the
 affine algebra ${\wh{gl_N}}({\Bbb C})$ of level $M$ and
${\wh{gl_M}}({\Bbb C})$ of
level $N$. Let
 ${\Cal G}_i'$ be the derived algebras of ${\Cal G}_i$. Then we have the
so
called dual pair property given in [F1]: 
 $
 [{\Cal G}_1', {\Cal G}_2']=(0).    
$
 However, clearly, we have  $[{\Cal L}_1', {\Cal L}_2']\not= (0)$.
                                                              
\bigskip
 
\noindent{\bf Acknowledgements }
 
\medskip
 
SB and YG are  supported by grants from  the
Natural Sciences and Engineering Research Council of Canada.
ST is  supported by a grant from the National
Natural Science Foundation of
China. The authors would like to thank the Fields Institute  for its
hospitality
during the preparation of this work.

\enddocument
  \loadeufm
\loadeusm
\loadbold
\input amstex
\documentstyle{amsppt}
\NoRunningHeads
 
\magnification=1200
\pagewidth{32pc}
\pageheight{42pc}
\vcorrection{1.2pc}

\define\wh{\widehat}

\define\bc{\Bbb C}

\define\bbz{{\Bbb Z}^{\nu+1}}
\define\bbl{{\hat{\Cal L}_{{\Bbb C}_Q}}}

\define\vep{\varepsilon}

\topmatter
\title A Unified View of some Vertex Operator Constructions
\endtitle
\author Stephen Berman\footnotemark"$^{1}$", Yun
Gao\footnotemark"$^{2}$", Shaobin Tan
\footnotemark"$^{3}$"
\endauthor
\footnotetext"$^1$"{ Department of Mathematics and Statistics, University
of Saskatchewan,
Saskatoon, Canada
S7N 5E6; berman\@math.usask.ca}
\footnotetext "$^2$"{Department of Mathematics and Statistics,
 York University, Toronto, Canada M3J 1P3; ygao\@ yorku.ca}
\footnotetext"$^3$"{ Department of Mathematics, Xiamen University, Xiamen,
China 361005; tans\@jingxian.xmu.edu.cn} 

 
\abstract
We present a general vertex operator construction based on the Fock space 
for an affine Lie algebras of type $A$.
This construction allows us to give a unified treatment for both the
homogeneous and principle
realizations of the affine Lie algebras $\hat{gl}_N$ as well as  for some extended affine Lie algebras 
coordinatized by certain quantum tori.
 
\endabstract
 
\endtopmatter \document
 
\subhead\S 0. Introduction\endsubhead
 
\medskip
 This paper presents a unified view of certain vertex operator constructions for some of the extended affine Lie 
algebras 
(EALA's for short) which are coordinatized by certain quantum tori. Recall that for the affine Kac-Moody Lie algebras 
vertex operator representations were developed in [LW] and [KKLW] for the principal realizations and in [FK],[S] 
in the homogeneous realization. Our motivation comes from the paper [F1] of I. Frenkel, where he presented a unified 
construction for both the principal and homogeneous realizations of the affine Lie algebras of type $A^{(1)}.$ 
This is accomplished by using the affine algebra $\hat{gl}_M$ rather than $\hat {sl}_M.$ Moreover Frenkel used a 
Clifford  algebra structure, which was inherent in his situation, to define a new type of normal ordering which then 
led to his unified view. The Clifford structure had been studied before in the works [F3,4] and [KP].

   The structure theory of EALA's has been developed over the last ten years(see [H-KT], [BGK],[AABGP] and [ABGP]). 
Roughly speaking these Lie algebras are generalizations of both the affine Kac-Moody Lie algebras and the finite 
dimensional simple Lie algebras over the complex numbers which admit Laurent like coordinates  in a finite number 
of variables. It turns out that algebras of different types admit different types of coordinates. For example, 
those of type $A_l$ admit the non-commutative quantum torus as coordinates (see [M],[BGK]).
Representations for these Lie algebras over quantum tori have been constructed in
[JK], [G-KK], [G1,2,3], [BS], [VV].  When $l=1$ there are even Jordan 
algebras which serve as coordinates of EALAs.  Representations for this kind of
Lie algebra have been initiated by [T1]. 
Perhaps the  examples which have attracted the most attention so far are the toroidal algebras 
which have the commutative associative Laurent polynomials as their coordinates. The toroidal 
algebras have been studied since the mid 80's both in terms of their structure theory as well as their representation 
theory (see [F2], [MRY], [Y], [W], [BC], [FJW],[T2], and among others). For our purposes we 
want to mention that vertex operator representations have played a predominant role in much of this work. Indeed, in 
[G1,2], one finds both homogeneous and principal realizations given for many of the EALA's  of type $A$. The principal
realization for those EALAs was also implicitly given by [G-KK] in studying the so-called $\Gamma$-conformal algebras.  
Our goal 
in this work is to unify the various approaches and show how they all follow from the same type of approach. Of course, 
the work in the affine case, namely [F1],  shed light on doing this.

   Working with a standard type of Fock space we are able to define a general type of vertex operator which depends on 
a non-zero scalar from $\bc$ and to then  compute the commutator of two of these operators. This is presented in the 
second section of this paper while, in the first section, we give the basics on the Lie algebras, which are all of 
type $A$, which we will later go on to find representations for. Already in Section one it is evident that there is 
somewhat of a unified picture for these algebras. When we define our vertex operators in Section two the reader will 
see that we are using a Clifford algebra structure to define the normal ordering we are using, just as was done in [F1]. 
In the third section we introduce some Lie algebras associated to certain choices of subgroups, $G$, of non-zero 
complex numbers as well as the choice of a positive integer $M$. These Lie algebras are spanned by the moments of our 
vertex operators and hence, by  construction, we automatically have a representation for this Lie algebra. We show the 
representations we have are completely reducible and find the irreducible components. In the fourth and final section 
we show how certain choices of the group $G$ and the integer $M$ lead to representations of the algebras of section one. 
Thus, we recover both the principal and homogeneous representations for the affine algebras of type $A$ as well as those 
for the EALA's studied in [BS], [G1,2], [G-KK]. It is from seeing the various applications in this fourth section that one 
understands the unification of our treatment. Finally we want to emphasize that this unified treatment would not have 
been possible without first knowing the particular special cases of this result.

\bigskip
 
\subhead \S 1. Preliminaries\endsubhead
 
\medskip
 
In this section we shall review some of the  basics on Lie algebras coordinatized
by quantum tori. We present this from a general point of view which unifies our treatment. For notation we always denote the integers, positive
 and negative integers respectively by ${\Bbb Z}$,
${\Bbb Z}_+$, and ${\Bbb Z}_-$.
 
Let $\frak g$ be any associative $\Bbb C$-algebra with a symmetric
bilinear form $(\cdot,\cdot)$: $\frak g$$\times \frak g$$\to {\Bbb C}$
such that $(xy,z)=(x,yz)$ for $x,y,z\in$ $\frak g$. Let
$A=\oplus_{\alpha\in {\Bbb Z}^{\nu+1}} A_{\alpha}$, $\nu\ge 0$, be any
${\Bbb Z}^{\nu+1}$-graded associative algebra such that
dim$A_{\alpha}<\infty$ for all $\alpha\in {\Bbb Z}^{\nu+1}$. Fix a base
$(x_{i
\alpha})_{i\in I_{\alpha}}$ of $A_{\alpha}$, where $I_{\alpha}$ is the
index set corresponding to the subspace $A_{\alpha}$. Let 
 $d_0,d_1,\cdots,d_{
\nu}$ be degree derivations of $A$ such that $d_ix=\alpha_ix$ for $x\in
A_{\alpha}$, $i=0,1,\cdots,\nu$ and $\alpha=(\alpha_0,
\cdots, \alpha_{\nu})\in {\Bbb Z}^{\nu+1}.$ We define a $\Bbb C$-linear
map $\phi:$ $A\to {\Bbb C}$ by linear extension of
$$
  \phi(x_{i\alpha})=\cases 1 &\text{if}\;\;{\alpha=(0,\cdots, 0)}\\
0 &\text {otherwise}.\endcases\tag 1.1
$$
for $i\in I_{\alpha},$ $\alpha\in {\Bbb Z}^{\nu+1}$. The tensor product
${\frak g}\otimes_{\Bbb C}A$, with the canonical product $(x\otimes a)
(y\otimes b)=xy\otimes ab$ for $x,y\in {\frak g}, a,b\in A$, is also an
associative algebra. Moreover with the commutator product
$[x\otimes a, y\otimes b]_{loop}=(x\otimes a)
(y\otimes b)-
(y\otimes b)(x\otimes a),$ ${\frak g}\otimes_{\Bbb C}A$ forms a ${\Bbb
Z}^{\nu+1}$-graded Lie algebra. We call this algebra a loop type
Lie algebra. Consider the vector space
$$
{\hat {\frak g}}_A:=({\frak g}\otimes_{\Bbb C}A)\oplus {\Cal C}\tag 1.2
$$
where $\Cal C$=$\oplus_{0\le i\le \nu}{\Bbb C}c_i$ is a
$\nu+1$-dimensional
vector space. There is an alternating bilinear map $[\cdot,\cdot]
:$ ${\hat {\frak g}}_A\times {\hat {\frak g}}_A\to {\hat {\frak g}}_A$
determined by the conditions:
$$
[c_i,{\hat {\frak g}}_A]=0
$$
$$
[x\otimes a, y\otimes b]=[x\otimes a, y\otimes b]_{loop}+(x,y)\sum_{0\le
i\le\nu}\phi((d_ia)b)c_i
$$
for $x,y\in {\frak g}, a,b\in A$ and $i=0,1,\cdots, \nu$. It is
straightforward to check that ${\hat {\frak g}}_A$ is a Lie algebra.
Indeed
there is an exact sequence of Lie algebras with canonical maps
$$
0\to \oplus_{0\le i\le\nu}{\Bbb C}c_i\to {\hat {\frak g}}_A\to 
{\frak g}\otimes_{\Bbb C}A\to 0, 
$$
and so we have that ${\hat {\frak g}}_A$ is a central extension of the loop type Lie
algebra ${\frak g}\otimes_{\Bbb C}A$.
 
  Let $M_{\infty}(\bc)=$span$_{\bc}\{E_{ij}|\;\;1\le i,j <\infty\}$ be
the infinite matrix algebra, where $E_{ij}$ is the infinite matrix with a $1$
in the $(i,j)$-entry and zero's elsewhere. We also let 
$M_n(\bc)=$span$_{\bc}\{E_{ij}|\;\;1\le i,j\le n\}$. This subspace
of $M_{\infty}(\bc)$ for $n\ge 1$ is isomorphic to the usual matrix algebra of $n \times n$ 
matrices with entries in $\bc$ .

Let $Q=(q_{ij})$ be a $(\nu+1)\times(\nu+1)$ matrix with entries
$q_{ij}\in
{\Bbb C}^{\times}$ satisfying $q_{ii}=1$ and $q_{ij}=q_{ji}^{-1}$ for
$0\le i,j\le \nu$. The quantum torus associated with the matrix
$Q$ is a unital associative ${\Bbb C}$-algebra
${\Bbb C}_Q :={\Bbb C}_Q[t_0^{\pm 1},\cdots,t_{\nu}^{\pm 1}]$
with generators $t_0^{\pm 1},\cdots,t_{\nu}^{\pm 1}$ and relations
 $
t_it_i^{-1}=t_i^{-1}t_i=1, \quad  t_it_j=q_{ij}t_jt_i,
$
for $0\le i,j\le \nu$. If $Q$ is $2\times 2$ matrix, so then $\nu=1$, 
the matrix $Q=(q_{ij})$ is determined by a single $q=q_{10}$. In this case we often simply
denote
${\Bbb C}_Q={\Bbb C}_Q[t_0^{\pm 1},t_1^{\pm 1}]$ by ${\Bbb C}_{q}$.   
Choose the bilinear
form on $M_n(\Bbb C)$ to be the trace form. Set $A=\bc_Q[t^{\pm
1}_0,\cdots,
t^{\pm 1}_{\nu}]$,
 with the ${\Bbb Z}^{\nu+1}$-gradation $A=\oplus_{\alpha\in
\bbz}A_{\alpha}$, where
the subspace $A_{\alpha}$ is spanned by
$t^{\alpha}=t_0^{\alpha_0}t_1^{\alpha_1}\cdots t_{\nu}^{\alpha_{\nu}}$
for $\alpha=
(\alpha_0,\cdots, \alpha_{\nu})\in\bbz$.      Define $\sigma_Q:$
$\bbz\times\bbz\to \bc$ by
$$
\sigma_Q(\alpha,\beta)=\prod_{0\le i<j\le\nu}q_{ji}^{\alpha_j\beta_i}\tag
1.3
$$
for $\alpha=(\alpha_0,\cdots,\alpha_{\nu})$,
$\beta=(\beta_0,\cdots,\beta_{\nu})\in\bbz$. Then we have
 $
t^{\alpha}t^{\beta}=\sigma_Q(\alpha,\beta)t^{\alpha+\beta}.$

The proof of the following Lemma is clear.

\medskip

\proclaim
{\bf Lemma 1.4} Let $m,n\ge 1$ be integers. Then there is a Lie algebra
isomorphism
$$
(M_m(\bc)\otimes M_n(\bc))^{\wedge}_{\bc_Q}\cong
(M_{mn}(\bc))^{\wedge}_{\bc_Q}
$$
which is given by
$$
E_{ij}\otimes E_{kl}\otimes t^{\alpha}\mapsto E_{(i-1)n+k,(j-1)n+l}\otimes
t^{\alpha}
$$
$$
c_s\mapsto c_s,\;\;\;\; s=0,1,\cdots, \nu
$$
for $\alpha=(\alpha_0,\cdots,\alpha_{\nu})\in\bbz$, $1\le i,j\le m, 1\le
k,l\le n.$
\endproclaim
 
 Let $\bbl$ be the Lie subalgebra of $(M_m(\bc)\otimes
M_n(\bc))^{\wedge}_{\bc_Q} $ generated by elements of the form
$ E_{ij}\otimes E_{kl}\otimes t_0^{\alpha_0(n-1)+l-k}t^{\alpha}$ for $1\le
i,j\le m$, $1\le k,l\le n$ and $\alpha=
(\alpha_0,\cdots,\alpha_{\nu})\in \bbz$. The following result gives the structure of $\bbl .$

\proclaim 
{\bf Proposition 1.5} 
$$
 \bbl\cong(M_m(\bc)\otimes
M_n(\bc))^{\wedge}_{\bc_{Q^{\ast}}}  
$$
 where ${\bc_Q}={\bc_Q}[t_0^{\pm
1},\cdots, t_{\nu}^{\pm 1}]
$ with $Q=(q_{ij})$, and ${\bc_{Q^{\ast}}}={\bc_{Q^{\ast}}}[\tau_0^{\pm
1},\cdots,
\tau_{\nu}^{\pm
1}]
$ with $Q^{\ast}=(q^{\ast}_{ij})$ such that $q^{\ast}_{ij}=q_{ij}$ if
$i,j\not= 0$, and $q^{\ast}_{ij}=q^n_{ij}$ if $i=0$ or $j=0$.
 
\endproclaim
Proof. Define a linear map $f:$ $(M_m(\bc)\otimes
M_n(\bc))^{\wedge}_{\bc_{Q^{\ast}}}  \to \bbl$ by
$$
E_{ij}\otimes E_{kl}\otimes \tau^{\alpha}\mapsto (\prod_{1\le
s\le\nu}q_{s0}^{l\alpha_s})E_{ij}\otimes
E_{kl}\otimes
t_0^{(n-1)\alpha_0+l-k}t^{\alpha}-k\delta_{ij}\delta_{kl}\delta_{\alpha,{
0}}c_0
 $$
$$
c_0\mapsto nc_0,\;\;\;\;\;\; c_s\mapsto c_s, \;\;\;s=1,2,\cdots, \nu
$$
for $\alpha=(\alpha_0,\cdots, \alpha_{\nu})\in\bbz$, $1\le i,j \le m$ and
$1\le k,l\le n$. Let ${\bar
{\alpha}}=(n\alpha_0+l-k,\alpha_1,\cdots,\alpha_{\nu})$, and 
${\bar {\alpha}'}=(n\alpha_0'+l'-k',\alpha_1',\cdots,\alpha_{\nu}')\in
\bbz$. Using the identity
$$
\sigma_{Q}({\bar {\alpha}},{\bar
{\alpha}'})=\sigma_{Q^{\ast}}(\alpha,\alpha')\prod_{1\le j\le
\nu}q_{j0}^{\alpha_j(l'-k')}\tag 1.6
$$
one can easily check that the map $f$ is the desired Lie algebra
isomorphism.  {\qed}

Putting together the two previous results we get the following identification of $\bbl$.

\medskip
\proclaim
 {\bf Corollary 1.7} 
$$
\bbl\cong(M_{mn}(\bc))^{\wedge}_{\bc_{Q^{\ast}}}  
$$
where $Q$ and $Q^{\ast}$ are given in the previous proposition .
{\qed}
\endproclaim
\medskip
 
Let $\xi=\xi_n$ be an $n$-th primitive root of unity and let $E,F\in
M_n(\bc)$ be defined by saying
$$
E=E_{12}+E_{23}+\cdots +E_{n-1,n}+E_{n1}, \quad
F=\sum_{i=1}^nE_{ii}(\xi^{i-1}).
\tag 1.8
$$
Then the following fact is well-known.
 
\proclaim{Lemma 1.9} The set of matrices $\{F^iE^j\}_{1 \leq i, j\leq
n}$, forms a basis of
the 
matrix algebra $M_n(\bc)$ (so a basis of the general linear Lie algebra
${gl}_n(\Bbb{C})$). Moreover,
$$EF=\xi FE, \quad E^n=F^n=Id_n, \tag 1.10$$
and
$$
E_{ij}={\frac 
1 n}\sum^{n-1}_{k=0}\xi^{k(1-i)}F^kE^{j-i},\;\;\;\;F^iE^j=\sum^n_{l=1}
\xi^{i(
l-1)}E_{l,{\overline {l+j}}}\tag 1.11
$$
for $1\le i,j\le n$, where, for notation, we are letting ${\bar l}$ denoted the unique integer, 
$l$, in $\{ 1, 2, \cdots, n \}$ such that ${\bar l}=l($mod $n)$.
\endproclaim
 
\medskip

Note that
$$
E_{ij}\otimes E_{kl}\otimes t_0^{\alpha_0(n-1)+l-k}t^{\alpha}=
\sum_{s=0}^{n-1}\xi^{s(1-k)}E_{ij}\otimes F^sE^{l-k}\otimes
t_0^{\alpha_0n+l-k}t_1^{\alpha_1}\cdots t_{\nu}^{\alpha_{\nu}}
$$
$$
=\sum_{s=0}^{n-1}\xi^{s(1-k)}E_{ij}\otimes F^sE^{\alpha_0'}\otimes
t_0^{\alpha_0'}t_1^{\alpha_1}\cdots t_{\nu}^{\alpha_{\nu}}
$$
where $\alpha_0'=\alpha_0n+l-k$, for $1\le i,j\le m$, $1\le k,l\le n$, and
$\alpha=(\alpha_0,\alpha_1,\cdots, \alpha_{\nu})\in
\bbz.$ From this one sees that the Lie subalgebra $\bbl$ of
$(M_m(\bc)\otimes M_n(\bc))^{\wedge}_{\bc_Q}  $ has a basis
of the form
$$
E_{ij}\otimes F^kE^{l_0}\otimes t_0^{l_0}\cdots
t_{\nu}^{l_{\nu}},\;\;\;\; c_0,\; c_1,\cdots, c_{\nu}\tag 1.12
$$
where $1\le i,j\le m$, $0\le k\le n-1$ and $l_0,l_1,\cdots,l_{\nu}\in
{\Bbb Z}$. Moreover, the commutation relation of $\bbl$ are
determined by
$$
[E_{ij}\otimes F^kE^{\alpha_0}\otimes t^{\alpha},
E_{i'j'}\otimes F^{k'}E^{\alpha_0'}\otimes t^{\alpha'}]\tag 1.13
$$
$$
=\delta_{ji'}\xi^{\alpha_0k'}\sigma_{Q}(\alpha,\alpha')E_{ij'}\otimes 
F^{k+k'}E^{
\alpha_0+\alpha_0'}\otimes t^{\alpha+\alpha'}
$$
$$
-\delta_{j'i}\xi^{\alpha_0'k}\sigma_{Q}(\alpha',\alpha)E_{i'j}\otimes 
F^{k+k'}E^{
\alpha_0+\alpha_0'}\otimes t^{\alpha+\alpha'}
+n\delta_{ji'}\delta_{ij'}\delta_{{\overline {k+k'}},{
0}}\delta_{\alpha+\alpha', 0}\xi^{\alpha_0k'}\sum_{0\le
s\le\nu}\alpha_sc_s
$$
for $1\le i,i',,j,,j'\le m,$ $0\le k,k'\le n-1$,
$\alpha=(\alpha_0,\cdots,\alpha_{\nu}),$
$\alpha'=(\alpha_0',\cdots,\alpha_{\nu}')
\in \bbz,$ as well as the fact that the elements $c_0,\cdots, c_{\nu}$ are central in $\bbl$.
 
 From now on we will identity the Lie algebra $(M_m(\bc)\otimes
M_n(\bc))^{\wedge}_{\bc_Q}$ with $(M_{mn}(\bc))^{\wedge}_{\bc_Q}$,
and also identify $\bbl$ with   $(M_{mn}(\bc))^{\wedge}_{\bc_{Q^{\ast}}}$, where
$Q=(q_{ij}),$ $Q^{\ast}=(q^{\ast}_{ij})$ and as above $q^{\ast}_{ij}=q_{ij}
$ if $i,j\not= 0$, and $q^{\ast}_{ij}=q^n_{ij}$ if $i=0$ or $j=0$. For 
simplicity we will write $a\alpha=(a\alpha_1,\cdots, a\alpha_{\nu})$ for 
$a\in{\Bbb Z}$ and $\alpha\in {\Bbb Z}^{\nu}$, also we will write
$$
q_0^{\alpha}=q_{10}^{\alpha_1}\cdots q_{\nu 0}^{\alpha_{\nu}}
$$
for $q_0=(q_{10},\cdots, q_{\nu 0})\in {\Bbb C}^{\nu}$.

The
Lie algebra structure (1.13) of $\bbl$ can be described by formal
power series identities. For this purpose we let $z,z_1,z_2$ be formal
variables. For $1\le i,j\le m, $  $0\le k\le n-1$ and
$\alpha=(\alpha_1, \cdots, \alpha_{\nu})\in {\Bbb Z}^{\nu}$, we set
$$
X_{ij}^k(\alpha,z)=\sum_{l\in \Bbb Z}(E_{ij}\otimes F^k E^l\otimes
t_0^lt_1^{\alpha_1}\cdots t_{\nu}^{\alpha_{\nu}})z^{-l}\in {\bbl}[[
z,z^{-1}]],\tag 1.14
$$
and $\delta(z)=\sum_{l\in {\Bbb Z}}z^l$, $(D\delta)(z)=\sum_{l\in \Bbb
Z}lz^l$. Then the algebra structure of $\bbl$ is described by the
following lemma.
 
\proclaim
{\bf Lemma 1.15} Let $1\le i,j,i',j'\le m$, $0\le k,k'\le n-1$,
$\alpha=(\alpha_1,\cdots, \alpha_{\nu})$, $\alpha'
=(\alpha_1',\cdots, \alpha_{\nu}')\in {\Bbb Z}^{\nu}$. Then the following
power series identity  is equivalent to (1.13) 
$$
[X_{ij}^k(\alpha, z_1), X_{i'j'}^{k'}(\alpha',
z_2)]=\delta_{ji'}\sigma(\alpha,\alpha')X_{ij'}^{\overline
{k+k'}}(\alpha+\alpha',\xi^{-k'}z_1)\delta({\frac
{\xi^{k'}z_2} {z_1q^{\alpha}_0}})\tag 1.16
$$
$$
-\delta_{j'i}\sigma(\alpha',\alpha)X_{i'j}^{\overline
{k+k'}}(\alpha+\alpha',\xi^{-k}z_1)\delta({\frac
{z_2q_0^{\alpha'}}{\xi^{k}z_1} })
$$
$$
+n\delta_{ji'}\delta_{ij'}\delta_{\overline 
{k+k'},0}\delta_{\alpha+\alpha',
{ 0}}\sigma(\alpha,\alpha')\{(D\delta)({\frac
{\xi^{k'}z_2} {z_1q_0^{\alpha}}})c_0+\delta({\frac
{\xi^{k'}z_2} {z_1q_0^{\alpha}}})\sum_{1\le s\le
\nu}\alpha_sc_s\}
$$
where ${\bar k}=k($mod$ \ n)$ and $ k \in\{0,1,\cdots, n-1\}$.
 \endproclaim

As  very special cases, one  chooses $n=1,$ $\nu=0$, then $\bbl$ is just
the affine algebra ${\widehat {gl}_m}(\bc)$ in the so-called homogeneous
picture; while if one chooses $m=1,$ $\nu=0$ then $\bbl$ is 
the affine algebra ${\widehat {gl}_n}(\bc)$ in the so-called
principal picture. In these two cases, the identity (1.16) can simply be written
as
follows 
$$
[X_{ij}^0(z_1), X_{kl}^0(z_2)]=X_{il}^0(z_1)\delta_{jk}\delta({\frac
{z_2} {z_1}})-
X_{kj}^0(z_2)\delta_{il}\delta({\frac {z_2}
{z_1}})+\delta_{jk}\delta_{il}(D\delta)({\frac {z_2} {z_1}})c_0,
\tag 1.17
$$
for $1\le i,j,k,l\le m$, for the first of these  and
$$
[X^i_{11}(z_1), X^j_{11}(z_2)]=X^{\overline {i+j}}_{11}(z_2)\delta({\frac
{\xi^jz_2} {z_1}})-
X^{\overline {i+j}}_{11}(z_1)\delta({\frac {\xi^iz_1}
{z_2}})+n\delta_{{\overline{i+j}},0}(D\delta)({\frac {\xi^jz_2} {z_1}}
)c_0.\tag 1.18
$$
for $0\le i,j\le n-1$,  ${\bar i}=i($mod$ \ n)$ and $ i \in \{0,1,\cdots, n-1\}$ for the second one.
 
  Moreover if we choose $n=1, \nu=1$, or $m=1, \nu=1$, and write
$q_{10}=q$, then $\bbl$ gives respectively the homogeneous
realization of the Lie algebra ${\widehat {gl}_m}(\bc_q)$, and the
principal
realization of ${\widehat {gl}_n}(\bc_{q^{n}})$. The
algebra structure of these two cases can
be described as follows
$$
  [X_{ij}^0(r,z_1), X_{kl}^0(s,z_2)]  
=  X_{il}^0(r+s,z_1)\delta_{jk}\delta({\frac
{z_2}{q^rz_1}})-
 X_{kj}^0(r+s,z_2)\delta_{il}\delta({\frac {q^sz_2}{z_1}}) \tag 1.19
$$
$$
 +\delta_{il}\delta_{jk}\delta_{r+s,0}((D\delta)({\frac
{z_2}{q^rz_1}})c_0+
 r\delta({\frac {z_2}{q^rz_1}})c_1) $$
 for $1\le i,j,k,l\le m$ and $r,s\in {\Bbb Z}$, for the first and
 $$
  [X^i_{11}(r,z_1), X^j_{11}(s,z_2)]= X^{\overline
{i+j}}_{11}(r+s,\xi^{-j}z_1)\delta({\frac
{\xi^jz_2}{q^rz_1}})-
 X^{\overline {i+j}}_{11}(r+s,\xi^{-i}z_2)\delta({\frac
{q^sz_2}{\xi^iz_1}})
\tag 1.20
 $$
$$
 +n\delta_{r+s,0}\delta_{\overline{i+j},0}((D\delta)({\frac
{\xi^jz_2}{q^rz_1}})c_0+r
 \delta({\frac {\xi^jz_2}{q^rz_1}})c_1)
 $$
 for $1\le i,j\le n$, $r,s\in {\Bbb Z}$, and ${\overline{i+j}}=i+j$(mod
$n$) for the second.
 
Finally if we choose $\nu=1$,  $m,n\ge 1$, write $q_{10}=q$, then
$\bbl$ is isomorphic to the affine Lie algebra
${\widehat {gl}_{mn}}(\bc_{q^n})$, which contains the special cases
mentioned above. The algebra structure is as follows.
$$
[X^k_{ij}(r,z_1),X^{k'}_{i'j'}(r',z_2)]=\delta_{ji'}X_{ij'}^{\overline{k+k'}}
(r+r',\xi^{-k'}z_1)\delta({\frac {\xi^{k'}z_2} {q^rz_1}})\tag 1.21
$$
$$
-\delta_{j'i}X_{i'j}^{\overline{k+k'}}
(r+r',\xi^{-k}z_2)\delta({\frac {q^{r'}z_2} {\xi^{k}z_1}})
$$
$$
+n\delta_{ji'}\delta_{ij'}\delta_{{\overline {k+k'}},0}\delta_{r+r',0}\{
(D\delta)({\frac {\xi^{k'}z_2} {q^rz_1}})c_0+r\delta({\frac {\xi^{k'}z_2}
{q^rz_1}})c_1\}.
$$
In subsequent sections we are going to give
 irreducible representations for a class of Lie algebras which include the
Lie algebras mentioned  above.
 
\medskip
\medskip
 
\subhead \S 2.  Fock Space and Vertex Operators \endsubhead
 
\medskip
 
In this section, we shall define the Fock space we need and  construct a family
of vertex
operators acting on it. Then we go on to derive the commutation relations between these vertex operators in various
situations.
Some of these commutation relations were implicitly worked out in [G1].
However, we will use the
ideas from [F1] to tie a Clifford algebra structure to our vertex operators. This makes our approach
very natural and concise.
 
Let $\vep_1,\dots,\vep_M$ $(M\ge 1)$ be symbols. We form lattices
$$
\Gamma_M=\oplus^M_{i=1}{\Bbb Z}\vep_i,\;\;\;\;Q_M=\oplus^{M-1}_{i=1}{\Bbb
Z}(\vep_i-\vep_{i+1}),
\tag 2.1 $$
with a symmetric bilinear form $(\vep_i,\vep_j)=\delta_{ij}$. We also
extend this bilinear form to the ${\Bbb
C}$-vector space
$$ H_M:={\Bbb C}\otimes \Gamma_M.\tag 2.2 $$
 
 For each $k \in \Bbb Z$ we take a copy of  $H_M$ with basis labeled by  $\vep_i(k)$  for $1\leq i \leq M,
k\in \Bbb Z$. That is, $\vep_i(k)$, is to be a  copy of $\vep_i$. We form a Lie algebra
$$
{\Cal H}_M={\hbox {span}}_{\Bbb C}\{\vep_i(k),\; c|1\le i\le
M,\;k\in{\Bbb Z}\},
\tag 2.3 $$
with the Lie product
$$[\alpha(k), \beta(l)]=k(\alpha,\beta)\delta_{k+l,0}c, \tag 2.4 $$
for $\alpha,\beta\in H_M,$
$k,l
\in {\Bbb Z}$, where $c$ is a central element. Let
$${\Cal H}^{\pm}_M={\hbox
{span}}\{\vep_i(k)|k\in{\Bbb Z}_{\pm},
1\le i\le M\}. \tag 2.5$$
 Then $${\hat {\Cal H}}_M={\Cal H}^{+}_M+{\Bbb C}c+{\Cal H}^{-}_M \tag
2.6$$
forms a Heisenberg
subalgebra of ${\Cal H}_M$.
 
Let ${\Cal S}({\Cal H}^{-}_M)$ be the symmetric algebra over the abelian
algebra ${\Cal H}^-_M$ and lat
$${\Bbb C}[\Gamma_M]:= \oplus_{\alpha\in \Gamma_M}{\Bbb
C}e^{\alpha}, \tag 2.7$$
 be the group algebra over $\Gamma_M$ twisted by a 2-cocycle so that 
$e^{\alpha}e^{\beta}=\epsilon(\alpha,
 \beta)e^{\alpha+\beta}$ for $\alpha,\beta\in\Gamma_M$. The cocycle 
$$\epsilon:\; \Gamma_M\times \Gamma_M\to \{\pm 1\}, \tag 2.8$$ is defined
by setting
$$\epsilon(\vep_i,\vep_j)=1  \text{ if }
i\le j,  \quad \epsilon(\vep_i,\vep_j)=-1 \text{ if }
i>j, \tag 2.9$$
 and $$\epsilon(\sum_im_i\vep_i,
\sum_jn_j\vep_j)=\prod_{i,j}(\epsilon(\vep_i,\vep_j))^{m_in_j}, \tag
2.10$$
for $m_i,n_j\in{\Bbb Z}$. One can easily check the following result.
 
\proclaim{Lemma 2.11} $\epsilon$ is bi-multiplicative on $\Gamma_M$.
Moreover.
$$\epsilon(\alpha,\alpha)=(-1)^{\frac{(\alpha,\alpha)}{2}},\;\;
\epsilon(\alpha,\beta)\epsilon(\beta,\alpha)
=(-1)^{(\alpha,\beta)}$$ for $\alpha,\beta\in Q_M$.
\endproclaim

Now we define the Fock space
$$
V_M={\Cal S}({\Cal H}^-_{M})\otimes {\Bbb C}[\Gamma_M]
\tag 2.12$$
which affords representations for both the Lie algebra ${\Cal H}_M$ and
the group algebra ${\Bbb C}[\Gamma_M
]$ with the following actions:
 
\medskip
 
$\;\;\;\;\vep_i(k).u\otimes e^{\beta}=k({\frac
{\partial}{\partial\vep_i(-k)}}u)\otimes e^{\beta},$ for $k\in{\Bbb
Z}_+$,
 
$\;\;\;\;\vep_i(k).u\otimes e^{\beta}=(\vep_i(k)u)\otimes e^{\beta},$ for
$k\in{\Bbb
Z}_-$,
 
$\;\;\;\;\vep_i(0).u\otimes e^{\beta}=(\vep_i,\beta)u\otimes e^{\beta},$
 
$\;\;\;\;c.u\otimes e^{\beta}=u\otimes e^{\beta},$ and
$e^{\alpha}.u\otimes e^{\beta}=\epsilon(\alpha,\beta)
u\otimes e^{\alpha+\beta},$
 
\medskip
 
\noindent for $\alpha,\beta\in\Gamma_M,$ $1\le i\le M$, and $u\in {\Cal
S}({\Cal H}^-_{M})$.
 For $\alpha\in\Gamma_M$, we define (we are using the standard notation from [FLM])
$$
\alpha(z)=\sum_{k\in{\Bbb Z}}\alpha(k)z^{-k}\in ({\hbox
{End}}V_M)[[z,z^{-1}]]
\tag 2.13$$
and
$$
E^{\pm}(\alpha,z)={\hbox {exp}}(\sum_{k\in{\Bbb Z}_{\pm}}{\frac
{\alpha(k)} {k}}z^{-k})\in ({\hbox {
End}}V_M)[[z,z^{-1}]].
\tag 2.14$$
Then the follow lemma is straightforward.
 
\proclaim{Lemma 2.15} For $\alpha,\beta\in\Gamma_M$, $a,b\in {\Bbb
C}^{\times}:={\Bbb
C}\setminus\{0\}$, we have
 
\medskip
 
$\;\;\;\;E^{\pm}(0,az)=1,$ $[\alpha(k), E^{+}(\beta,az)]=0$, if $k\ge 0$,
  
$\;\;\;\;[\alpha(z_1), E^{\pm}(\beta,
az_2)]=-(\alpha,\beta)E^{\pm}(\beta, az_2)\sum_{k\in{\Bbb Z}_{\mp}}({\frac {az_2}{z_1}})^k$,
 
$\;\;\;\;E^{\pm}(\alpha, az)E^{\pm}(\beta, bz)={\hbox
{exp}}(\sum_{k\in{\Bbb Z}_{\pm}}{\frac 1 k}(a^{-k}
\alpha(k)+b^{-k}\beta(k))z^{-k}),$
 
$\;\;\;\;E^{+}(\alpha, az_1)E^{-}(\beta, bz_2)=E^{-}(\beta,
bz_2)E^{+}(\alpha, az_1)(1-{\frac {bz_2}{az_1}})
^{(\alpha,\beta)}.$
\endproclaim

Let $v=\alpha_1(-1)^{k_1}\cdots\alpha_r(-r)^{k_r}\otimes e^{\beta}\in
V_M,$ we define a degree operator
$d_0$ of $V_M$ by setting
$$d_0v=(-\sum^r_{i=1}ik_i-{\frac 1 2}(\beta,\beta))v.\tag 2.16$$
 If $a$ is any non-zero
complex number we define operators
$$
z^{\alpha}.u\otimes e^{\beta}=z^{(\alpha,
\beta)}u\otimes
e^{\beta},\;\;\;\;\;
a^{\alpha}.u\otimes
e^{\beta}=a^{(\alpha, \beta)}u\otimes e^{\beta}\tag
2.17$$
for $\alpha,\beta\in \Gamma_M$, $u\in {\Cal S}({\Cal H}^-_M).$ Then $a^{\alpha}$ is just
the evaluation map, at $a$, 
of
the operator $z^\alpha$.
The following result is well-known.
 
\proclaim{Lemma 2.18}
$$\align &
[d_0, E^{\pm}(\alpha, az)]=-D_zE^{\pm}(\alpha, az)=(\sum_{k\in {\Bbb
Z}_{\pm}}\alpha(k)(
az)^{-k})E^{\pm}(\alpha, az), \\
& [\alpha(0), z^{\beta}]=0,\;\;\;\;\;\; z^{\alpha}
e^{\beta}=z^{(\alpha,\beta)}e^{\beta}
z^{\alpha}
\endalign $$
for $\alpha, \beta\in\Gamma_M$, $a\in {\Bbb C}^{\times},$ where
$D_z={\frac d {dz}}$.
 
\endproclaim
 
\medskip
 
We will have need to raise some of the  complex numbers which arise in our construction below to various powers and care must be taken with this. We thus set up the notation we use for this now.  For any complex number $a\not= 0$, there is a unique real number 
 $\theta\in [0,2\pi)$ such that $a=|a|e^{\theta\sqrt {-1}}$. We
define 
$$
Lna=\theta\sqrt{-1}+\ln |a|
$$
Viewing ${\Bbb
C}^{\times}=\Bbb{C}\setminus\{0\}$ as multiplicative group, we
call a subgroup $G$ of $
{\Bbb C}^{\times}$ {\bf  admissible}  if $G=T\times F$, where 
$T=<\xi>$ is a cyclic group of finite order $|T|$ and
$F=<q_j|j\in J>$ is a free abelian group with free generators $q_j, j \in J$. 
For $a=\xi^{-n_0}q^{n_1}_{i_1}\cdots q^{n_k}_{i_k}\in G$, where 
$n_0,n_1,\cdots, n_k\in {\Bbb Z}, i_1,\cdots, i_k\in J$, $0\le
n_0\le |T|-1$. We define
$$
a^r=e^{r(-n_0Ln\xi+n_1Lnq_{i_1}+\cdots +n_kLnq_{i_k})}\tag 2.19
 $$
for $r\in \Bbb C$.

   Recall  the
definition of limit of formal power series from [FLM]. Let $V$ be a
vector space over $\Bbb C$. Let
$$f(z_1,z_2)=\sum_{i,j\in{\Bbb Z}}a_{ij}z^i_1z^j_2\in V[[z^{\pm
1}_1,z^{\pm 1}_2]], $$
we say the limit, $\lim_{
z_2\to z_1}f(z_1,z_2)$, exists if, for any $l\in {\Bbb Z}$,
$a_{i,l-i}=0$
whenever $|i|>>0$, and
 write $$
\lim_{z_2\to z_1}f(z_1,z_2)=f(z_1,z_1)=\sum_{l\in {\Bbb Z}}(\sum_{i\in
{\Bbb
Z}}a_{i,l-i})z_1^l.\tag 2.20$$
For our purposes we
 need another notion of limit as well. Let $f(x,z)=\sum_{i\in {\Bbb
Z}}C_i(x)z^i$, where
$C_i(x)=\sum_jc_{ij}(x)v_j\in
{\Bbb C}(x)\otimes V$, and $j$ runs over a finite set for each 
fixed $i\in
{\Bbb Z}$, and
$c_{ij}(x)\in \Bbb{C}(x)$ are complex rational  functions. We say
the limit,
$\lim_{x\to
a}f(x,z)$, exists if the function $c_{ij}(x)$ has a usual limit
at the point
$a\in {\Bbb C}$ for all $i,j\in {\Bbb
Z}$, and write $$\lim_{x\to
a}f(x,z)=f(a,z)=\sum_i(\sum_jc_{ij}(a)v_j)z^j. \tag 2.21$$
 
\proclaim{Lemma 2.22}  For $1\le i\le M$, we have
$$\align & \lim_{a\to 1}{\frac 1 {1-a}}(a^{-\vep_i}-1)=\vep_i(0), \tag
2.23\\
& \lim_{a\to 1}{\frac 1 {1-a}}(\sum_{k\in {\Bbb Z}_{\pm}}{\frac
{\vep_i(k)} k}(az)^{-k}-
\sum_{k\in {\Bbb Z}_{\pm}}{\frac {\vep_i(k)} k}z^{-k})=\sum_{k\in {\Bbb
Z}_{\pm}}\vep_i(k)z^{-k}.
\tag 2.24\endalign$$
\endproclaim
\demo{Proof} For any $v=u\otimes e^{\beta}\in V_M$, let
$m=-(\vep_i,\beta)\in {\Bbb Z}$. Then
$$
{\frac 1 {1-a}}(a^{-\vep_i}-1).v={\frac {a^m-1}{1-a}}v,
$$
 and
$$
\lim_{a\to 1}{\frac 1 {1-a}}(a^{-\vep_i}-1).v=-mv=\vep_i(0).v
$$
as required. The second identity is clear. \qed \enddemo
 
\proclaim{Corollary 2.25}
$$
\lim_{a\to 1}{\frac 1 {1-a}}(E^{\pm}(\vep_i, az)-E^{\pm}(\vep_i,
z))=E^{\pm}(\vep_i, z)
\sum_{k\in {\Bbb Z}_{\pm}}\vep_i(k)z^{-k}.
$$
\endproclaim
\demo{Proof} Note that
$$
E^{\pm}(\vep_i, az)-E^{\pm}(\vep_i, z)
=\sum^{\infty}_{l=1}{\frac 1 {l!}}[(\sum_{k\in {\Bbb Z}_{\pm}}
{\frac {\vep_i(k)} {k}}(az)^{-k})^l-(\sum_{k\in {\Bbb Z}_{\pm}}
{\frac {\vep_i(k)} {k}}z^{-k})^l],
$$
and $A^l-B^l=(A-B)\sum^{l-1}_{j=0}A^{l-1-j}B^j$, we obtain by applying
the previous lemma
$$
\lim_{a\to 1}{\frac 1 {1-a}}(E^{\pm}(\vep_i, az)-E^{\pm}(\vep_i, z))
$$
$$
=\sum_{l=1}^{\infty}{\frac 1 {l!}}[l(\sum_{k\in {\Bbb Z}_{\pm}}
{\frac {\vep_i(k)} {k}}z^{-k})^{l-1}]\sum_{k\in {\Bbb
Z}_{\pm}}\vep_i(k)z^{-k}
=E^{\pm}(\vep_i, z)\sum_{k\in {\Bbb Z}_{\pm}}\vep_i(k)z^{-k}.\qed
$$
\enddemo
 
\proclaim{Corollary 2.26}
$$
\lim_{a\to 1}{\frac 1 {1-a}}(E^{\pm}(-\vep_i, z)E^{\pm}(\vep_i, az)-1)=
\sum_{k\in {\Bbb Z}_{\pm}}\vep_i(k)z^{-k}.
$$
\endproclaim
\demo{Proof} This follows from the fact that
$$
E^{\pm}(-\vep_i, z)E^{\pm}(\vep_i, az)-1=E^{\pm}(-\vep_i,
z)(E^{\pm}(\vep_i, az)-E^{\pm}(\vep_i, z)),
$$
and the previous corollary.\qed
\enddemo
 
\medskip
 
For $\alpha\in\Gamma_M$, we define
$$
X(\alpha,z)=E^-(-\alpha,
z)E^+(-\alpha,z)e^{\alpha}z^{\alpha}z^{(\alpha,\alpha)/2}.
\tag 2.27$$
 
We may formally write
$$X(\alpha,z)=\sum_{k\in {\Bbb
Z}+(\alpha,\alpha)/2}x_k(\alpha)z^{-k},\tag 2.28$$
 where $x_k(
\alpha)\in {\hbox {End}}(V_M)$ for $k\in {\Bbb Z}+(\alpha,\alpha)/2$. It
is known from [F1]
that
, if
$(\alpha,\alpha)=1$, the operators $\{x_k(\alpha),x_k(-\alpha)|k\in {\Bbb
Z}+{\frac 1 2}\}$ generate a
Clifford algebra with the relations
 
$$
\{x_k(\alpha),x_{-l}(-\alpha)\}=\delta_{kl},\;\{x_k(\alpha),
x_{l}(\alpha)\}=0,\;
\{x_k(-\alpha),x_{l}(-\alpha)\}=0
\tag 2.29
$$
 
\noindent for all $k,l\in {\Bbb Z}+{\frac 1 2}$. Related to this Clifford
structure,
we define the following  normal ordering (see [F1] and [G3]):
 
$$
:x_k(\vep_i)x_{-l}(-\vep_j):=x_k(\vep_i)x_{-l}(-\vep_j)-\delta_{ij}\delta_{kl}
\theta(k) \tag
2.30
$$
 
\noindent for $k,l\in {\Bbb Z}+{\frac 1 2},$ $1\le i,j\le M$, where
$\theta(k)=0$ if $k<0$, $\theta(k)=1$ if
$k>0$. By applying (2.15) and (2.18), one can easily prove the following
result
 
\proclaim{Lemma 2.31} For $1\le i, j\le M$, and $a\in{\Bbb C}^{\times}$,
we have
$$
 :X(\vep_i,z_1)X(-\vep_j,az_2):
 $$
 $$
= (1-{\frac {az_2}{z_1}})^{-\delta_{ij}}{\frac
{(az_2)^{\delta_{ij}/2}}
{z_1^{\delta_{ij}/2}}} \big(\epsilon(\vep_i,\vep_j)z_1^{\frac
{(\vep_i,\vep_i-\vep_j)}{2}}\cdot (az_2)^{-\frac
{(\vep_j,\vep_i-\vep_j)}{2}}
$$
$$
 \cdot e^{\vep_i-\vep_j}z_1^{\vep_i}(az_2)^{-\vep_j}
E^-(-\vep_i,z_1)E^-(\vep_j,az_2)E^+(-\vep_i,z_1)E^+
(\vep_j,az_2)-\delta_{ij}
). 
$$ 
In particular, if $i\not= j$, then
$$ \align & :X(\vep_i,z_1)X(-\vep_j,az_2):=\epsilon(\vep_i,\vep_j)z_1^{\frac 1
2}(az_2)^{\frac 1 2}
e^{\vep_i-\vep_j}z_1^{\vep_i}(az_2)^{-\vep_j}\tag 2.32 \\
& \;\;\;\;\;\;\;\;\;\;\cdot
E^-(-\vep_i,z_1)E^-(\vep_j,az_2)E^+(-\vep_i,z_1)E^+(\vep_j,az_2),
\endalign $$
and, if $i=j$, then
$$ \align &
(1-{\frac {az_2}{z_1}}):X(\vep_i,z_1)X(-\vep_i,az_2): \tag 2.33 \\
= & {\frac {(az_2)^{\frac 1 2}}{z_1^{\frac 1 2}}}\big(({\frac {z_1}
{az_2}})^{\vep_i}
E^-(-\vep_i,z_1)E^-(\vep_i,az_2)E^+(-\vep_i,z_1)E^+(\vep_i,az_2)-1\big).
\endalign $$
\endproclaim
 
\proclaim{Proposition 2.34} For $1\le i,j\le M$, and $a\in {\Bbb
C}^{\times}$, we have
$$ \align
 & :X(\vep_i,z)X(-\vep_j,az): \\
=& \cases & \epsilon(\vep_i,\vep_j)z^{\frac 1 2}(az)^{\frac 1 2}
e^{\vep_i-\vep_j}z^{\vep_i}(az)^{-\vep_j}
E^-(-\vep_i,z)E^-(\vep_j,az)E^+(-\vep_i,z)E^+(\vep_j,az) \\
& \text{if}\;\; i\not=
j\\
&
\vep_i(z)
\;\;\;\;\;\;\;\;\;\;\;\;\;\;\;\;\;\;\;\;\;\;\;\;\;\;\;\;\;\;\;\;\;
\text{if}\;\;  i=j, a=1\\
& {\frac
{a^{1/2}}{1-a}}(a^{-\vep_i}E^-(-\vep_i,z)E^-(\vep_i,az)E^+(-\vep_i,z)
E^+(\vep_i,az)-1)\\
& \text{if}\;\; i=j, a\not= 1.\endcases
\endalign $$
 \endproclaim
 
\demo{Proof} Taking the limit $z_2 \to z_1$ in (2.32) and (2.33) gives
the first and third
identities. The second identity follows from Lemma 2.22, Corollary 2.26 and 
the
third identity by taking
the limit $a\to 1$. \qed \enddemo
 
\noindent{\bf Remark 2.35.} Note that the second identity in Proposition
2.34 was given in
[F1].
 
\medskip
 
\proclaim{Definition 2.36} 
For $a\in {\Bbb C}^{\times}$, $1\le i,j\le M$, we define
 $
X_{ij}(a,z)=:X(\vep_i,z)X(-\vep_j,az):
$
\endproclaim
 \medskip
 
Now we can state our main theorem of this section.
 
\proclaim{Theorem 2.37} 
For $a,b\in {\Bbb C}^{\times}$ and $1\le i,j,k,l\le M$, we have
 
 (i) if $ab\not= 1$, then
$$\align &
[X_{ij}(a,z_1), X_{kl}(b,z_2)]=X_{il}(ab,z_1)\delta_{jk}\delta({\frac {z_2}
{az_1}})-X_{kj}(ab,z_2)\delta_{il}\delta({\frac {z_1} {bz_2}}) \tag 2.38
\\
& \;\;\;\;\;\;\;\;+{\frac {a^{\frac 1 2}b^{\frac 1
2}}{1-ab}}\delta_{il}\delta_{jk}(\delta({\frac {z_2}
{az_1}})
-\delta({\frac {z_1} {bz_2}}))c,
\endalign $$
 
(ii) if $ab=1$, then
$$ \align
 & [X_{ij}(a,z_1),
X_{kl}(b,z_2)]=(X_{il}(1,z_1)\delta_{jk}-X_{kj}(1,z_2)\delta_{il})\delta({\frac
{z_2} {az_1}})
+\delta_{il}\delta_{jk}(D\delta)({\frac {z_2} {az_1}})c.\tag 2.39 \endalign $$
 
\endproclaim
 
The proof of Theorem 2.37 will be carried out in several steps. 
In what follows we will
freely use the following two lemmas. (2.40) can be found in [FLM] and
[K],  and (2.43) can be found
in [J], [G1, 2] or [BS].
 
\proclaim{Lemma 2.40} Let $Y(z_1,z_2)$ be a formal power series in
$z_1,z_2$ with coefficients
in a vector
space, such that  $\lim_{z_2\to z_1}f(z_1,z_2)$ exists. Then
$$
Y(z_1,z_2)\delta({\frac {az_1}{z_2}})=Y(z_1,az_1)\delta({\frac
{az_1}{z_2}}), \tag 2.41
$$
$$
Y(z_1,z_2)(D\delta)({\frac {z_2}{az_1}})=Y(z_1,az_1)(D\delta)({\frac
{z_2}{az_1}})-
(D_{z_2}Y(z_1,z_2))\delta({\frac {z_2}{az_1}}),\tag 2.42
$$
for $a\in \Bbb{C}^\times$.
\endproclaim

\proclaim{Lemma 2.43} Suppose $a,b\in {\Bbb C}^{\times}$. Then
$$
(1-{\frac {z_2}{az_1}})^{-1}(1-{\frac {bz_2}{z_1}})^{-1}-{\frac
{az_1}{z_2}}{\frac {z_1}{
bz_2}}(1-{\frac {z_1}{bz_2}})^{-1}(1-{\frac {az_1}{z_2}})^{-1}
$$
$$
=\cases (1-ab)^{-1}{\frac {az_1}{z_2}}(\delta({\frac
{az_1}{z_2}})-\delta({\frac {z_1}{bz_2}})) &
\text{if}\;\; ab\not=
1\\
{\frac {az_1}{z_2}}(D\delta)({\frac {z_2}{az_1}}) &\text{if}\;\;
ab=1.\endcases$$
 
\endproclaim
 
\medskip
 
Now we divide the proof for Theorem 2.37 into four different cases.

\medskip
 
{\bf Case 1.} $i\not= j,$ $k\not= l$.
 
\medskip
 
We obtain, by applying (2.15) and (2.18), that
$$
\align & [X_{ij}(a, z_1),X_{kl}(b,
z_2)]\tag 2.44\\
= & \epsilon(\vep_i,\vep_j)\epsilon(\vep_k,\vep_l)a^{\frac 1 2}b^{\frac 1
2}
e^{\vep_i-\vep_j}e^{\vep_k-\vep_l}a^{-\vep_j}b^{-\vep_l}z_1^{\vep_i-\vep_j}
z_2^
{\vep_k-\vep_l}z_1z_2
 \\
& \cdot E^-(-\vep_i,z_1)E^-(-\vep_k,z_2)E^-(\vep_j,az_1)E^-(\vep_l,bz_2)
\\
& \cdot E^+(-\vep_i,z_1)E^+(-\vep_k,z_2)
E^+(\vep_j,az_1)E^+(\vep_l,bz_2)P(z_1,z_2)
\endalign $$
where
$$\align &
P(z_1,z_2)\\
=&a^{-(\vep_j,\vep_k-\vep_l)}z_1^{(\vep_i-\vep_j,\vep_k-\vep_l)}(1-{\frac
{z_2}{z_1}})^{\delta_
{ik}}(1-{\frac {z_2}{az_1}})^{-\delta_
{jk}}(1-{\frac {bz_2}{z_1}})^{-\delta_
{il}}(1-{\frac {bz_2}{az_1}})^{\delta_
{jl}} \endalign
$$
$$
-(-1)^{(\vep_i-\vep_j,\vep_k-\vep_l)}b^{-(\vep_l,\vep_i-\vep_j)}
z_2^{(\vep_k-\vep_l,\vep_i-\vep_j)}(1-
{\frac {z_1}{z_2}})^{\delta_
{ik}}(1-{\frac {z_1}{bz_2}})^{-\delta_
{il}}(1-{\frac {az_1}{z_2}})^{-\delta_
{jk}}(1-{\frac {az_1}{bz_2}})^{\delta_
{jl}}
$$
$$
=a^{-(\vep_j,\vep_k-\vep_l)}z_1^{(\vep_i-\vep_j,\vep_k-\vep_l)}(1-{\frac
{z_2}{z_1}})^{\delta_
{ik}}(1-{\frac {bz_2}{az_1}})^{\delta_
{jl}}\left((1-{\frac {z_2}{az_1}})^{-\delta_
{jk}}(1-{\frac {bz_2}{z_1}})^{-\delta_
{il}}\right.
$$
$$
-\left.(-1)^{\delta_{il}+\delta_{jk}}({\frac
{z_2}{az_1}})^{-\delta_{jk}}({\frac {bz_2}{z_1}})^{-\delta_{il}}
(1-{\frac {az_1}{z_2}})^{-\delta_
{jk}}(1-{\frac {z_1}{bz_2}})^{-\delta_
{il}}\right).
$$
Applying Lemma 2.43 we have two subcases.
 
\medskip
 
{\it Subcase 1}. If $ab\not=1$, then
$$
P(z_1,z_2)=\cases 0 &\text{if}\;\;i\not= l, j\not= k\\
{\frac 1{(1-ab)z_1z_2}}(\delta({\frac {z_2}{az_1}})-\delta({\frac
{z_1}{bz_2}})) &\text{if}\;\;i=l, j=k\\
z_1^{-1}\delta({\frac {z_1}{bz_2}}) &\text{if}\;\;i=l, j\not= k\\
(az_1)^{-1}\delta({\frac {z_2}{az_1}}) &\text{if}\;\;i\not= l, j=
k.\endcases
$$
Now (2.38) follows from (2.44) and Lemma 2.40.
 
\medskip
 
{\it Subcase 2}. If $ab=1$, then
$$
P(z_1,z_2)=a^{-(\vep_j,\vep_k-\vep_l)}z_1^{(\vep_i-\vep_j,\vep_k-\vep_l)}
(1-{\frac {z_2}{z_1}})^{\delta_
{ik}}(1-{\frac {bz_2}{az_1}})^{\delta_
{jl}}((1-{\frac {z_2}{az_1}})^{-\delta_
{jk}-\delta_{il}})
$$
$$
-(-1)^{\delta_{il}+\delta_{jk}}({\frac
{z_2}{az_1}})^{-\delta_{jk}-\delta_{il}}
(1-{\frac {az_1}{z_2}})^{-\delta_
{jk}-\delta_{il}}
$$
$$
=\cases 0 &\text{if}\;\;i\not= l, j\not= k\\
{\frac 1{z_1z_2}}(D\delta)({\frac {z_2}{az_1}}) &\text{if}\;\;i=l, j=k\\
z_1^{-1}\delta({\frac {z_2}{az_1}}) &\text{if}\;\;i=l, j\not= k\\
(az_1)^{-1}\delta({\frac {z_2}{az_1}}) &\text{if}\;\;i\not= l, j=
k,\endcases
$$
which yields (2.39) by applying Lemma 2.40.
 
\medskip
 
{\bf Case 2}. $i= j,$ $k\not= l$, and $ab\not= 1$.
 
\medskip
 
If $a=1$, then
$$
X_{ij}(a, z_1)=\vep_i(z_1).
$$
Since
$$\align
& [\vep_i(z_1),
e^{\vep_k-\vep_l}]=(\delta_{ik}-\delta_{il})e^{\vep_k-\vep_l},\\
& [\vep_i(z_1),
E^{\pm}(-\vep_k,z_2)]=\sum_{n\in{\Bbb Z}_{\pm}}\delta_{ik}({\frac
{z_1}{z_2}})^nE^{\pm}(-\vep_k,z_2),
\endalign $$
and
$$
[\vep_i(z_1),E^{\pm}(\vep_l,bz_2)]=-\sum_{n\in{\Bbb
Z}_{\pm}}\delta_{il}({\frac {z_1}{bz_2}})^n
E^{\pm}(\vep_l,bz_2),
$$
 we have
$$
  [X_{ij}(a, z_1),X_{kl}(b, z_2)]=[\vep_i(z_1),X_{kl}(b, z_2)]=
X_{kl}(b, z_2)
$$
$$\cdot\left(\delta_{ik}-\delta_{il}+
\delta_{ik}\sum_{n\in{\Bbb Z}_{-}}({\frac {z_1}{z_2}})^n+
\delta_{ik}\sum_{n\in{\Bbb Z}_{+}}({\frac
{z_1}{z_2}})^n-\delta_{il}\sum_{n\in{\Bbb Z}_{-}}
({\frac {z_1}{bz_2}})^n
-\delta_{il}\sum_{n\in{\Bbb Z}_{+}}({\frac {z_1}{bz_2}})^n\right)
$$
$$ \align
= & X_{kl}(b, z_2)(\delta_{ik}\delta({\frac
{z_2}{z_1}})-\delta_{il}\delta({\frac
{z_1}{bz_2}})) \\
= &X_{kl}(b, z_1)\delta_{jk}\delta({\frac {z_2}{z_1}})-X_{kj}(b,
z_2)\delta_{il}\delta({\frac
{z_1}{bz_2}}),
\endalign $$
as needed.
 
If $a\not= 1$, then
$$
X_{ij}(a, z_1)={\frac
{a^{1/2}}{1-a}}(a^{-\vep_i}E^-(-\vep_i,z_1)E^-(\vep_i,az_1)
E^+(-\vep_i,z_1)E^+(\vep_i,az_1)-1),
$$
Applying (2.15) and (2.18), we get
$$\align
& [X_{ij}(a, z_1),X_{kl}(b, z_2)]\tag 2.45 \\
=& {\frac {a^{1/2}}{1-a}}b^{\frac 1
2}\epsilon(\vep_k,\vep_l)e^{\vep_k-\vep_l}
a^{-\vep_i}b^{-\vep_l}z_2^{\vep_k-\vep_l}z_2 \\
& \cdot
E^-(-\vep_i,z_1)E^-(\vep_i,az_1)E^-(-\vep_k,z_2)E^-(\vep_l,bz_2)\\
& \cdot E^+(-\vep_i,z_1)E^+(\vep_i,az_1)
E^+(-\vep_k,z_2)E^+(\vep_l,bz_2)Q(z_1,z_2),
\endalign $$
where
$$\align &
Q(z_1,z_2)=a^{-(\vep_j,\vep_k-\vep_l)}(1-{\frac {z_2}{z_1}})^{\delta_
{ik}}(1-{\frac {z_2}{az_1}})^{-\delta_
{ik}}(1-{\frac {bz_2}{z_1}})^{-\delta_
{il}}(1-{\frac {bz_2}{az_1}})^{\delta_
{il}}\endalign
$$
$$
-(1-
{\frac {z_1}{z_2}})^{\delta_
{ik}}(1-{\frac {z_1}{bz_2}})^{-\delta_
{il}}(1-{\frac {az_1}{z_2}})^{-\delta_
{ik}}(1-{\frac {az_1}{bz_2}})^{\delta_
{il}}
$$
$$
=a^{\delta_{il}-\delta_{ik}}(1-{\frac {z_2}{z_1}})^{\delta_
{ik}}(1-{\frac {bz_2}{az_1}})^{\delta_
{il}}\left((1-{\frac {z_2}{az_1}})^{-\delta_
{ik}}(1-{\frac {bz_2}{z_1}})^{-\delta_
{il}}\right.
$$
$$
\left.-(-1)^{\delta_{ik}+\delta_{il}}({\frac {az_1}{z_2}})^{\delta_{ik}}
({\frac {z_1}{bz_2}})^{\delta_{il}}(1-{\frac {z_1}{bz_2}})^{-\delta_
{il}}(1-{\frac {az_1}{z_2}})^{-\delta_
{ik}}\right)
$$
$$
=\cases 0 &\text{if}\;\;i\not= k, i\not= l\\
{\frac {1-a}a}\delta({\frac {z_2}{az_1}}) &\text{if}\;\;i= k, i\not= l\\
(a-1)\delta({\frac {z_1}{bz_2}}) &\text{if}\;\;i\not= k, i= l\endcases
$$
 We thus have (2.38) by applying Lemma 2.40.
 
\medskip
 
{\bf Case 3}. $i= j,$ $k\not= l$, and $ab= 1$.
 
\medskip
 
 As we did in Case 2, we  treat the case $a=1$ and $a\not= 1$
separately. First for
$a=1$ (so $b=1$),
 then
 $$
 X_{ij}(a,z_1)=\vep_i(z_1),
 $$
 $$
 X_{kl}(b,z_2)=\epsilon(\vep_k,\vep_l)
e^{\vep_k-\vep_l}z_2^{\vep_k-\vep_l}z_2
E^-(-\vep_k,z_2)E^-(\vep_l,z_2)E^+(-\vep_k,z_2)E^+(\vep_l,z_2),
$$
and we have
$$\align &
[X_{ij}(a,z_1),X_{kl}(b,z_2)]=[\vep_i(z_1),X_{kl}(1,z_2)]\\
=& X_{kl}(1,z_2)(\delta_{ik}\delta({\frac
{z_2}
{z_1}})-\delta_{il}\delta({\frac {z_2}
{z_1}}))\\
=&(X_{il}(1,z_1)\delta_{jk}-X_{kj}(1,z_2)\delta_{il})\delta({\frac {z_2}
{z_1}}),\endalign
$$
as expected.
 
Next assume that $a\not= 1$. By applying similar arguments as  in
Case 2, we have
$$
\align & [X_{ij}(a,z_1),X_{kl}(b,z_2)]\tag 2.46\\
=& {\frac 1
{1-a}}\epsilon(\vep_k,\vep_l)e^{\vep_k-\vep_l}a^{-\vep_i}a^{\vep_l}
z_2^{\vep_k-\vep_l}z_2\\
& \cdot E^-(-\vep_i,z_1)
 E^-(\vep_i,az_1)E^-(-\vep_k,z_2)E^-(\vep_l,bz_2)
\\
& \cdot E^+(-\vep_i,z_1)E^+(\vep_i,az_1)
 E^+(-\vep_k,z_2)E^+(\vep_l,bz_2)R(z_1,z_2)
\endalign $$
where
$$\align &
R(z_1,z_2)=a^{\delta_{il}-\delta_{ik}}(1-{\frac
{z_2}{z_1}})^{\delta_{ik}}(1-{\frac {bz_2}{az_1}})^{
\delta_{il}}\left((1-{\frac
{z_2}{az_1}})^{-\delta_{il}-\delta_{ik}}\right.
\endalign$$
$$
\left.-(-1)^{\delta_{il}+\delta_{ik}}({\frac
{az_1}{z_2}})^{\delta_{il}+\delta_{ik}}
(1-{\frac {az_1}{z_2}})^{-\delta_{il}-\delta_{ik}}\right)
$$
$$
=\cases 0 &\text{if}\;\;i\not= k, i\not= l\\
(b-1)\delta({\frac {z_2}{az_1}}) &\text{if}\;\;i= k, i\not= l\\
(a-1)\delta({\frac {z_2}{az_1}}) &\text{if}\;\;i\not= k, i= l.\endcases
$$
Substituting this into (2.46) and applying Lemma 2.40, we get (2.39).
 
\medskip
 
{\bf Case 4}. $i=j, k=l$.
 
\medskip
 
Note that
$$
X_{ij}(a,z_1)
$$
$$
=\cases {\frac
{a^{1/2}}{1-a}}(a^{-\vep_i}E^-(-\vep_i,z_1)E^-(\vep_i,az_1)E^+(-\vep_i,
z_1)E^+(\vep_i,az_1)-1) &
\text{if}\;\; i=j, a\not= 1\\
\vep_i(z_1) &\text{if}\;\;  i=j, a=1
\endcases
$$
$$
X_{kl}(b,z_2)
$$
$$
=\cases {\frac
{b^{1/2}}{1-b}}(b^{-\vep_k}E^-(-\vep_k,z_2)E^-(\vep_k,bz_2)E^+(-\vep_k,
z_2)E^+(\vep_k,bz_2)-1) &
\text{if}\;\; k=l, b\not= 1\\
\vep_k(z_2) &\text{if}\;\;  k=l, b=1.
\endcases
$$
 
We consider three subcases. First we assume that $a=b=1$, then
$$[X_{ij}(a,z_1),
X_{kl}(b,z_2)]  = [\vep_i(z_1),\vep_k(z_2)] =\delta_{ik}(D\delta)({\frac
{z_2}{z_1}}),$$ as desired.
 
Next we assume $a= 1$, $b\not= 1$, then
$$
[X_{ij}(a,z_1),X_{kl}(b,z_2)]=[\vep_i(z_1),X_{kl}(b,z_2)]
$$
$$
={\frac
{b^{1/2}}{1-b}}b^{-\vep_k}E^-(-\vep_k,z_2)E^-(\vep_k,bz_2)E^+(-\vep_k,
z_2)E^+(\vep_k,bz_2)
$$
$$\cdot\left(\delta_{ik}\sum_{n\in {\Bbb Z}_{-}}({\frac
{z_2}{z_1}})^n
+
\delta_{ik}\sum_{n\in {\Bbb Z}_{+}}({\frac
{z_2}{z_1}})^n-\delta_{ik}\sum_{n\in {\Bbb Z}_{-}}(
{\frac
{z_1}{bz_2}})^n-
\delta_{ik}\sum_{n\in {\Bbb Z}_{+}}({\frac {z_1}{bz_2}})^n\right)
$$
$$
={\frac
{b^{1/2}}{1-b}}b^{-\vep_k}E^-(-\vep_k,z_2)E^-(\vep_k,bz_2)E^+(-\vep_k,
z_2)E^+(\vep_k,bz_2)\delta_{ik}(\delta({\frac {z_2}{z_1}})-\delta({\frac
{z_1}{bz_2}}))
$$
$$
=X_{il}(ab,z_1)\delta_{jk}\delta({\frac {z_2}
{az_1}})-X_{kj}(ab,z_2)\delta_{il}\delta({\frac {z_1} {bz_2}})
+{\frac {a^{\frac 1 2}b^{\frac 1
2}}{1-ab}}\delta_{il}\delta_{jk}(\delta({\frac {z_2} {az_1}})
-\delta({\frac {z_1} {bz_2}})),
$$
as required.
 
Finally we assume $a\not= 1$, $b\not= 1$, then
$$
[X_{ij}(a,z_1),X_{kl}(b,z_2)]\tag 2.47
$$
$$
={\frac {a^{1/2}}{1-a}}{\frac {b^{1/2}}{1-b}}a^{-\vep_i}b^{-\vep_k}
E^-(-\vep_i,z_1)E^-(-\vep_k,z_2)E^-(\vep_i,az_1)
$$
$$
\cdot E^-(\vep_k,bz_2)E^+(-\vep_i,z_1)E^+(-\vep_k,z_2)E^+(\vep_i,az_1)E^+(\vep_k,bz_2
S(z_1,z_2),
$$
where
$$\align &
S(z_1,z_2)=(1-{\frac {z_2}{z_1}})^{\delta_{ik}}(1-{\frac
{bz_2}{z_1}})^{-\delta_{ik}}
(1-{\frac {z_2}{az_1}})^{-\delta_{ik}}
(1-{\frac {bz_2}{az_1}})^{\delta_{ik}}
\endalign $$
$$
-(1-{\frac {z_1}{z_2}})^{\delta_{ik}}(1-{\frac
{z_1}{bz_2}})^{-\delta_{ik}}
(1-{\frac {az_1}{z_2}})^{-\delta_{ik}}
(1-{\frac {az_1}{bz_2}})^{\delta_{ik}}
$$
$$
=(1-{\frac {z_2}{z_1}})^{\delta_{ik}}(1-{\frac
{bz_2}{az_1}})^{\delta_{ik}}
\left((1-{\frac {bz_2}{z_1}})^{-\delta_{ik}}
(1-{\frac {z_2}{az_1}})^{-\delta_{ik}}\right.
$$
$$
\left.-({\frac {az_1}{z_2}})^{\delta_{ik}}({\frac
{z_1}{bz_2}})^{\delta_{ik}}
(1-{\frac {z_1}{bz_2}})^{-\delta_{ik}}
(1-{\frac {az_1}{z_2}})^{-\delta_{ik}}\right)
$$
$$
=\cases (1-{\frac {z_2}{z_1}})(1-{\frac {bz_2}{az_1}}){\frac
{az_1}{z_2}}(D\delta)({\frac {z_2}{az_1}})
 &\text{if}\;\; i=k, ab= 1\\
0 &\text{if}\;\;i\not=k\\
{\frac 1 {1-ab}}(1-{\frac {z_2}{z_1}})(1-{\frac {bz_2}{az_1}}){\frac
{az_1}{z_2}}(\delta
({\frac {z_2}{az_1}})-\delta
({\frac {z_1}{bz_2}})) &\text{if}\;\; i=k, ab\not= 1.
\endcases
$$
Substitute $S(z_1, z_2)$ back into (2.47) and apply Lemma 2.40 to get
(2.38) and (2.39).
 This now completes the proof of Theorem 2.37.
 
\medskip
 
\noindent{\bf Remark 2.48} If $ab\not= 1$, we have
$$
\delta({\frac {z_2} {az_1}})
-\delta({\frac {z_1} {bz_2}})=\sum_{m\in{\Bbb Z}}(1-(ab)^m)({\frac {z_2}
{az_1}})^m,
$$
thus
$$
{\frac {a^{\frac 1 2}b^{\frac 1
2}}{1-ab}}\delta_{il}\delta_{jk}(\delta({\frac {z_2} {az_1}})
-\delta({\frac {z_1} {bz_2}}))=a^{\frac 1 2}b^{\frac 1
2}\delta_{il}\delta_{jk}
\sum_{m\in{\Bbb Z}}{\frac {1-(ab)^m}{1-ab}}({\frac {z_2} {az_1}})^m
$$
$$
=a^{\frac 1 2}b^{\frac 1 2}\delta_{il}\delta_{jk}(\sum_{m\in{\Bbb
Z}_+}(\sum^{m-1}_{s=0}(ab)
^s)
({\frac {z_2} {az_1}})^m+\sum_{m\in {\Bbb
Z}_-}(-\sum^{-m}_{s=1}(ab)^{-s})
({\frac {z_2} {az_1}})^m),
$$
which gives
$$
\lim_{b\to a^{-1}}{\frac {a^{\frac 1 2}b^{\frac 1
2}}{1-ab}}\delta_{il}\delta_{jk}
(\delta({\frac {z_2} {az_1}})
-\delta({\frac {z_1} {bz_2}}))
$$
$$
=\delta_{il}\delta_{jk}(\sum_{m\in{\Bbb Z}_+}m
({\frac {z_2} {az_1}})^m+\sum_{m\in {\Bbb Z}_-}m({\frac {z_2}
{az_1}})^m)=
\delta_{il}\delta_{jk}
(D\delta)({\frac {z_2} {az_1}}).
$$
This indicates that the second identity of Theorem 2.37 can be obtained
from the first one by   
taking the limit
as $b\to a^{-1}$.
 
\medskip
 
\subhead \S 3. Lie Algebras and Representations \endsubhead
 
 \medskip

In this section we are going to define a class of Lie algebras from our
vertex operators which will correspond to  admissible subgroups of $\bc^{\times}$. Indeed for some
choice of the admissible 
group $G$ and positive integer $M$, the Lie algebra ${\Cal G}(G,M)$
(defined below) of operators,
which act on the
Fock space $V_M$, will give realizations of some infinite dimensional Lie
algebras studied in Section 1.  This will include the affine
algebra ${\wh {gl_M}}({\Bbb C})$ in both the principal and homogeneous
pictures as well as  some Lie algebras
with quantum torus coordinates. 
Towards this end, we first introduce some new notation for the vertex operators constructed in the
proceeding section.
 
\medskip
 
\proclaim{Definition 3.1}  For $a,b\in {\Bbb C}^{\times}$, $1\le i,j\le
M$, we set
$X_{ij}(a,b,z):=X_{ij}(
 a^{-1}b, az)$, and write
 $$
 X_{ij}(a,b,z)=\sum_{k\in {\Bbb Z}}x_{ij}(k,a,b)z^{-k} \tag 3.2 $$
 where $x_{ij}(k,a,b)\in$ End$V_M$.
\endproclaim
 
With this notation Theorem 2.37 can be re-written as follows.
 
\medskip
 
\proclaim{Theorem 3.3}  
Let  $a_1,a_2,b_1,b_2\in {\Bbb C}^{\times}$, and $1\le
i,j,k,l\le M$. We have
 
\medskip
 
 (i) if $a_1a_2\not= b_1b_2$, then
$$\align &
[X_{ij}(a_1,b_1,z_1), X_{kl}(a_2,b_2,z_2)]\\
=& X_{il}(a_1,{\frac
{b_1b_2}{a_2}},z_1)\delta_{jk}\delta({\frac {a_2z_2} {b_1z_1}})
-
X_{kj}(a_2,{\frac {b_1b_2}{a_1}},z_2)\delta_{il}
\delta({\frac {a_1z_1} {b_2z_2}})\\
& +\frac{(a_1^{-1}b_1)^{\frac 1 2}(a_2^{-1}b_2)^{frac 1
2}}{1-a_1^{-1}b_1a_2^{-1}b_2}
\delta_{il}\delta_{jk}(\delta({\frac {a_2z_2} {b_1z_1}})
-\delta({\frac {a_1z_1} {b_2z_2}}))c,
\endalign $$
 
(ii) if $a_1a_2= b_1b_2$, then
$$
[X_{ij}(a_1,b_1,z_1), X_{kl}(a_2,b_2,z_2)]
$$
$$
=X_{il}(a_1,{\frac {b_1b_2}{a_2}},z_1)\delta_{jk}\delta({\frac {a_2z_2}
{b_1z_1}})
-
X_{kj}(a_2,{\frac {b_1b_2}{a_1}},z_2)\delta_{il}\delta({\frac {a_1z_1}
{b_2z_2}})
+
\delta_{il}\delta_{jk}(D\delta)({\frac {a_2z_2} {b_1z_1}})c.$$
 
\endproclaim
 
\medskip

Fix an integer $M\ge 1$ and an admissible subgroup $G$ of ${\Bbb
C}^{\times}$.
 Let
${\Cal G}(G,M)$ be the vector space spanned by $c$ and all of the
coefficients of the vertex operators
$X_{ij}(a,b,z)$ for all $1\le i,j\le M$, and $a,b\in G$. Then
we have the following result.
 
\medskip
 
\proclaim{Theorem 3.4} ${\Cal G}(G,M)$
 forms a Lie algebra of operators acting on the Fock space  $V_M$.
Moreover
 $$
 V_M=\oplus_{k\in{\Bbb Z}}V_M^{(k)}
 $$
 where $V_M^{(k)}=e^{k\vep_M+Q_M}\otimes{\Cal S}({\Cal H}^-_M)$, and
$V_M^{(k)}$ is an irreducible
 ${\Cal G}(G,M)$-module.
 \endproclaim
 
\demo{Proof}  It is obvious from Theorem 3.3 that ${\Cal G}(G,M)$ is a Lie algebra and that $V_M^{(k)}$ is a ${\Cal
G}(G,M)$-module.
To see it is
irreducible,
 we note that the Heisenberg algebra ${\wh {\Cal H}}_M\subset$ ${\Cal
H}_M$, and
 ${\Cal H}_M$ is spanned by the coefficient operators of the vertex
operators $X_{ii}(1,1,z)$ for $
 1\le i\le M$. This then implies that,
if $W$ is a non-zero submodule of $V_M^{(k)}$, we can choose
a non-zero element of the form $v=e^{k\vep_M+\alpha}\otimes 1\in
W$ for some $\alpha\in Q_M$.
 Moreover, it is easy to check that
 $$
 x_{ij}(n_{ij}-1,1,1).v=\epsilon(\vep_i,\vep_j)\epsilon(\vep_i-\vep_j,
 k\vep_M+\alpha)e^{k\vep_M+\alpha+\vep_i-\vep_j}
 $$
 for all $1\le i\not= j\le M$, where $n_{ij}=(\vep_j-\vep_i,
k\vep_M+\alpha)\in {\Bbb
Z}$. Therefore $
 e^{k\vep_M+\beta}\otimes 1\in W$ for any $\beta\in Q_M$. This
thus gives $W=V_M^{(k)}
 $ as needed.\qed
 
\enddemo
 
 \medskip
 
\noindent{\bf Remark 3.5} Note that the coefficients of the vertex
operators $X_{ij}(a, z)$ and $X_{ij}(a, bz)$ for any given $a, b\in G$
span the same space. Thus $\Cal{G}(G, M)$ is spanned by $c$ and the
coefficients of the operators $X_{ij}(a, z)$ for $1\leq i, j\leq M, a\in
G$. Therefore, it follows from Theorem 4.25 in [G1] that $\Cal{G}(G, M)$
is an affinization of the Lie algebra $gl_M(\Cal{R}[t, t^{-1}; \tau])$,
where
$\Cal{R}= \Bbb{C}[G]$ is the group algebra and $\Cal{R}[t, t^{-1}; \tau]$
is the skew Laurent
polynomial ring.
 
Recall definition (2.19). We extend the cocycle map $\epsilon:$ $H_M\times
H_M\to\{\bc^{\times}\}$ by
defining
$$
\epsilon(\sum r_i\vep_i,\sum
s_i\vep_i)=\prod_{i,j}(\epsilon(\vep_i,\vep_j))^{r_is_j}\tag 3.6
$$
for $r_i,s_i\in \bc$. It is obvious that
$$
\epsilon(\alpha+\beta,\gamma)=\epsilon(\alpha,\gamma)\epsilon(\beta,\gamma),\;\
\;\;\epsilon(\alpha, \beta+\gamma)
=\epsilon(\alpha,\beta)\epsilon(\alpha,\gamma)
$$
for $\alpha,\beta,\gamma\in H_M$. Moreover, if we restrict $\epsilon$ to
$\Gamma_M\times\Gamma_M$, then $\epsilon$ gives us
the 2-cocycle defined in previous section.
 
  Let $\alpha_i=\vep_i-\vep_{i+1},$ $i=1,\cdots, M-1$, and
$\alpha_M=\vep_1+\cdots+\vep_M$. Then $Q_M=\oplus_{i=1}^{M-1}{\Bbb
Z}\alpha_i$, and $\Gamma_M=\oplus_{i=1}^{M}{\Bbb Z}\alpha_i$. Let
$Q^0_M=\{\alpha\in\bc\otimes_{\Bbb Z}Q_M|\;\;(\alpha, Q_M)\in {\Bbb Z}
\}$ be the dual of the lattice $Q_M$, and set
$$
L^0_M=\{\alpha\in H_M=\bc\otimes_{\Bbb Z}\Gamma_M|\;\;(\alpha, Q_M)\in
{\Bbb Z}
\}.\tag 3.7
$$
Then $L^0_M=Q^0_M\oplus {\bc}\alpha_M$. Let $I=L^0_M/Q_M$, then we have a
$Q_M$-coset decomposition of $L^0_M$
$$
L^0_M=\oplus_{i\in I}(\lambda_i+Q_M)\tag 3.8
$$
for some $\lambda_i\in L^0_M.$

\medskip
 \proclaim
{\bf Proposition 3.9} The Lie algebra of operators, ${\Cal G}(G,M)$,  acts on
the space $V^{(\lambda_i)}_M=$ ${\Cal S}({\Cal H}^-_M)\otimes
{\bc}[Q_M+\lambda_i]$, and $V^{(\lambda_i)}_M$ affords an irreducible
representation of ${\Cal G}(G,M)$ for $i\in I=L^0_M/Q_M$.
Moreover $V^{(\lambda_i)}_M\cong V^{(\lambda_j)}_M$ if and only if $i=j$.
 \endproclaim

\demo{Proof} For $u\otimes e^{\alpha+\lambda_i}\in V^{(\lambda_i)}_M$, $u\in
{\Cal S}({\Cal H}^-_M),$
$\alpha\in Q_M$, we note that
$$
z^{\beta}.(u\otimes
e^{\alpha+\lambda_i})=z^{(\beta,\alpha+\lambda_i)}u\otimes
e^{\alpha+\lambda_i},\;\;\;\;a^{\gamma}.(u\otimes 
e^{\alpha+\lambda_i})=
a^{(\gamma,\alpha+\lambda_i)}u\otimes e^{\alpha+\lambda_i}
$$
$$
e^{\beta}.(u\otimes
e^{\alpha+\lambda_i})=\epsilon(\beta,\alpha+\lambda_i)u\otimes
e^{\alpha+\beta+\lambda_i}
$$
where $(\beta,\alpha+\lambda_i)\in{\Bbb Z},$ $
(\gamma,\alpha+\lambda_i)\in {\bc},$ $\epsilon(\beta,\alpha+\lambda_i)
\in{\bc}^{\times}$ for $\beta\in Q_M,$ $\gamma\in \Gamma_M$ and $a\in G$.
This implies that $X_{ij}(a,b,z).(
u\otimes e^{\alpha+\lambda_i})\in V^{(\lambda_i)}_M[[z,z^{-1}]],$ and so
the Lie algebra ${\Cal G}(G,M)$
acts on the space $V^{(\lambda_i)}_M$. The irreducibility of 
$V^{(\lambda_i)}_M$, for $i\in I$, follows from a similar argument
as  in Theorem 3.4. The last part of this proposition is clear. \qed
 \enddemo

\bigskip
 
\subhead \S 4. Applications \endsubhead
 
 \medskip

 In this section we assume  the admissible subgroup $G$ has the form
$G=T\times F\subset 
{\bc}^{\times}$, where $T=<\xi>$ is generated
by a root of unity $\xi$, and $F$ is a free group with a finite number generators.
First, let $G=\{1\}$, and $M\ge 2$ be any integer. Then the Lie algebra
${\Cal G}(G,M)$ is generated
by the coefficients of the vertex operators $X_{ij}(1,1,z)$ for all $1\le
i,j\le M$. Moreover
from Theorem 3.3 , we see that
$$\align &
[X_{ij}(1,1,z_1),
X_{kl}(1,1,z_2)] \tag 4.1\\
=& X_{il}(1,1,z_1)\delta_{jk}\delta({\frac {z_2}{z_1}})-
X_{kj}(1,1,z_2)\delta_{il}\delta({\frac
{z_1}{z_2}})+\delta_{il}\delta_{jk}(D\delta)({\frac
{z_2}{z_1}}) \endalign$$
for $1\le i,j\le M$. Comparing (4.1) with  (1.17), we obtain the following
result which was originally due to
[F1], see also [FK] and [S].
 
\medskip
 
\proclaim{Corollary 4.2} Let $G=\{1\}$ and $M\ge 2$. Then ${\Cal G}(G,M)$
gives a
representation of the affine
algebra ${\wh {gl_M}}({\Bbb C})$ in the homogeneous picture on the Fock
space $V_M$, and the representation
is given by the mapping:
$$\align &
E_{ij}\otimes t^k_0\mapsto x_{ij}(k,1,1), \\
& c_0\mapsto c
\endalign $$
for $1\le i,j\le M$ and $k\in {\Bbb Z}$.
\endproclaim
 
\medskip
 
Next we choose $G$ to be a cyclic group of order $N\ge 2$ with generator
$\xi=\xi_N$, and take $M=1$.  Note that
$$X_{11}(\xi^i,\xi^j,z)=X_{11}(\xi^{i-j-1}, \xi^{-1}, \xi^{j+1}z)$$
 for
$0\le i,j\le N-1$. This
implies that the Lie algebra ${\Cal G}(G,M)$ is generated
by the coefficients of the vertex operators
$X_{11}(\xi^{i-1},\xi^{-1},z)$ for $0\le
i\le N-1$.
 From Theorem 3.3 we have
 $$
 [X_{11}(\xi^{i-1},\xi^{-1},z_1), X_{11}(\xi^{j-1},\xi^{-1},z_2)]
 $$
 $$
 =\cases X_{11}(\xi^{i-1},\xi^{-j-1}, z_1)\delta({\frac
{\xi^{j}z_2}{z_1}})-
 X_{11}(\xi^{j-1},\xi^{-i-1}, z_2)\delta({\frac
{\xi^{i}z_1}{z_2}})
+
 (D\delta)({\frac {\xi^{j}z_2}{z_1}})c\\
 \;\;\;\;\;\;\;\;\;\;\;\;\;\;\;\;\;\;\;\;\;\;\;\;\;\;\;\;\text{if}\;
i+j=0(\text{mod}N)\\
 X_{11}(\xi^{i-1},\xi^{-j-1}, z_1)\delta({\frac
{\xi^{j}z_2}{z_1}})-
 X_{11}(\xi^{j-1},\xi^{-i-1}, z_2)\delta({\frac
{\xi^{i}z_1}{z_2}})+\\
 \;\;\;\;\;\;\;\;\;\;\;\;\;\;\;\;\;\;\;\;\;\;\;
{\frac {e^{-{\frac {i+j}2}Ln\xi}}{\xi^{i+j}-1}}
 (\delta({\frac {\xi^{j}z_2}{z_1}})-\delta({\frac {\xi^{i}z_1}{z_2}}))
 \;\;\;\;\;\text{if}\; i+j\not=
 0(\text{mod}N).\endcases
 $$
 Recalling the definition of $X_{ij}(a,b,z)$, we have
 $$
 X_{11}(\xi^{i-1},\xi^{-j-1},
z_1)\delta({\frac{\xi^{j}z_2}{z_1}})=
 X_{11}(\xi^{i-1},\xi^{-j-1},
\xi^{j}z_2)\delta({\frac{\xi^{j}z_2}{z_1}})
$$
$$
 =X_{11}(\xi^{i+j-1},\xi^{-1}, z_2)\delta({\frac
{\xi^{j}z_2}{z_1}}),
 $$
while
 $$
 X_{11}(\xi^{j-1},\xi^{-i-1}, z_2)\delta({\frac
{\xi^{i}z_1}{z_2}})
=
 X_{11}(\xi^{i-1},\xi^{-j-1}, \xi^{j}z_1)\delta({\frac
{\xi^{i}z_1}{z_2}})
 $$
$$
=X_{11}(\xi^{i+j-1},\xi^{-1}, z_1)\delta({\frac
{\xi^{i}z_1}{z_2}})
 $$
 for $0\le i,j\le N-1$. Therefore we get
 $$
 [X_{11}(\xi^{i-1},\xi^{-1},z_1), X_{11}(\xi^{j-1},\xi^{-1},z_2)] \tag
4.3
 $$
 $$
 =\cases X_{11}(\xi^{i+j-1},\xi^{-1}, z_2)\delta({\frac
{\xi^{j}z_2}{z_1}})-
 X_{11}(\xi^{i+j-1},\xi^{-1}, z_1)\delta({\frac
{\xi^{i}z_1}{z_2}})+
 (D\delta)({\frac {\xi^{j}z_2}{z_1}})c\\
  \;\;\;\;\;\;\;\;\;\;\;\;\;\;\;\;\;\;\;\;\;\;\;\;\;\;\;\;\text{if}\;
i+j=0(\text{mod}N)\\
 X_{11}(\xi^{i+j-1},\xi^{-1}, z_2)\delta({\frac
{\xi^{j}z_2}{z_1}})-
 X_{11}(\xi^{i+j-1},\xi^{-1}, z_1)\delta({\frac
{\xi^{i}z_1}{z_2}})+\\
 \;\;\;\;\;\;\;\;\;\;\;\;\;\;\;\;\;\;\;\;\;\;\;\;
{\frac {e^{-{\frac
 {i+j}2}Ln\xi}}{\xi^{i+j}-1}}
 (\delta({\frac {\xi^{j}z_2}{z_1}})-\delta({\frac {\xi^{i}z_1}{z_2}})
  \;\;\;\;\text{if}\; i+j\not=
 0(\text{mod}N)c.\endcases
 $$
 Comparing this with the identity (1.18), we obtain
 the following result which was originally due to [F1] and [KKLW].
 
\medskip
 
 \proclaim{Corollary 4.4}  Let $M=1$, and let $G$ be the group generated by
$\xi$, where $\xi$
is a $N$-th
 primitive root of unity for $N\ge 2$. Then the Lie algebra ${\Cal
G}(G,M)$ gives a representation
 of the affine algebra ${\wh {gl_N}}({\Bbb C})$ on the Fock space $V_1$
in the principal picture. The
 representation is given by the mapping
 $$\align &
 F^iE^k\otimes t^k_0\mapsto
 \cases x_{11}(k, \xi^{i-1},\xi^{-1})+
{\frac {e^{{\frac i
2}Ln\xi}}{\xi^i-1}}\delta_{k0}c &\text{if}\; 1\le i\le
 N-1\\x_{11}(k, \xi^{-1},\xi^{-1}) &\text{if}\; i=0\endcases, \\
& c_0\mapsto
{\frac cN}
  \endalign $$
  for $k\in {\Bbb Z}$.
  \endproclaim
 
  \medskip
 
Next we choose $M\ge 2$ and $G=<q>$ where $q\not= 0$ is not a root of
unity. Note that
$X_{ij}(a,b,z)=X_{ij}(1,a^{-1}b,az)$
for $a,b\in G$. We see that the Lie
algebra ${\Cal G}(G,M)$ is generated by the coefficients of the vertex
operators of the form $X_{ij}(1,
q^r,z)$ for all $r\in {\Bbb Z}$ and $1\le i,j\le M$. We apply Theorem 3.3
to obtain
 
\medskip
 
(i) if $r+s=0$, then
$$
[X_{ij}(1,q^r,z_1),X_{kl}(1,q^s,z_2)] \tag 4.5
$$
$$
= X_{il}(1,q^{r+s},z_1)\delta_{jk}\delta({\frac
{z_2}{q^rz_1}})-
X_{kj}(1,q^{r+s},z_2)\delta_{il}\delta({\frac {z_1}{q^sz_2}})
+\delta_{il}\delta_{jk}(D\delta)({\frac {z_2}{q^rz_1}})c,
$$
 
(ii) if $r+s\not= 0$, then
$$\align &
[X_{ij}(1,q^r,z_1),X_{kl}(1,q^s,z_2)]\tag 4.6\\
=& X_{il}(1,q^{r+s},z_1)\delta_{jk}\delta({\frac
{z_2}{q^rz_1}})-
X_{kj}(1,q^{r+s},z_2)\delta_{il}\delta({\frac {z_1}{q^sz_2}})\\
&+{\frac {q^{\frac
{r+s}2}}{1-q^{r+s}}}\delta_{il}\delta_{jk}(\delta({\frac {z_2}{q^rz_1}})-
\delta({\frac {z_1}{q^sz_2}}))c.
\endalign $$
Comparing the above two identities with the identity (1.19), we derive the
following result which
was given in [G1].
 
\medskip
 
\proclaim{Corollary 4.7}  Let $M\ge 2$, and $G=<q>$ be the group
generated by  $q\not= 0$
and
$q$ is not a root
of unity. Then the Lie algebra  ${\Cal G}(G,M)$ of operators, acting on
$V_M$, gives a representation of
the Lie algebra ${\wh {gl_M}}({\Bbb C}_q)$.
 The representation is given by the mapping
$$
E_{ij}\otimes t^m_0t^r_1\mapsto \cases x_{ij}(m, 1, q^r)+{\frac {q^{\frac
{r}2}}{1-q^{r}}}\delta_{ij}\delta_{m0}c
 &\text{if}\; r\not= 0\\
 x_{ij}(m,1,1) &\text{if}\; r=0\endcases
 $$
 $$
 c_0\mapsto c,\;\;\;\;\;\;c_1\mapsto 0\;\;\;\;\;\;\;\;\;\;\;\;\;\;\;
 $$
 for $m,r\in {\Bbb Z}$.\
\endproclaim
 
\medskip
 
\noindent{\bf Remark 4.8}. The representation of the Lie algebra ${\wh
{gl_M}}({\Bbb C}_Q)$
given in (4.7) is called
 the homogeneous realization. This is because of the fact that we are using
the homogeneous gradation. Moreover the
 algebra ${\Cal G}(G,M)$ contains a subalgebra of ${\Cal G}(<1>,M)$ which
is
generated by the operators $
 x_{ij}(m,1,1)$ for $1\le i,j\le M$ and $m\in {\Bbb Z}$, and it is clear
that this subalgebra is nothing but the
 affine algebra ${\wh {gl_M}}({\Bbb C})$ in the homogeneous picture.
 
\medskip
 
 Similarly,  we may have the  principal realization of the Lie algebra
${\wh {gl_N}}({\Bbb
C}_Q)$. For this purpose, we
 choose the group $G=<\xi,q>$ where $q\not= 0$ is not a root of unity and
$\xi$ is a $N$-th primitive
 root of unity. Let $M=1$. Then the Lie algebra ${\Cal G}(G,M)$ is
generated by the coefficients of the
 vertex operators of the form $X_{11}(\xi^{i-1},\xi^{-1}q^r,z)$ for all
$r\in {\Bbb Z}$ and $0\le i\le N-1
 $. From Theorem 3.3 we have
 
\medskip
 
 (i) if $r+s=0$ and ${\overline {i+j}}=0$(mod$N$), then
 $$\align &
 [X_{11}(\xi^{i-1},\xi^{-1}q^r,z_1), X_{11}(\xi^{j-1},\xi^{-1}q^s,z_2)]\\
=&
 X_{11}(\xi^{i+j-1},\xi^{-1}q^{r+s},\xi^{-j}z_1)\delta({\frac
{\xi^jz_2}{q^rz_1}}) \\
 & -
 X_{11}(\xi^{i+j-1},\xi^{-1}q^{r+s},\xi^{-i}z_2)\delta({\frac
{\xi^iz_1}{q^sz_2}})
 +(D\delta)({\frac {\xi^jz_2}{q^rz_1}})c,
\endalign  $$
 
 (ii) if $r+s\not= 0$ or ${\overline {i+j}}\not= 0$(mod$N$), then
 $$\align &
 [X_{11}(\xi^{i-1},\xi^{-1}q^r,z_1), X_{11}(\xi^{j-1},\xi^{-1}q^s,z_2)]\\
=&
 X_{11}(\xi^{i+j-1},\xi^{-1}q^{r+s},\xi^{-j}z_1)\delta({\frac
{\xi^jz_2}{q^rz_1}})
 -
 X_{11}(\xi^{i+j-1},\xi^{-1}q^{r+s},\xi^{-i}z_2)\delta({\frac
{\xi^iz_1}{q^sz_2}})
 \\
& +{\frac {e^{{\frac {i+j} 2}Ln\xi}q^{\frac {r+s}2}}{\xi^{i+j}-q^{r+s}}}
 (\delta({\frac {\xi^jz_2}{q^rz_1}})-\delta({\frac {\xi^iz_1}{q^sz_2}}))c.
 \endalign $$
 Thus if we write
 $$
 {\bar X}_{11}(\xi^{i-1},\xi^{-1}q^r,z)=\cases
X_{11}(\xi^{-1},\xi^{-1},z) &\text{if}\; r=0, i=
 0(\text{mod}N)\\
 X_{11}(\xi^{i-1},\xi^{-1}q^r,z)+
 {\frac {e^{{\frac i 2}Ln\xi}q^{\frac {r}2}}{\xi^{i}-q^{r}}}c 
&\text{otherwise}\endcases
 $$
 for $r\in {\Bbb Z}$ and $0\le i\le N-1$, then
 the above two identities can be written as one identity
 $$
 [{\bar X}_{11}(\xi^{i-1},\xi^{-1}q^r,z_1),
{\bar X}_{11}(\xi^{j-1},\xi^{-1}q^s,z_2)] \tag 4.9
$$
$$
=
 {\bar X}_{11}(\xi^{i+j-1},\xi^{-1}q^{r+s},\xi^{-j}z_1)\delta({\frac
{\xi^jz_2}{q^rz_1}})
  -
 {\bar X}_{11}(\xi^{i+j-1},\xi^{-1}q^{r+s},\xi^{-i}z_2)\delta({\frac
{\xi^iz_1}{q^sz_2}}) 
$$
$$
 +\delta_{r+s,0}\delta_{{\overline {i+j}},0}(D\delta)({\frac
{\xi^jz_2}{q^rz_1}})c.
$$
 Comparing this with the identity (1.20), we get the following result which
was given in [BS] for the
$N=2$ case and in [G2] for arbitrary $N$.
 
 \medskip
 
 \proclaim{Corollary 4.10} Let $M=1$ and $G=<\xi,q>$ be an admissible
subgroup
 of ${\Bbb C}^{\times}$
 generated
by $q$
with
 $q\not= 0$  not a root of unity, and $\xi$ is a $N$-th primitive
 root of unity for $N\ge 2$. Then the Lie algebra ${\Cal G}(G,M)$ of
operators, acting on $V_1$, gives a representation of
 the algebra ${\wh {gl_N}}({\Bbb C}_{q^N})$.
The representation is given by
the mapping
$$
F^iE^{m}\otimes t^m_0t^r_1\mapsto \cases x_{11}(m,\xi^{i-1},\xi^{-1}q^r
)+
\delta_{m,0}{\frac {e^{{\frac i 2}Ln\xi}q^{\frac {r}2}}{\xi^i-q^{r}}}c
 &\text{if}\; i\not= 0(\text{mod}N)\;\;\text{or}\;\;r\not= 0\\
 x_{11}(m,\xi^{-1},\xi^{-1}) &\text{if}\;
i=0(\text{mod}N)\;\;\text{and}\;\;r=0\endcases
 $$
 $$
 c_0\mapsto {\frac cN},\;\;\;\;\;\;c_1\mapsto
0\;\;\;\;\;\;\;\;\;\;\;\;\;\;\;
 $$
 for $r\in {\Bbb Z}$ and $0\le i\le N-1$.\endproclaim
 
 \medskip
 
 \noindent{\bf Remark 4.11}. In general, let $M,N\ge 2$ be integers,  $
G=<\xi,q_1,\cdots, q_{\nu}> 
$
an admissible subgroup of
${\Bbb
C}^{\times}$ with finitely
many generators, where $q_1,\cdots, q_{\nu}$ are the free
generators of $G$ and $\xi$ is
 an $N$-th root of unity.  Then the Lie algebra ${\Cal G}(G,M)$  of
operators, acting on the
 Fock space $V_M$, gives a representation to the Lie algebra ${\wh
{gl_{MN}}}({\Bbb C}_Q)$ where the quantum
 torus ${\Bbb C}_Q={\Bbb C}_Q[t_0^{\pm 1},t_1^{\pm 1},\cdots,t_{\nu}^{\pm
1}]$ is determined by the matrix
 $Q=(q_{ij})_{(\nu+1)\times(\nu+1)}$ with $q_{i0}=q_i^N, q_{0i}=q^{-N}_i$
for $1\le i\le\nu$, and $q_{ij}=1$
 for all other values of $i,j$.
 
 In particular,
if $\nu=1$, that is $G=<\xi,q>$, then the Lie algebra ${\Cal G}(G,M)$
  is generated by the coefficient operators of the vertex operators $
X_{ij}(\xi^{k-1},\xi^{-1}q^r,z)$ for $1\le i,j\le M$, $1\le k\le N-1$ and 
$r\in \Bbb Z$. Moreover, from Theorem 3.3, we have

(i) if ${\overline {k+k'}}\not= 0$ or $r+s\not= 0$, then
$$\align &
[X_{ij}(\xi^{k-1},\xi^{-1}q^r,z_1), X_{i'j'}(\xi^{k'-1},\xi^{-1}
q^{r'},z_2)]\\
=& \delta_{ji'}X_{ij'}(\xi^{k-1},\xi^{-1-k'}q^{r+r'},z_1)\delta({\frac{\xi^{k'}z_2} {q^rz_1}})\\
&-\delta_{j'i}X_{i'j}(\xi^{k'-1},\xi^{-1-k}q^{r+r'},z_2))\delta({\frac
{\xi^{k}z_1} {q^{r'}z_2}})\\
&+{\frac {e^{{\frac {k+k'}{2}}Ln\xi}q^{\frac {r+r'} 2}}{\xi^{k+k'}-q^{r+r'}}}\delta_{ji'}
\delta_{j'i}(\delta({\frac{\xi^{k'}z_2} {q^rz_1}})-\delta({\frac{\xi^{k}z_1} {q^{r'}z_2}}))c\\
=&\delta_{ji'}X_{ij'}(\xi^{k+k'-1},\xi^{-1}q^{r+r'},\xi^{-k'}z_1)\delta({\frac{\xi^{k'}z_2} {q^rz_1}})\\
&-\delta_{j'i}X_{i'j}(\xi^{k+k'-1},\xi^{-1}q^{r+r'},\xi^{-k}z_2))\delta({\frac{\xi^{k}z_1} {q^{r'}z_2}})\\
&+{\frac {e^{{\frac {k+k'}{2}}Ln\xi}q^{\frac {r+r'} 2}}{\xi^{k+k'}-q^{r+r'}}}\delta_{ji'}\delta_{j'i}
(\delta({\frac{\xi^{k'}z_2} {q^rz_1}})-\delta({\frac{\xi^{k}z_1} {q^{r'}z_2}}))c,
\endalign
$$
and

(ii) if ${\overline {k+k'}}= 0$ and $r+s= 0$, then
$$\align &
[X_{ij}(\xi^{k-1},\xi^{-1}q^r,z_1), X_{i'j'}(\xi^{k'-1},\xi^{-1}q^{r'},z_2)]\\
=& \delta_{ji'}X_{ij'}(\xi^{k-1},\xi^{-1-k'}q^{r+r'},z_1)\delta({\frac  
{\xi^{k'}z_2} {q^rz_1}})\\
&-
\delta_{j'i}X_{i'j}(\xi^{k'-1},\xi^{-1-k}q^{r+r'},z_2))\delta({\frac
{\xi^{k}z_1} {q^{r'}z_2}})
+
\delta_{ji'}\delta_{j'i}(D\delta)({\frac
{\xi^{k'}z_2} {q^rz_1}})c\\
=&\delta_{ji'}X_{ij'}(\xi^{-1},\xi^{-1},\xi^{-k'}z_1)\delta({\frac
{\xi^{k'}z_2} {q^rz_1}})\\
&-
\delta_{j'i}X_{i'j}(\xi^{-1},\xi^{-1},\xi^{-k}z_2))\delta({\frac
{\xi^{k}z_1} {q^{r'}z_2}})
+\delta_{ji'}\delta_{j'i}(D\delta)({\frac
{\xi^{k'}z_2} {q^rz_1}})c.
\endalign $$
Set
$$
{\overline X}_{ij}(\xi^{k-1},\xi^{-1}q^r,z)
$$
$$= 
 \cases X_{ij}(\xi^{k-1},\xi^{-1}q^r,z )+
\delta_{ij}{\frac {e^{{\frac k 2}Ln\xi}q^{\frac
{r}2}}{\xi^k-q^{r}}}c
 &\text{if}\; k\not= 0\text{ (mod}N)\;\;\text{or}\;\;r\not= 0\\
 X_{ij}(\xi^{-1}, \xi^{-1},z) &\text{if}\; k=0\text{
(mod}N)\;\;\text{and}\;\;r=0\endcases
 $$
Then we have
$$
[{\overline X}_{ij}(\xi^{k-1},\xi^{-1}q^r,z_1),
{\overline X}_{i'j'}(\xi^{k'-1},\xi^{-1}q^{r'},z_2)]
$$
$$
=\delta_{ji'}{\overline
X}_{ij'}(\xi^{k+k'-1},\xi^{-1}q^{r+r'},\xi^{-k'}z_1)\delta({\frac
{\xi^{k'}z_2} {q^rz_1}})
$$
$$-
\delta_{j'i}{\overline
X}_{i'j}(\xi^{k+k'-1},\xi^{-1}q^{r+r'},\xi^{-k}z_2))\delta({\frac
{\xi^{k}z_1} {q^{r'}z_2}})
+\delta_{ji'}\delta_{j'i}\delta_{\overline
{k+k'},0}\delta_{r+r',0}(D\delta)({\frac
{\xi^{k'}z_2} {q^rz_1}})c.
$$
Comparing this identity with  identity (1.21), we get

\proclaim
{\bf Corollary 4.12}  Let $M,N\ge 2$, and let $G=<\xi,q>$ be an admissible
subgroup
 of ${\Bbb C}^{\times}$
 generated
by $q$
with
 $q\not= 0$  not a root of unity, and $\xi$ an $N$-th primitive
 root of unity. Then the Lie algebra ${\Cal G}(G,M)$ of
operators, acting on $V_1$, gives a representation of
 the algebra ${\wh {gl_{MN}}}({\Bbb C}_{q^N})$, and the representation is
given
by
the mapping
$$
E_{ij}\otimes F^kE^{m}\otimes t^m_0t^r_1
$$
$$
\mapsto \cases
x_{ij}(m,\xi^{k-1},\xi^{-1}q^r
)+
\delta_{m,0}\delta_{ij}
{\frac {e^{{\frac k 2}Ln\xi}q^{\frac
{r}2}}{\xi^i-q^{r}}}c
 &\text{if}\; k\not= 0(\text{mod}N)c\;\;\text{or}\;\;r\not= 0\\
 x_{ij}(m,\xi^{-1},\xi^{-1}) &\text{if}\;
k=0(\text{mod}N)\;\;\text{and}\;\;r=0\endcases
 $$
 $$
 c_0\mapsto {\frac cN},\;\;\;\;\;\;c_1\mapsto
0\;\;\;\;\;\;\;\;\;\;\;\;\;\;\;
 $$
 for $r\in {\Bbb Z}$ and $0\le k\le N-1$.\endproclaim
 
 \medskip

The Lie algebra ${\Cal G}(G,M)$ given in the previous corollary contains 
two interesting 
  subalgebras which give representations to the Lie algebras ${\wh
{gl_{N}}}({\Bbb C}_{q^N})$ and
 ${\wh {gl_{M}}}({\Bbb C}_q)$ . Moreover we will see that
these two subalgebras 
 contain subalgebras that give representations to the affine algebras
${\wh {gl_{N}}}({\Bbb C})$ of level $M$ and
 ${\wh {gl_{M}}}({\Bbb C})$ of level $N$ respectively.
 Indeed, for $a,b\in G=<\xi,q>$, let
 $$
 Y(a,b,z)=\sum_{k=1}^MX_{kk}(a,b,z) \tag 4.13
 $$
 and formally write
$$Y(a,b,z)=\sum_{m\in {\Bbb Z}}y(m,a,b)z^{-m}.\tag 4.14$$
 
 Let ${\Cal L}_1$ be the Lie algebra
 generated by all of the coefficients of $Y(a,b,z)$ for $a,b\in G$.
 We note that $$Y(\xi^iq^r,\xi^jq^s,z)=
 Y(\xi^{i-j-1},\xi^{-1}q^{s-r},\xi^{j+1}q^{r}z),\tag 4.15$$
 so ${\Cal L}_1$ is
indeed generated by the coefficients
 of the vertex operators with the form $Y(\xi^{i-1},\xi^{-1}q^{r},z)$ for
$r\in {\Bbb Z}$ and $0\le i\le
 N-1
 $. Moreover, applying Theorem 3.3, we have, if $r+s\not= 0$ or
${\overline {i+j}}\not= 0$(mod$N$)
 
 $$ 
 [Y(\xi^{i-1},\xi^{-1}q^{r},z_1),Y(\xi^{j-1},\xi^{-1}q^{s},z_2)]\tag
4.16
$$
$$
=
 \sum_{k=1}^M[X_{kk}(\xi^{i-1},\xi^{-1}q^{r},z_1),X_{kk}(\xi^{j-1},
 \xi^{-1}q^{s}
,z_2)] 
 $$
 $$
 =\sum_{k=1}^M\left(X_{kk}(\xi^{i+j-1},\xi^{-1}q^{r+s},\xi^{-j}z_1)
 \delta({\frac 
{\xi^jz_2}{q^rz_1}})-
 X_{kk}(\xi^{i+j-1},\xi^{-1}q^{r+s},\xi^{-j}z_2)\delta({\frac
{\xi^iz_1}{q^sz_2}})\right.
 $$
 $$
 \left.+{\frac {e^{{\frac {i+j} 2}Ln\xi}q^{\frac
{r+s}2}}{\xi^{i+j}-q^{r+s}}}(\delta({\frac {\xi^jz_2}{q^rz_1}})-
 \delta({\frac {\xi^iz_1}{q^sz_2}}))\right)
 $$
 $$
 =Y(\xi^{i+j-1},\xi^{-1}q^{r+s},\xi^{-j}z_1)\delta({\frac
{\xi^jz_2}{q^rz_1}})-
 Y(\xi^{i+j-1},\xi^{-1}q^{r+s},\xi^{-j}z_2)\delta({\frac
{\xi^iz_1}{q^sz_2}})
 $$
 $$
 +{\frac {e^{{\frac {i_j}2}Ln\xi}q^{\frac
{r+s}2}}{\xi^{i+j}-q^{r+s}}}(\delta({\frac {\xi^jz_2}{q^rz_1}})-
 \delta({\frac {\xi^iz_1}{q^sz_2}})),
 $$
 while if $r+s= 0$ and ${\overline {i+j}}= 0$(mod$N$), then
 $$
 [Y(\xi^{i-1},\xi^{-1}q^{r},z_1),Y(\xi^{j-1},\xi^{-1}q^{s},z_2)]
\tag 4.17
$$
$$
=
 \sum_{k=1}^M[X_{kk}(\xi^{i-1},\xi^{-1}q^{r},z_1),X_{kk}(\xi^{j-1},\xi^{-1}
 q^s
,z_2)]
$$
 $$
 =\sum_{k=1}^M\left(X_{kk}(\xi^{-1},\xi^{-1},\xi^{-j}z_1)\delta({\frac
{\xi^jz_2}{q^rz_1}})-
 X_{kk}(\xi^{-1},\xi^{-1},\xi^{-j}z_2)\delta({\frac {\xi^iz_1}{q^sz_2}})
 +(D\delta)({\frac {\xi^jz_2}{q^rz_1}})\right)
 $$
 $$
 =Y(\xi^{-1},\xi^{-1},\xi^{-j}z_1)\delta({\frac {\xi^jz_2}{q^rz_1}})-
 Y(\xi^{-1},\xi^{-1},\xi^{-j}z_2)\delta({\frac {\xi^iz_1}{q^sz_2}})
 +M(D\delta)({\frac {\xi^jz_2}{q^rz_1}}).
 $$
 Therefore, if we define
 $$
 {\bar Y}(\xi^{i-1},\xi^{-1}q^{r},z)=\cases
Y(\xi^{i-1},\xi^{-1}q^{r},z)+M
 {\frac {e^{{\frac i2}Ln\xi}q^{\frac {r}2}}{\xi^{i}-q^{r}}}c
&\text{if}\;r\not= 0\;\text{or}\;{\bar i}\not=
 0\\
 Y(\xi^{-1},\xi^{-1},z) &\text{if}\;r= 0\;\text{and}\;{\bar i}=
 0,\endcases
 $$
 then we can rewrite the  two identities (4.16) and (4.17) into just
one identity
 $$
 [{\bar Y}(\xi^{i-1},\xi^{-1}q^{r},z_1),{\bar
Y}(\xi^{j-1},\xi^{-1}q^{s},z_2)]\tag 4.18
$$
$$
 ={\bar Y}(\xi^{-1},\xi^{-1},\xi^{-j}z_1)\delta({\frac
{\xi^jz_2}{q^rz_1}})
 {\bar Y}(\xi^{-1},\xi^{-1},\xi^{-j}z_2)\delta({\frac
{\xi^iz_1}{q^sz_2}})
 +M(D\delta)({\frac {\xi^jz_2}{q^rz_1}})c.
 $$
 Therefore we have the following result
 
 \medskip
 
 \proclaim{Proposition 4.19}  The Lie algebra ${\Cal L}_1$ of operators
acting on $V_M$ gives a
representation of the
 Lie algebra ${\wh {gl_{N}}}({\Bbb C}_{q^N})$, 
 and the
 representation
 is given by the mapping
 $$
 F^iE^{m}\otimes t^m_0t^r_1\mapsto \cases y(m,\xi^{i-1},\xi^{-1}q^r )+
M{\frac {e^{{\frac i2}Ln\xi}q^{\frac {r}2}}{\xi^i-q^{r}}}\delta_{m,0}c
 &\text{if}\; i\not= 0\text{ (mod}N)\;\;\text{or}\;\;r\not= 0\\
 y(m,1,1) &\text{if}\; i=0\text{ (mod}N)\;\;\text{and}\;\;r=0\endcases
 $$
 $$
 c_0\mapsto Mc,\;\;\;\;\;\;c_1\mapsto 0\;\;\;\;\;\;\;\;\;\;\;\;\;\;\;
 $$
 for $r\in {\Bbb Z}$ and $0\le i\le N-1$.\endproclaim
 
 \medskip
 
 Recall from (4.7), the Fock space $V_M$ affords a representation of the
Lie algebra ${\Cal G}(<q>,M)\subset
 $ ${\Cal G}(<\xi,q>,M)$, where $\xi, q$ are given in (4.10), and
$${\Cal G}(<q>,M)=\text{ span }\{c \text{ and } x_{ij}(m,1,q^r)|\text{
for
} m,r\in {\Bbb Z}, 1\le i,j
 \le M\}. $$
 Now we define a subalgebra of ${\Cal G}(<q>,M)\subset$ ${\Cal
G}(<\xi,q>,M)$
 $$
 {\Cal L}_2=\text{span}\{c\;\; \text{and}\;\; x_{ij}(Nm,1,q^r)|\;
\text{for} \;\;
 m,r\in {\Bbb Z}, 1\le i,j
 \le M\}.
 $$
 Then we have
 
 \medskip
 
 \proclaim{Proposition 4.20} $ {\Cal L}_2$ forms a Lie subalgebra of
${\Cal G}(<q>,M)$, and $ {\Cal
L}_2$
 is also isomorphic to
 ${\Cal G}(<q>,M)$ via the isomorphism given by
 $$
 x_{ij}(m,1,q^r)\mapsto x_{ij}(Nm,1,q^r),\;\;\;\;\;\;c\mapsto Nc.
 $$
 Therefore $ {\Cal L}_2$ gives
 a  representation of ${\wh{gl_{N}}}({\Bbb C}_{q^N})$.
 \endproclaim
 
 \medskip
 
 \proclaim{Proposition 4.21}  For $m,n,r,s,\in {\Bbb Z}$ and $i\not= 0$
(mod $N$), $1\le k\not= l\le
M$, we
have
 $$\align &
 [y(m,\xi^{i-1},\xi^{-1}q^r), x_{kl}(Nn,1,q^s)]
\tag 4.22\\
=& (q^{rNn}-q^{sm})x_{kl}(m+Nn,\xi^{i-1},\xi^{-1}q^{r+s}).
 \endalign
$$
\endproclaim
 \demo{Proof} We apply Theorem 3.3 to obtain
 $$\align &
 [Y(\xi^{i-1},\xi^{-1}q^r, z_1), X_{kl}(1,q^s,z_2)]=
 [\sum_{j=1}^MX_{jj}(\xi^{i-1},\xi^{-1}q^r, z_1), X_{kl}(1,q^s,z_2)]\\
=&\sum_{j=1}^M\{X_{jl}(\xi^{i-1},\xi^{-1}q^{r+s},z_1)\delta_{jk}\delta({\frac
{z_2}{\xi^{-1}q^rz_1}})
-X_{kj}(1,\xi^{-i}q^{r+s},z_2)\delta_{jl}\delta({\frac
{\xi^{i-1}z_1}{q^sz_2}})\}\\
=& X_{kl}(\xi^{i-1},\xi^{-1}q^{r+s},z_1)\delta({\frac
{z_2}{\xi^{-1}q^rz_1}})
-X_{kl}(1,\xi^{-i}q^{r+s},z_2)\delta({\frac {\xi^{i-1}z_1}{q^sz_2}})\\
=& X_{kl}(\xi^{i-1},\xi^{-1}q^{r+s},z_1)\delta({\frac
{z_2}{\xi^{-1}q^rz_1}})
- X_{kl}(\xi^{i-1},\xi^{-1}q^{r+s},q^{-s}z_1) \delta({\frac
{\xi^{i-1}z_1}{q^sz_2}}).
\endalign
$$
This then gives
$$
[y(m,\xi^{i-1},\xi^{-1}q^r),
x_{kl}(n,1,q^s)]=\xi^{-n}(q^{rn}-q^{sm}\xi^{in})
x_{kl}(m+n,\xi^{i-1},\xi^{-1}q^{r+s}).
$$ which immediately implies (4.22).\qed
 
\enddemo
 
 \medskip
 
\noindent{\bf Remark 4.23} Let ${\Cal G}_i\subset {\Cal L}_i\subset {\Cal
G}(<\xi,q>,M)$, $i=1,2$, be such
that
 $$
 {\Cal
G}_1=\text{span}\{c\;\; \text{and}\;\; y(m,\xi^{i-1},\xi^{-1})|\; \text{
for }
 m\in {\Bbb Z}, 0\le i
 \le N-1\}, 
 $$
 $$
 {\Cal G}_2=\text{span}\{c\;\; \text{and}\;\; x_{ij}(Nm,1,1)|\; \text{
for }
\;\;
 m\in {\Bbb Z}, 1\le i,j
 \le M\}.
 $$
  Then the two subalgebras ${\Cal G}_1, {\Cal G}_2$ of ${\Cal
G}(<\xi,q>,M)$
  respectively give representations of the
 affine algebra ${\wh{gl_N}}({\Bbb C})$ of level $M$ and
${\wh{gl_M}}({\Bbb C})$ of
level $N$. Let
 ${\Cal G}_i'$ be the derived algebras of ${\Cal G}_i$. Then we have the
so
called dual pair property given in [F1]: 
 $
 [{\Cal G}_1', {\Cal G}_2']=(0).    
$
 However, clearly, we have  $[{\Cal L}_1', {\Cal L}_2']\not= (0)$.
                                                              
\bigskip
 
\noindent{\bf Acknowledgements }
 
\medskip
 
SB and YG are  supported by grants from  the
Natural Sciences and Engineering Research Council of Canada.
ST is  supported by a grant from the National
Natural Science Foundation of
China. The authors would like to thank the Fields Institute  for its
hospitality
during the preparation of this work.

\enddocument